\documentclass[a4paper, 10pt]{amsart}
\pdfoutput=1
\usepackage[foot]{amsaddr}
\usepackage{geometry}
\usepackage{lineno}
\usepackage{bm}
\usepackage{amsfonts}
\usepackage{amsmath}
\usepackage{amssymb}
\usepackage{amsthm}
\usepackage[dvipsnames]{xcolor}
\usepackage{graphicx}
\usepackage{caption}
\usepackage{float}
\usepackage{subcaption}
\usepackage{overpic}
\usepackage{url}
\usepackage{comment}

\graphicspath{{images/}}
\def\u{{\bm u}}

\begin{document}
\title[]{Non-intrusive PODI-ROM for patient-specific aortic blood flow in presence of a LVAD device}

\author{M. Girfoglio\textsuperscript{1,*}}
\email{mgirfogl@sissa.it}
\author{F. Ballarin\textsuperscript{1}}
\email{fballarin@sissa.it}
\author{G. Infantino\textsuperscript{2}}
\email{giuseppe.infantino@studenti.polito.it}
\author{F. Nicolò\textsuperscript{3}}
\email{nicolo\_francy84@hotmail.it}
\author{A. Montalto\textsuperscript{3}}
\email{andrea.montalto@libero.it}
\author{G. Rozza\textsuperscript{1}}
\email{grozza@sissa.it}
\author{R. Scrofani\textsuperscript{4}}
\email{roberto.scrofani@asst-fbf-sacco.it}
\author{M. Comisso\textsuperscript{3}}
\email{marina.comisso@gmail.com}
\author{F. Musumeci\textsuperscript{3}}
\email{fr.musumeci@gmail.com}
\thanks{\textsuperscript{*}Corresponding Author.}
\address{\textsuperscript{1}SISSA, Scuola Internazionale Superiore di Studi Avanzati, Area di Matematica, mathLab Trieste, Italy.}
\address{\textsuperscript{2}Politecnico di Torino, Collegio di Ingegneria Matematica, Modelli Matematici e Simulazioni Numeriche, Torino, Italy}
\address{\textsuperscript{3}Azienda Ospedaliera San Camillo, Unità Operativa Complessa di Cardiochirurgia e Chirurgia dei Trapianti, Roma, Italy}
\address{\textsuperscript{4}Azienda Ospedaliera FBF Luigi Sacco, Dipartimento di Cardiochirurgia, Milano, Italy}


\date{}

\dedicatory{}

\begin{abstract}
Left ventricular assist devices (LVADs) are used to provide haemodynamic support to patients with critical cardiac failure. Severe complications can occur because of the modifications of the blood flow in the aortic region. In this work, the effect of a continuous flow LVAD device on the aortic flow is investigated by means of a non-intrusive reduced order model (ROM) built using the proper orthogonal decomposition with interpolation (PODI) method. The full order model (FOM) is represented by the incompressible Navier-Stokes equations discretized by using a Finite Volume (FV) technique, coupled with three-element Windkessel models to enforce outlet boundary conditions in a multi-scale approach. A patient-specific framework is proposed: a personalized geometry reconstructed from Computed Tomography (CT) images is used and the individualization of the coefficients of the three-element Windkessel models is based on experimental data provided by the Right Heart Catheterization (RCH) and Echocardiography (ECHO) tests. Pre-surgery configuration is also considered at FOM level in order to further validate the model. A parametric study with respect to the LVAD flow rate is considered. The accuracy of the reduced order model is assessed against results obtained with the full order model.
\end{abstract}

\maketitle
\section{Introduction}\label{sec:intro}
Left ventricular assist devices (LVADs) provide full or partial mechanical circulatory support to the left ventricle of the heart. LVADs are increasingly used as both bridge to transplantation (BTT) and destination therapy (DT), for the treatment of patients with advanced heart failure (HF) refractory to maximal medical therapy \cite{A, B, C}. The use of LVADs has been associated with an increased risk of thrombus formation in the aortic region because of the formation of stagnation points and recirculation zones; indeed, while first generation devices provided pulsatile flows, current LVADs produce continuous flow (cf-LVADs). This constant flow to the aortic root may lead to decreased excursion or even complete closure of the aortic valve (AV), particularly at high pump speeds. The resultant stasis in the aortic root forms a nidus for clot formation. Aortic root thrombosis has been recognized as a major complication of cf-LVAD therapy which frequently necessitates device exchange in eligible patients to restore forward flow and prevent embolic stroke \cite{D, E, F, G}. 

Several works deal with the computational investigation of the hemodynamics in the aortic region in the presence of a LVAD device, both in a single configuration \cite{Bonnemain2012} and varying of physical (LVAD flow rate \cite{Bazilevs2009, Mazzitelli2016}) and geometrical (cannula angle \cite{Inci2012, Quaini2011, Karmonik2012, Karmonik2014, Osorio2013} and anastomosis position \cite{May-Newman2004, May-Newman2006, Karmonik2014b, Caruso2015, Quaini2011, Karmonik2012, Mazzitelli2016, Brown2011, Aliseda2017, Bonnemain2013, Bazilevs2009, Brown2011, Zbigniew2018}) parameters. 
In all these works, high fidelity full order models (FOMs) are used, based on either finite element and finite volume simulations.

Reduced order models (ROMs) (see, e.g., \cite{quarteroniRB2016, hesthaven2015certified, bennerParSys, ModelOrderReduction}) have been proposed as an efficient tool to approximate full order systems by significantly reducing the computational cost required to obtain numerical solutions in a parametric setting. The basic idea on which ROM is based is that often the parametric dependence of the problem at hand has an intrinsic dimension much lower than the number of degrees of freedom of the discretized system. In oreder to reach this dimensionality reduction, a database of several solutions is first collected by solving the original high fidelity model for different physical and/or geometrical parameters (\emph{offline phase}). Then, all the solutions are combined and compressed in order to build the space onto which we can project the solution manifold and efficiently compute the solutions for the new parameters (\emph{online phase}). 


In this work, a data-driven reduced order model based on the proper orthogonal decomposition with
interpolation (PODI) method \cite{Bui-Thanh2004} is used for the investigation of the modifications to aortic blood flow patterns induced by the presence of the outflow cannula of a LVAD device. PODI is an equation-free method based on proper orthogonal decomposition capable to build a reduced order model without any knowledge about the equations of the original problem. Therefore, this method does not require information about the full order formulation and not even modifications to the numerical solver. Data-driven non-intrusive ROMs have been widely used within industrial applications framework (see, e.g. \cite{Tezzele2018}). 
In our application, the parametrization is applied with respect to the flow rate provided by the LVAD device. The open-source finite volume solver OpenFOAM \cite{Weller1998} is used to generate the FOM solutions which are then used as a training set for the ROM. In order to obtain a model able to reproduce clinical configurations, geometry is reconstructed from patient-specific Computed Tomography (CT) images. Moreover, a multi-scale approach was adopted by coupling three-element Windkessel models \cite{Westerhof2008}, used as boundary conditions for the aorta model, and which coefficients are estimated by using experimental data provided by Right Heart Catheterization (RHC) and Echocardiography (ECHO) tests. Pre-surgery configuration is also considered in order to further validate the FOM and the estimation procedure. Therefore, a complete patient-specific framework is proposed. To the best of our knowledge, parametric ROM for modeling realistic aortic flow in presence of LVAD devices is introduced in this paper for the first time.

The work is organized as follows. In Sec. \ref {sec:mat_and_met} materials and methods are presented and in Sec. \ref{sec:results} the achieved results are introduced and discussed. Finally, in Sec. \ref{sec:conclusion} conclusions and perspectives are provided. 



\section{Clinical data and discretization by full order and reduced order models}\label{sec:mat_and_met}

\subsection{Clinical data}
In this work a patient, a 66 years old man, is considered. CT, RHC and ECHO tests have been carried out both in pre-surgery and post-surgery (i.e., after receiving the LVAD device) configuration. The LVAD implanted is the Heartmate 3$^\text{TM}$ Left Ventricular Assist System \cite{HeartMate}.

\subsection{Geometrical model}
\begin{figure}[]
\centering
 \begin{overpic}[width=0.4\textwidth]{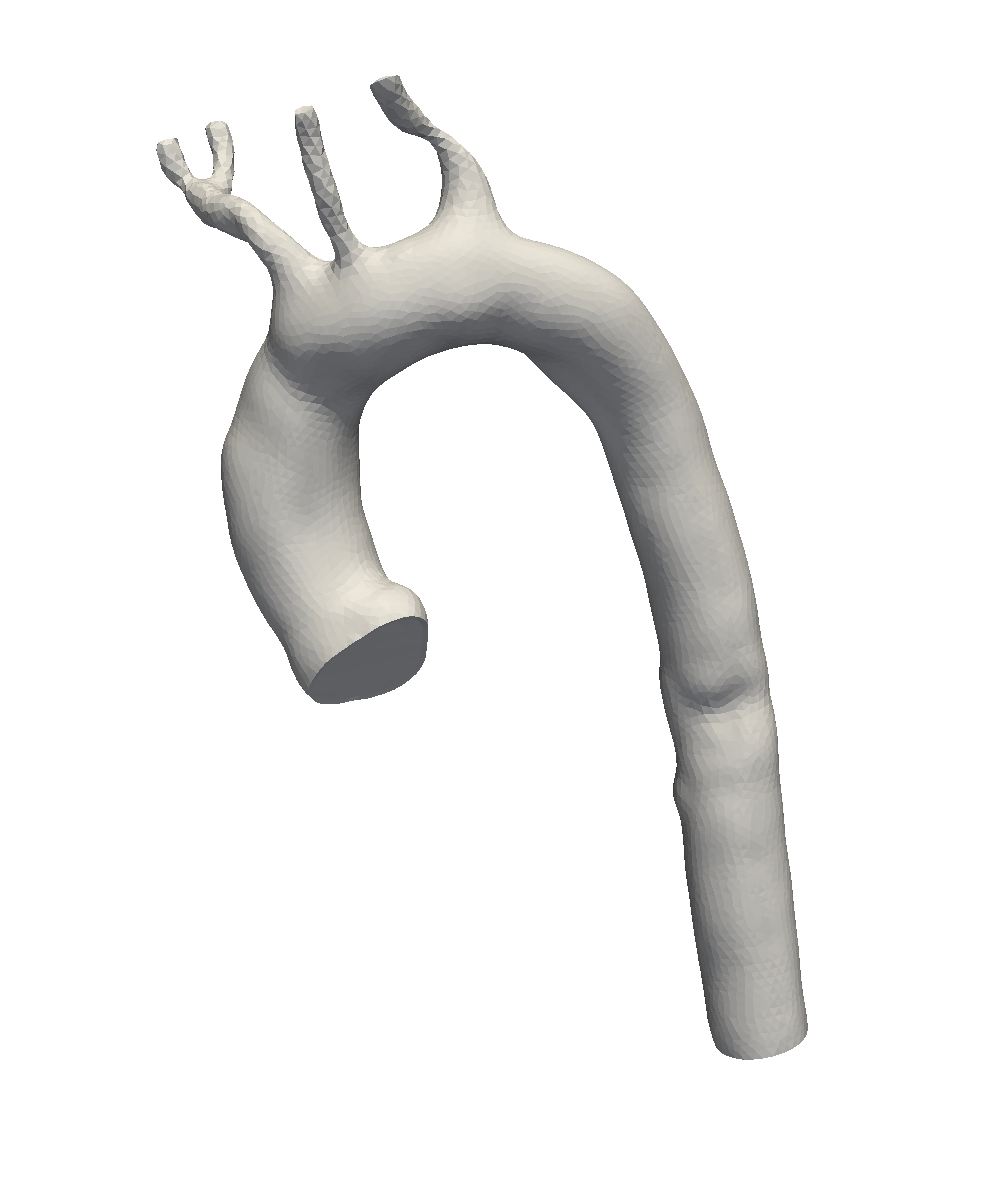}
\put(30,100){\small{a)}}
      \end{overpic}
 \begin{overpic}[width=0.4\textwidth]{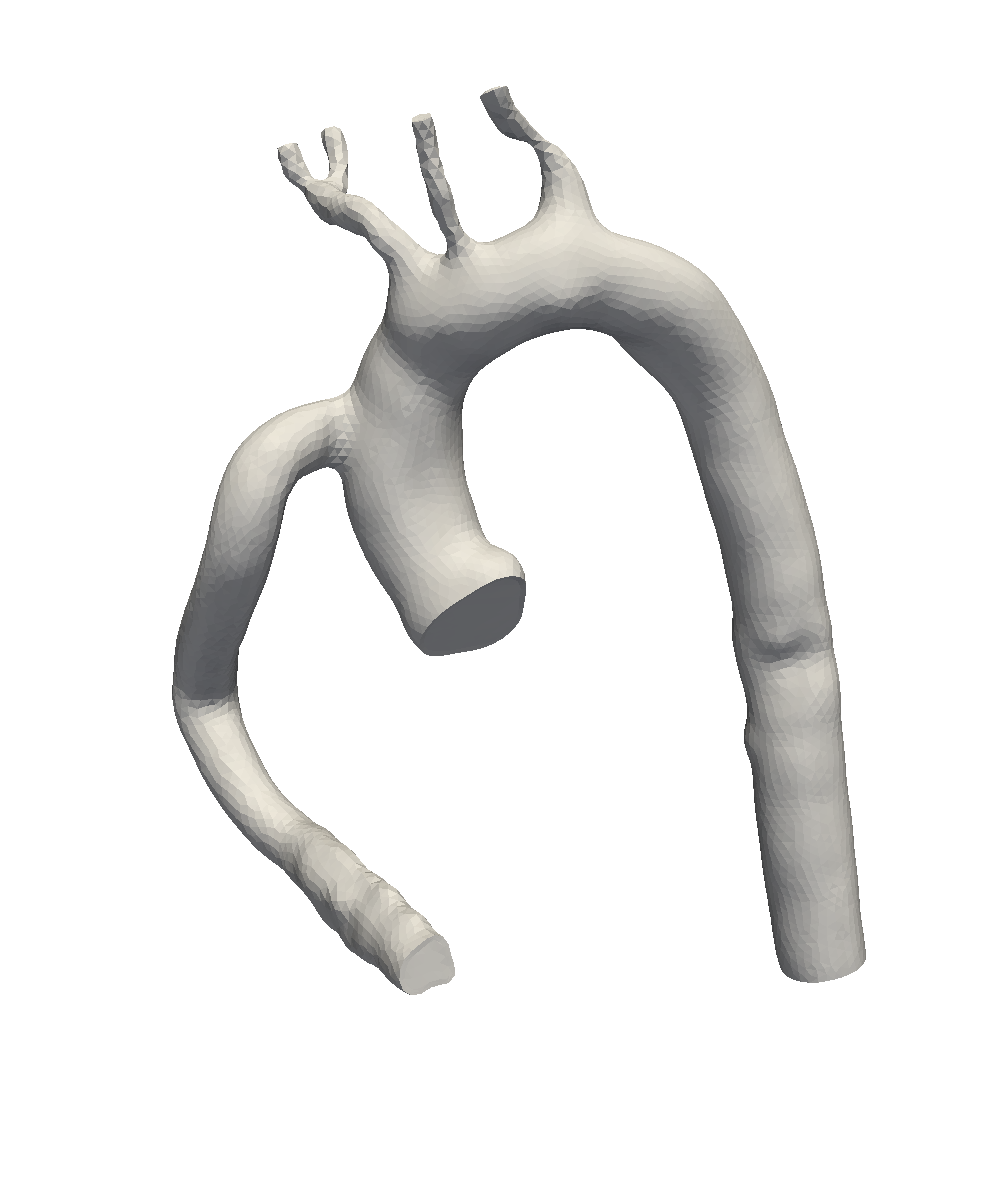}
\put(30,100){\small{b)}}
      \end{overpic}
\caption{Patient specific aorta models obtained from CT images: a) Pre-surgery configuration, b) Post-surgery configuration.}\label{fig:geom}
\end{figure}

Real patient-specific aorta models were reconstructed from CT images by using the open source medical image analysis software 3D Slicer\textsuperscript{\textregistered} (\url{http://www.slicer.org}). The models include the ascending aorta, brachiocephalic artery, right subclavian artery, right common carotid artery, left common carotid artery, left subclavian artery and descending aorta, and, in the post-surgery configuration, the outflow cannula of the LVAD device as well, as shown in Fig. \ref{fig:geom}.

\subsection{The full order model (FOM)}\label{sec:FOM}
In this Section we briefly introduce the mathematical model, i.e. the incompressible Navier-Stokes equations, with proper boundary conditions as well as the space and time discretization adopted. 

\subsubsection{The mathematical problem: Navier-Stokes equations}
We consider the motion of the blood in a time-independent domain $\Omega$
over a time interval of interest $(t_0,t^\star]$. The flow is described by the incompressible Navier-Stokes equations
\begin{align}
\rho\,\partial_t \bm{u} + \rho\,
\nabla \cdot \left(\bm{u} \otimes \bm{u}\right) -
\nabla \cdot \bm{\sigma} & = 0\quad \mbox{ in }\Omega \times (t_0,t^\star],\label{eq:ns-mom}\\
\nabla \cdot \bm{u} & = 0\quad\, \mbox{ in }\Omega \times (t_0,t^\star],\label{eq:ns-mass}
\end{align}
endowed with proper boundary conditions.
$\rho = 1060$ Kg/m$^3$ is the blood density, $\bm{u}$ is the blood velocity, $\partial_t$ denotes the time derivative, $\bm{\sigma}$ is the Cauchy stress tensor. 
 Equation (\ref{eq:ns-mom}) represents the conservation of the linear momentum, while eq. (\ref{eq:ns-mass}) represents the conservation of the mass. In this work, the blood is considered as a Newtonian fluid and $\bm{\sigma}$ can be written as
\begin{equation}\label{eq:newtonian}
\bm{\sigma}(\bm{u}, p) = -p \mathbf{I} +\mu (\nabla\bm{u} + \nabla\bm{u}^T),
\end{equation}
where $p$ is the pressure and $\mu = 0.004$ Pa $\cdot$ s is the blood \emph{dynamic} viscosity. For the sake of convenience, we also define the viscous stress tensor $\bm{\tau}$ as follows
\begin{equation}\label{eq:newtonian2}
\bm{\tau}(\bm{u}) = \mu (\nabla\bm{u} + \nabla\bm{u}^T).
\end{equation}
Notice that by plugging \eqref{eq:newtonian} into eq.~\eqref{eq:ns-mom}, eq.~\eqref{eq:ns-mom}
can be rewritten as
\begin{align}
\rho\, \partial_t \bm{u} + \rho\,
\nabla \cdot \left(\bm{u} \otimes \bm{u}\right) + \nabla p - \mu\,\Delta\bm{u}  = \bm{f}\mbox{ in }\Omega \times (t_0,t^\star].\label{eq:ns-lapls-1}
\end{align}

In order to investigate the blood flow patterns, we introduce the Wall Shear Stress (WSS) defined in the following way
\begin{equation}\label{eq:WSS}
WSS = \bm{\tau}_w \cdot \bm{n},
\end{equation}
where $\bm{n}$ is the unit normal vector and $\bm{\tau}_w$ is the tangential component of the wall viscous stress tensor. When the flow is pulsatile, it is useful to make reference to the Time Averaged WSS (TAWSS), 
\begin{equation}\label{eq:TAWSS}
TAWSS = \dfrac{1}{T} \int_0^T WSS\ dt.
\end{equation}

Finally, in order to characterize the flow regime under consideration, we define the Reynolds number as
\begin{equation}\label{eq:re}
Re = \frac{U L}{\nu},
\end{equation}
where $\nu=\mu/\rho$ is the \emph{kinematic} viscosity of the blood, and $U$ and $L$ are characteristic macroscopic velocity and length, respectively. For a blood flow in a cylindrical vessel, $U$ is the mean sectional velocity and $L$ is the diameter. 

\subsubsection{Boundary conditions}
Experimental measurements obtained by the RHC and ECHO tests are reported in Tables \ref{tab:data_pre} and \ref{tab:data_post} for pre-surgery and post-surgery configuration respectively. They are used in order to enforce realistic boundary conditions. For clinical reasons, four different tests are available for the post-surgery configurations whilst only one for the pre-surgery configuration. Note that RCH and ECHO tests provided measurements related to the polmonary circulation as well. However, these data are not reported because they do not affect the model used in this work, that deals with the systemic compartment only. In Table \ref{tab:area_k} we report the values of boundaries cross-sectional areas.

\begin{table}
\begin{center}
\begin{tabular}{|c|c|c|c|c|}
\hline
$PAS$ [mmHg] & $PAD$ [mmHg] & $PAM$ [mmHg] & $CO$ [l/min] & $SV$ [ml]       \\
\hline
 108 & 66 & 78  & 5.63 & 55 \\
\hline
\end{tabular}
\end{center}
\caption{Pre-surgery configuration: experimental data obtained by the RHC and ECHO tests. $PAS$ = systolic arterial pressure, $PAD$ = diastolic arterial pressure, $PAM$ = average arterial pressure, $CO$ = average cardiac flow rate, $SV$ = stroke volume.}\label{tab:data_pre}
\end{table}

\begin{table}
\begin{center}
\begin{tabular}{|c|c|c|c|}
\hline
& $PF$ [l/min]  & $\omega$ [rpm] & $PAM$ [mmHg]      \\
\hline
Test 1     & 4.1   & 5400  & 78  \\\hline
Test 2     & 4.2   & 5600  & 90       \\\hline
Test 3     & 4.5   & 6000  & 100      \\\hline
Test 4     & 5   & 5600  & 83     \\\hline
\end{tabular}
\end{center}
\caption{Post-surgery configuration: experimental data obtained by the RHC and ECHO tests. $PF$ = LVAD flow rate, $\omega$ = pump speed, $PAM$ = average arterial pressure.}\label{tab:data_post}
\end{table}

\begin{table}
\begin{center}
\begin{tabular}{|c|c|}
\hline
 & $A$ [cm$^2$]       \\\hline
Outflow cannula    & 1.3  \\\hline
Ascending aorta    & 6.42  \\\hline
Right subclavian artery      & 0.156       \\\hline
Right common carotid artery     & 0.246        \\\hline
Left common carotid artery      & 0.168       \\\hline
Left subclavian artery  &  0.446      \\\hline
Descending aorta    & 3.68  \\\hline
\end{tabular}
\end{center}
\caption{Values of boundaries cross-sectional areas.}\label{tab:area_k} 
\end{table}

\begin{figure}
\centering
\includegraphics[width=0.5\textwidth]{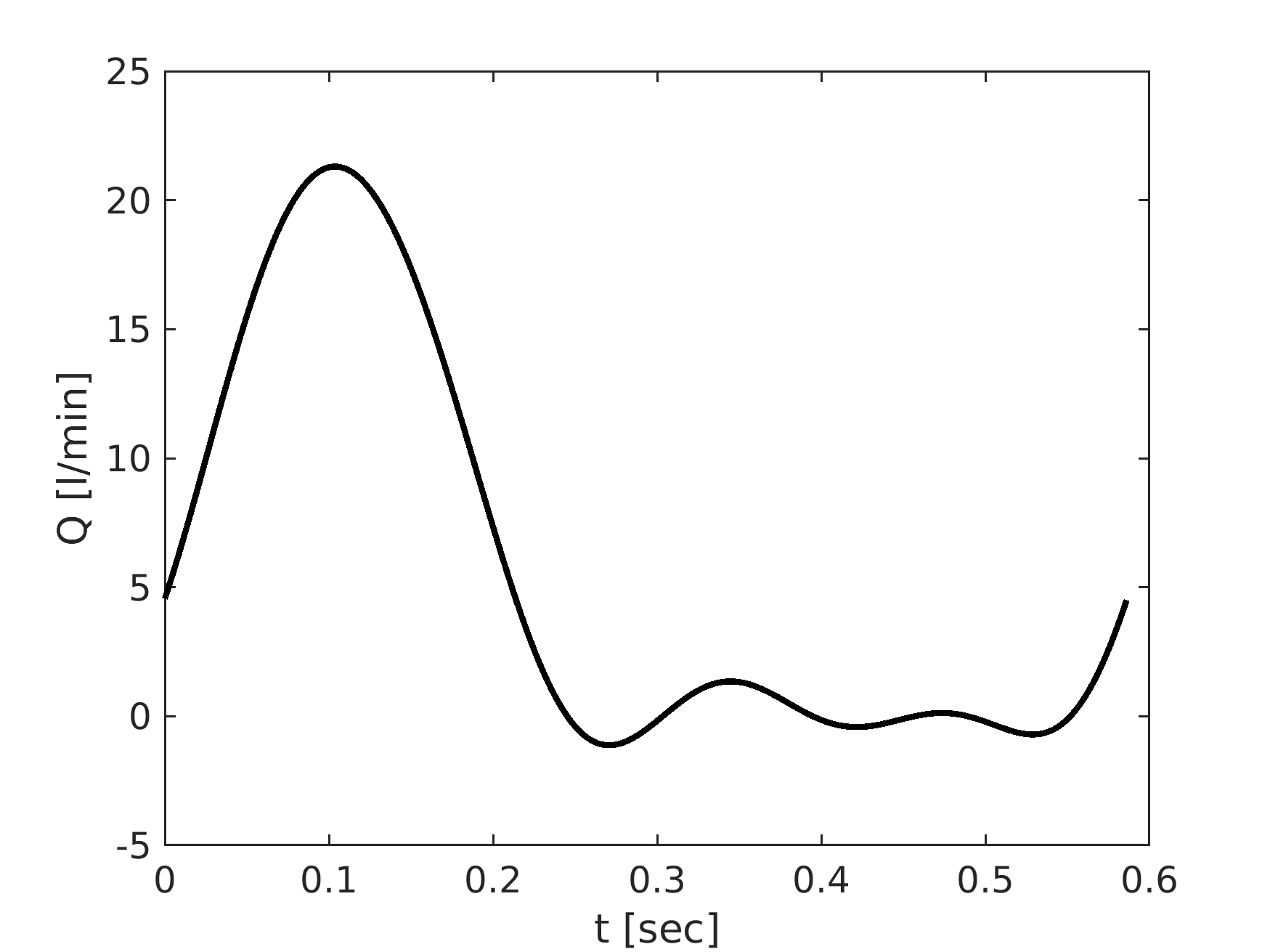}
\caption{Aortic inflow waveform enforced in the pre-surgery case.}
\label{fig:flowrate}
\end{figure}

In the pre-surgery configuration, a realistic flow rate $Q$ waveform was enforced on the ascending aorta section (Figure \ref{fig:flowrate}). The amplitude of the flow waveform has been set according to the average flow rate over the cardiac cycle, $CO$,
\begin{equation}
CO = \dfrac{1}{T} \int_0^T Q\ dt,
\end{equation}\label{eq:CO_average}
measured by the RHC test. The value of a period of the cardiac cycle, $T$, is obtained as
\begin{equation}\label{eq:T}
T = \dfrac{SV}{CO},
\end{equation}
where $SV$ is the stroke volume measured by the ECHO test.

\begin{figure}
\centering
\includegraphics[width=0.5\textwidth]{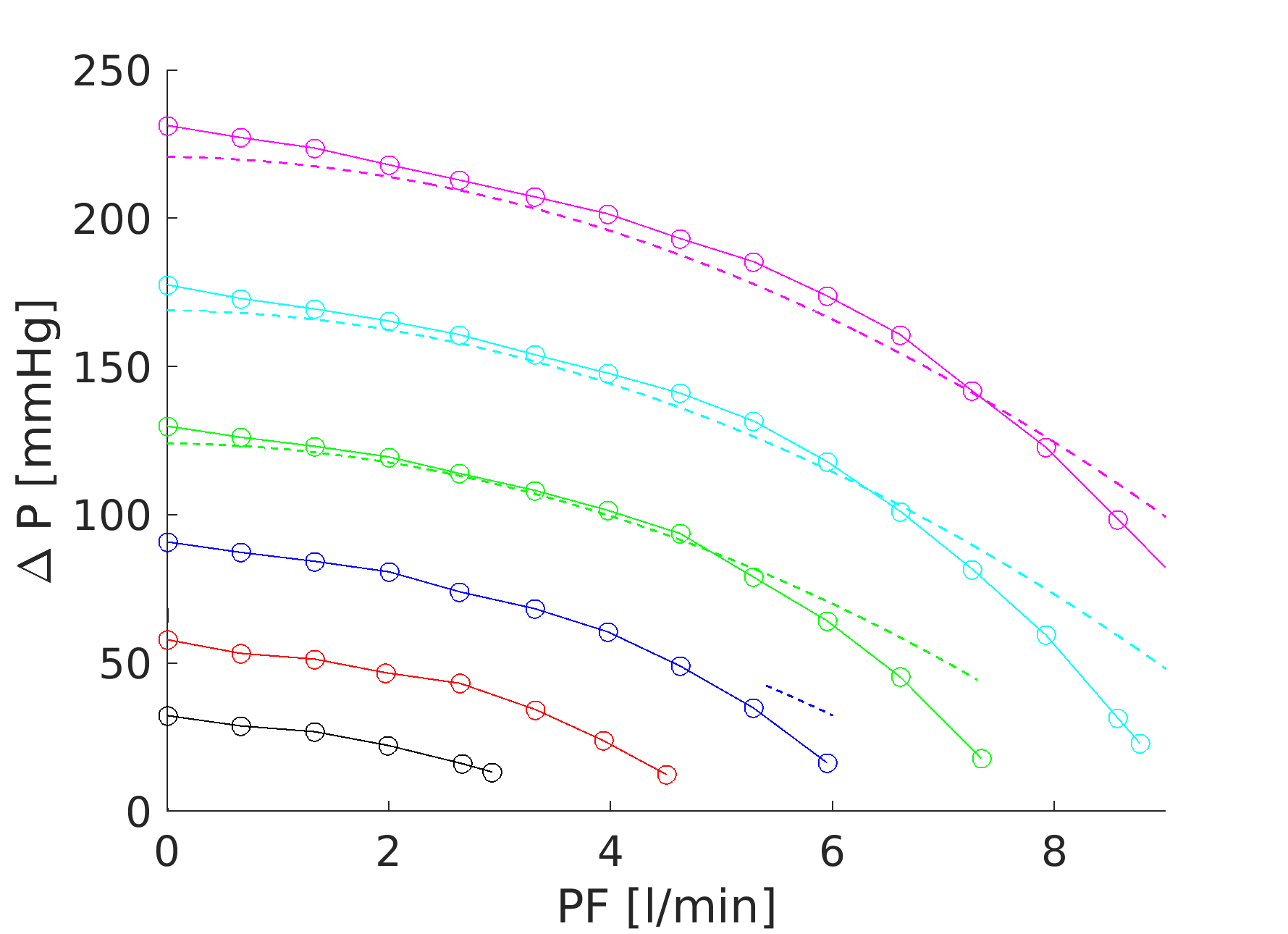}
\caption{Pressure head ($\Delta P$) - volume flow rate ($PF$) curves (continuous line with circles) and analytical fitting (dashed line) based on eq. \ref{eq:pump} for Heartmate 3$^\text{TM}$ \cite{HeartMate} pump at several pump speed values: $\omega = 3000$  rpm (black), $\omega = 4000$ rpm (red), $\omega = 5000$ rpm (blue), $\omega = 6000$ rpm (green), $\omega = 7000$ rpm (cyan), and $\omega = 8000$ rpm (magenta).}
\label{fig:pump_dynamics}
\end{figure}

On the other hand, in the post-surgery configuration, the LVAD flow rate, $PF$, has been used as inlet boundary condition applied to the outflow cannula section. Note that the aortic valve is closed during all the cardiac cycle, i.e. the cardiac flow rate is supplied by the LVAD device only and the ascending aorta section is treated as a wall. In Figure \ref{fig:pump_dynamics}, the pressure head ($\Delta P$)
- volume flow rate ($PF$) curves for the Heartmate 3$^\text{TM}$ Left Ventricular Assist System \cite{HeartMate} at several pump speed values $\omega$ are shown. The basic pump dynamics can, in principle, be described in the following way \cite{Shi2010}
\begin{equation}\label{eq:pump}
\Delta P = K_A \omega^2 + K_B \omega \cdot PF + K_C PF^2,
\end{equation}
where $K_A$, $K_B$ and $K_C$ are constants which depend on pump design. After some numerical experiments, we found that the coefficients given in Table \ref{tab:coeffs_fit} provide an acceptable fit as showed in Figure \ref{fig:pump_dynamics}. We note that a better agreement could be obtained by considering more complex equations in order to take into account nonlinear effects \cite{Shi2010}. Based on the analytical fitting \ref{eq:pump}, we can compute the values of $\Delta P$ for the all the tests under consideration (see Table \ref{tab:delta_P}). The analytical fitting \ref{eq:pump} will be also used in Sec. \ref{sec:res_ROM} in order to compute $\omega$ at given $PF$ and $\Delta P$.

\begin{table}
\begin{center}
\begin{tabular}{|c|c|c|}
\hline
$K_A$ [mmHg/rpm$^2$] & $K_B$ [mmHg $\cdot$ l/min/rpm] & $K_C$ [mmHg $\cdot$ l$^2$/rpm$^2$] \\
\hline
3.45e-6      & -5.9e-5  &  -1.45 \\\hline
\end{tabular}
\end{center}
\caption{Parameter settings for the pump dynamics (eq. \ref{eq:pump}).}\label{tab:coeffs_fit} 
\end{table}

\begin{table}
\begin{center}
\begin{tabular}{|c|c|}
\hline
  & $\Delta P$ [mmHg] \\
\hline
Test 1      &  75\\\hline
Test 2      &  81.3\\\hline
Test 3      &  93.3\\\hline
Test 4      &  70.4\\\hline
\end{tabular}
\end{center}
\caption{$\Delta P$ values based on eq. \ref{eq:pump} for all the tests under consideration.}\label{tab:delta_P} 
\end{table}

\begin{figure}
\centering
\includegraphics[width=0.4\textwidth]{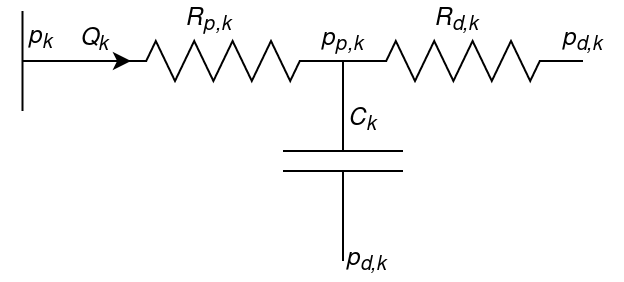}
\caption{Three-element Windkessel model for the generic outlet $k$}
\label{fig:RCR}
\end{figure}

Outflow boundary conditions were applied at each outlet of the model, right subclavian artery, right common carotid artery, left common carotid artery, left subclavian artery and descending aorta, by using a three-element Windkessel RCR model \cite{Westerhof2008}. The Windkessel model consists of a proximal resistance $R_{p,k}$, a compliance $C_k$, and a distal resistance $R_{d,k}$, for each outlet $k$ (Figure \ref{fig:RCR}). The downstream pressure, $p_k$, is expressed through the following DAE system: 
\begin{equation}\label{eq:RCR}
\begin{cases}
C_k \dfrac{dp_{p,k}}{dt} + \dfrac{p_{p,k} - p_{d,k}}{R_{d,k}} = Q_k, \\ \\
p_k - p_{p,k} = R_{p,k}Q_k, \\
\end{cases}
\end{equation}
where $Q_k$ is the flow rate, and $p_{p,k}$ and $p_{d,k}$ are the proximal and the distal pressure, respectively.
The total resistance, $R_k = R_{p,k} + R_{d,k}$, was evaluated according to the rules for a parallel circuit
\begin{equation}
R_k = RVS \dfrac{\sum_k A_k}{A_k},
\end{equation}
where $A_k$ is the cross-sectional area (see Table \ref{tab:area_k}) and $RVS$ is the systemic vascular resistance estimated as follows
\begin{equation}\label{eq:RVS}
RVS =
\begin{cases}
\dfrac{PAM}{CO}, & \mbox{in the pre-surgery case} \\ \\
\dfrac{PAM}{PF}, & \mbox{in the post-surgery case} \\
\end{cases}
\end{equation}
where $PAM$ is the average arterial pressure measured by the RHC test (see Tables \ref{tab:data_pre} and \ref{tab:data_post}). 
For each outlet $k$, we assumed \cite{Laskey1990}
\begin{equation}
\dfrac{R_{p,k}}{R_k} = 0.056.
\end{equation}\label{eq:ratio}
On the other hand, the aortic compliance, $C$, can be estimated as follows \cite{Bulpitt1999}:
\begin{equation}\label{eq:ratioC}
C = \dfrac{PAS - PAD}{SV}
\end{equation}
where $PAS$ and $PAD$ are the systolic and the diastolic pressure measured by the RCH test in the pre-surgery configuration, respectively (see Table \ref{tab:data_pre}). It should be noted that such value is also used in the post-surgery configuration. Finally, the compliance $C_k$ related to the outlet $k$ was evaluated according to the rules for a parallel circuit
\begin{equation}
C_k = C \dfrac{A_k}{\sum_k A_k}.
\end{equation}

Table \ref{tab:data_computed_pre} shows the values of $T$, $RVS$ and $C$ computed by using eqs. \eqref{eq:T}, \eqref{eq:RVS}, and \eqref{eq:ratioC}, respectively, for the pre-surgery configuration. Table \ref{tab:data_computed_post} shows the values of $RVS$ computed by using eq. \eqref{eq:RVS} for the post-surgery configuration. Finally, tables \ref{tab:bc_pre} and \ref{tab:bc_post} report the values of Windkessel coefficients for the pre-surgery and post-surgery configurations, respectively.

\begin{table}
\begin{center}
\begin{tabular}{|c|c|c|}
\hline
$T$ [s]  & $RVS$ [dyne $\cdot$ s/cm$^5$]  & $C$ [cm$^5$/dyne]      \\
\hline
0.586   & 1105  & 9.85e-4  \\
\hline
\end{tabular}
\end{center}
\caption{Pre-surgery configuration: quantities computed by the experimental data reported in Tab. \ref{tab:data_pre}. $T$ = period of the cardiac cycle (eq. \ref{eq:T}), $RVS$ = system vascular resistance (eq. \ref{eq:RVS}), $C$ = aortic compliance (eq. \ref{eq:ratioC}).}\label{tab:data_computed_pre}
\end{table}

\begin{table}
\begin{center}
\begin{tabular}{|c|c|}
\hline
    & $RVS$ [dyne $\cdot$ s/cm$^5$]       \\
\hline
Test 1   & 1522   \\\hline
Test 2   & 1714    \\\hline
Test 3   & 1778   \\\hline
Test 4   & 1328    \\\hline
\end{tabular}
\end{center}
\caption{Post-surgery configuration: system vascular resistance (eq. \ref{eq:RVS}) computed by the experimental data reported in Tab. \ref{tab:data_post}.}\label{tab:data_computed_post}
\end{table}

\begin{table}
\begin{center}
\begin{tabular}{|c|c|c|c|}
\hline
$k$  & $R_{p,k}$ [dyne $\cdot$ s/cm$^5$] & $R_{d,k}$ [dyne $\cdot$ s/cm$^5$]  & $C_k$ [cm$^5$/dyne]
\\\hline
Right subclavian artery      &  1.84e3 & 3.11e4 & 3.26e-5       \\\hline
Right common carotid artery       &  1.23e3 & 2.07e4 & 5.16e-5      \\\hline
Left common carotid artery        &  1.78e3 & 3.01e4 & 3.52e-5       \\\hline
Left subclavian artery   &  7.09e2 & 1.19e4 & 9.35e-5      \\\hline
Descending aorta     & 7.8e1  & 1.31e3 & 7.72e-4\\\hline
\end{tabular}
\end{center}
\caption{Pre-surgery configuration Windkessel coefficients: proximal resistance $R_{p,k}$, distal resistance $R_{d,k}$ and compliance $C_{k}$, for each outlet $k$.}\label{tab:bc_pre} 
\end{table}

\begin{table}
\begin{center}
\begin{tabular}{|c|c|c|c|}
\hline
& $k$  & $R_{p,k}$ [dyne $\cdot$ s/cm$^5$] & $R_{d,k}$ [dyne $\cdot$ s/cm$^5$]
\\\hline
 Test 1& Right subclavian artery      &  2.56e3 & 4.32e4       \\\cline{2-4}
       & Right common carotid artery       &  1.63e3 & 2.74e4       \\\cline{2-4}
       & Left common carotid artery        &  2.38e3 & 4e4       \\\cline{2-4}
       & Left subclavian artery   &  8.96e2 & 1.51e4       \\\cline{2-4}
       & Descending aorta     & 1.08e2  & 1.83e3  \\\hline
 Test 2      & Right subclavian artery       &  2.88e3 & 4.86e4        \\\cline{2-4}
       & Right common carotid artery       &   1.83e3 & 3.08e4        \\\cline{2-4}
       & Left common carotid artery        &  2.68e3 & 4.51e4        \\\cline{2-4}
       & Left subclavian artery   &  1.01e3 & 1.7e4        \\\cline{2-4}
       & Descending aorta      &  1.22e2 & 2.06e3 \\\hline
 Test 3      & Right subclavian artery       &  2.99e3 & 5.05e4        \\\cline{2-4}
       & Right common carotid artery       &  1.9e3 & 3.2e4       \\\cline{2-4}
       & Left common carotid artery       &  2.78e3 & 4.68e4       \\\cline{2-4}
       & Left subclavian artery  &  1.04e3 & 1.76e4       \\\cline{2-4}
       & Descending aorta     &  1.27e2 & 2.14e3 \\\hline
 Test 4      & Right subclavian artery       &  2.19e3 & 3.68e4        \\\cline{2-4}
       & Right common carotid artery       &  1.39e3 & 2.33e4       \\\cline{2-4}
       & Left common carotid artery       &  2.03e3 & 3.42e4       \\\cline{2-4}
       & Left subclavian artery  &  7.64e2 & 1.29e4      \\\cline{2-4}
       & Descending aorta    &  9.25e1 & 1.56e3 \\\hline
\end{tabular}
\end{center}
\caption{Post-surgery configuration Windkessel coefficients: proximal resistance $R_{p,k}$ and distal resistance $R_{d,k}$, for each outlet $k$.}\label{tab:bc_post} 
\end{table}

\subsubsection{Space and temporal discretization}\label{sec:disc}
For the space discretization of problems \eqref{eq:ns-mass}-\eqref{eq:ns-lapls-1}, we adopt the Finite Volume (FV) approximation that is derived directly from the integral form of the governing equations. We have used
the finite volume C++ library OpenFOAM\textsuperscript{\textregistered} \cite{Weller1998}. We partition the computational domain $\Omega$ (i.e., the patient-specific geometrical models in Figure \ref{fig:geom}) into cells or control volumes $\Omega_i$, with $i = 1, \dots, N_{c}$, where $N_{c}$ is the total number of cells in the mesh. Let  \textbf{A}$_j$ be the surface vector of each face of the control volume.

\begin{figure}[h!]
\centering
\includegraphics[width=0.5\textwidth]{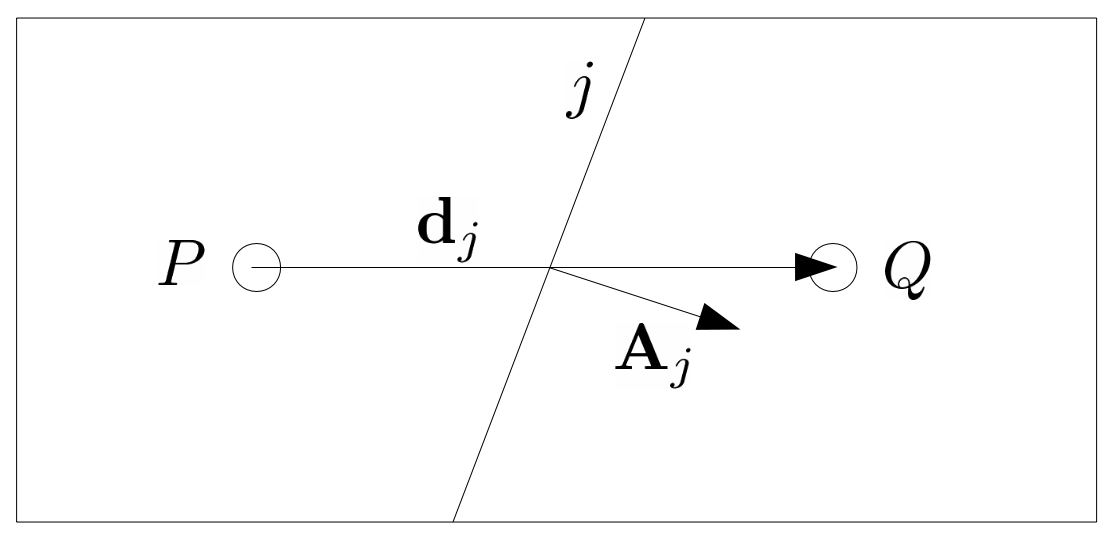}
\caption{Close-up view of two non-orthogonal control volumes in a 2D configuration.}
\label{fig:gradient_image}
\end{figure}

The integral form of eq.~\eqref{eq:ns-lapls-1} for each volume $\Omega_i$ is given by:

\begin{align}\label{eq:evolveFVtemp-1.1}
\rho \int_{\Omega_i} \dfrac{\partial \u}{\partial t} d\Omega + \rho\, \int_{\Omega_i} \nabla \cdot \left(\u \otimes \u\right) d\Omega - \mu \int_{\Omega_i} \Delta\u d\Omega + \int_{\Omega_i}\nabla p d\Omega  = 0.
\end{align}
By applying the Gauss-divergence theorem, eq.~\eqref{eq:evolveFVtemp-1.1} becomes:

\begin{align}\label{eq:evolveFV-1.1}
\rho \int_{\Omega_i} \dfrac{\partial \u}{\partial t}d\Omega + \rho\, \int_{\partial \Omega_i} \left(\u \otimes \u\right) \cdot d\textbf{A} - \mu \int_{\partial \Omega_i} \nabla\u \cdot d\textbf{A} + \int_{\partial \Omega_i}p d\textbf{A}  = 0.
\end{align}
Each term in eq.~\eqref{eq:evolveFV-1.1} is approximated as follows:

\begin{itemize}
\item[-] \textit{Gradient term}:

\begin{align}\label{eq:grad}
\int_{\partial \Omega_i}p d\textbf{A} \approx \sum_j^{} p_j \textbf{A}_j,
\end{align}
where $p_j$ is the value of the pressure relative to centroid of the $j^{\text{th}}$ face. 
The face center pressure values $p_j$ are obtained from the cell center values by means of a linear interpolation scheme. 


\item[-] \textit{Convective term}:
\begin{align}\label{eq:conv}
\int_{\partial \Omega_i} \left(\u \otimes \u\right) \cdot d\textbf{A} \approx \sum_j^{} \left(\u_j \otimes \u_j\right) \cdot \textbf{A}_j = \sum_j^{} \varphi_j \u_j, \quad \varphi_j = \u_j \cdot \textbf{A}_j,
\end{align}
where $\u_j$ is the fluid velocity relative to the centroid of each control volume face. In \eqref{eq:conv}, $\varphi_j$ is the convective flux associated to $\u$ through face $j$ of the control volume. The convective flux at the cell faces is obtained by a linear interpolation of the values from the adjacent cells. Also $\u$ needs to be approximated at cell face $j$ in order to get the face value $\u_j$. Different interpolation methods can be applied: central, upwind, second order upwind and blended differencing schemes \cite{jasak1996error}. In this work, we make use of a second order upwind scheme.
\item[-] \textit{Diffusion term}:
\begin{align}
\int_{\partial \Omega_i} \nabla\u \cdot d\textbf{A} \approx \sum_j^{} (\nabla\u)_j \cdot \textbf{A}_j, \nonumber
\end{align}
where $(\nabla\u)_j$ is the gradient of $\u$ at face $j$.
We are going to briefly explain how $(\nabla\u)_j$ is approximated with
second order accuracy on structured, orthogonal meshes. Let $P$ and $Q$ be two neighboring control volumes.
The term $(\nabla\u)_j$ is evaluated by subtracting
 the value of velocity at the cell centroid on the $P$-side of the face, denoted with $\u_P$,
 from the value of velocity at the centroid on the $Q$-side, denoted with $\u_Q$,
 and dividing by the magnitude of the distance vector $\textbf{d}_j$ connecting the two cell centroids:
\begin{align}
(\nabla\u)_j \cdot \textbf{A}_j = \dfrac{\u_Q - \u_P}{|\textbf{d}_j|} |\textbf{A}_j|. \nonumber
\end{align}
For non-structured, non-orthogonal meshes
(see Fig.~\ref{fig:gradient_image}), that are used in this work, an explicit non-orthogonal correction has to be added to the orthogonal component
in order to preserve second order accuracy. See \cite{jasak1996error} for details.
\end{itemize}

A partitioned approach has been used to deal with the pressure-velocity coupling. In particular a Poisson equation for pressure has been used. This is obtained by taking the divergence of the momentum equation \eqref{eq:ns-lapls-1} and exploiting the divergence free constraint \eqref{eq:ns-mass}:
\begin{equation}\label{eq:Poisson}
\Delta p = -\nabla \left(\u \otimes \u\right).
\end{equation}

The segregated algorithms available in OpenFOAM\textsuperscript{\textregistered} are SIMPLE \cite{SIMPLE} for steady-state problems, and PISO \cite{PISO} and PIMPLE \cite{PIMPLE} for transient problems. For this work, we choose the PISO algorithm. 

To discretize in time the equation \eqref{eq:evolveFV-1.1}, let $\Delta t \in \mathbb{R}$, $t^n = t_0 + n \Delta t$, with $n = 0, ..., N_T$ and $t^\star = t_0 + N_T \Delta t$. Moreover, we denote by $\u^n$ the approximation of the flow velocity at the time $t^n$. We adopt Backward Differentiation Formula of order 1 (BDF1),  see e.g. \cite{quarteroni2007numerical}. Given $\u^n$, for $n \geq 0$, we have, respectively, 
\begin{equation}\label{eq:BDF1_disc}
\partial_t \u \approx \dfrac{\u^{n+1} - \u^{n}}{\Delta t},
\end{equation}

Finally, a first-order scheme is also used for the discretization of the RCR Windkessel model \eqref{eq:RCR}: 
\begin{equation}\label{eq:RCR_disc}
\begin{cases}
C_k \dfrac{p_{p,k}^{n+1} - p_{p,k}^{n}}{\Delta t} + \dfrac{p_{p,k}^{n+1}}{R_{d,k}} = Q_k^{n}, \\ \\
p_k^{n+1} - p_{p,k}^{n+1} = R_{p,k}Q_k^{n}, \\
\end{cases}
\end{equation}
where we assumed $p_{d,k} = 0$.

\subsubsection{Multi-scale coupling}
The coupling process between the three-dimensional flow model and lumped Windkessel model can be summarized as follows:
\begin{enumerate}
\item At $t^n$, we know $\u^{n}$ and thus $Q_k^{n}$. Then we calculate $p_k^{n+1}$ by eq. \eqref{eq:RCR};
\item We solve the problem \eqref{eq:evolveFV-1.1}-\eqref{eq:Poisson} to obtain $\u^{n+1}$ and $Q_k^{n+1}$.
\end{enumerate}

\subsection{The reduced order model (ROM)}\label{sec:ROM}
The reduced order model we propose is the so-called \emph{proper orthogonal decomposition with interpolation}. In Sec. \ref{sec:ROM1} we provide a brief description of a such technique. 

\subsubsection{Proper orthogonal decomposition with interpolation}\label{sec:ROM1}
Proper orthogonal decomposition (POD) is a technique widely used within the reduced order modeling (ROM) framework for the study of parametric problems. POD allows to extract, from a set of high-dimensional snapshots, the basis minimizing the error between the original snapshots and their orthogonal projection. The data-driven approach here used is based only on data and does not require knowledge about the governing equations that describe the system. It is also non-intrusive, i.e. no modification of the simulation software is carried out. On the other hand, there are works that use non-intrusive methods that are not data-driven (see, e.g., \cite{Zou2017}). The original snapshots are projected onto the POD space in order to reduce their dimensionality. Then the solution manifold is approximated using an interpolation technique. Several examples of applications based on this so-called POD with interpolation (PODI) \cite{Bui-Thanh2004} techique can be found in literature, in a wide range of contexts: naval engineering problems \cite{demo2018shape, demo2019marine, demo2018b, Demo2019}, automotive \cite{Salmoiraghi2018, Dolci2016}, aeronautics \cite{Ripepi2018}. We also cite \cite{Forti2014} where a coupling with isogeometric analysis is performed.


We are going to describe briefly the computation of the POD modes. We consider a problem with $N$ degrees of freedom. Let $\bm{\varphi}_i$, with $i = 1, \dots , N_s$, be the snapshots related to a generic variable of interest collected by solving the high-fidelity problem, with different values of the input parameters $\bm{\pi}_i$, resulting in $N_s$ input-output pairs ($\bm{\pi}_i$, $\bm{\varphi}_i$). The snapshots matrix $\bm{S}$ is built arranging the snapshots as columns, such that $\bm{S} = [\bm{\varphi}_1, \bm{\varphi}_2, \dots, \bm{\varphi}_{N_s}]$. By applying the singular value decomposition to this matrix, we have:
\begin{equation}\label{eq:svd}
\bm{S} = \bm{U}\bm{\Sigma}\bm{V}^* \approx \bm{U}_k\bm{\Sigma}_k{\bm{V}_k}^*,
\end{equation}
where $\bm{U} \in \mathcal{A}^{N \times N_s}$ is the unitary matrix containing the left-singular vectors, $\bm{\Sigma} \in \mathcal{A}^{N_s \times N_s}$  is the diagonal matrix containing the singular values $\lambda_i$, and $\bm{V} \in \mathcal{A}^{N_s \times N_s}$, with the symbol $^*$  denoting the conjugate transpose. The leftsingular vectors, namely the columns of $\bm{U}$, are the so-called POD modes.
We can keep the first $k$ modes to span the optimal space with dimension $k$ to represent the snapshots. 
The matricies $\bm{U}_k \in \mathcal{A}^{N \times k}$, $\bm{\Sigma}_k \in \mathcal{A}^{k \times k}$, $\bm{V}_k \in \mathcal{A}^{N_s \times k}$ in Eq. \ref{eq:svd} are the truncated matrices with rank $k$. 

After constructing the POD space, we can project the original snapshots onto this space. We compute $\bm{C} \in \mathcal{R}^{k \times N_s}$ as $\bm{C} = {\bm{U}_k}^T \bm{S}$, where
the columns of $\bm{C}$ are the so-called modal coefficients. We express the input snapshots as a linear combination of the modes using such coefficients. Then, we have:
\begin{equation}\label{eq:svd3}
\bm{\varphi}_i = \sum_{j=1}^{N_s} \alpha_{ji}\phi_j \approx \sum_{j=1}^{k}  \alpha_{ji}\phi_j, \quad \forall \in [1,2, \dots, N_s],
\end{equation}
where $\alpha_{ji}$ are the elements of $\bm{C}$. Finally, we obtain the ($\bm{\pi}_i$, $\alpha_i$) pairs, for $i = 1, 2, \dots, N_s$, that sample the solution manifold in the parametric
space. We are able to interpolate the modal coefficients $\alpha$ and for any new parameter approximate the new coefficients. At the end, we compute
the high-dimensional solution by projecting back the (approximated) modal coefficients to the original space by using Equation \ref{eq:svd3}. We remark that the procedure can be repeated for several variables of interests. Furthermore, it is not necessary for such a variable to be an unknown of the original system (such as velocity and pressure); indeed, we will use the PODI technique not only for primal quantities, but also for derived quantities such as WSS.

Regarding the technical implementation of the PODI method, we use the Python package called EZyRB \cite{eazyrb}.





\section{Numerical results and discussion}\label{sec:results}
We validate the FOM model both for pre-surgery and post-surgery configuration in Sec. \ref{sec:FOM_valid}. Then, we investigate the performance of the ROM model in Sec. \ref{sec:res_ROM}.

\subsection{FOM validation}\label{sec:FOM_valid}

The number of PISO loops and non-orthogonal correctors has been fixed to 2 for all the simulations. The following solvers have provided a good compromise between stability, accuracy, and numerical cost. The linear algebraic system
associated with eq. \eqref{eq:evolveFV-1.1} is solved using an iterative solver with symmetric Gauss-Seidel smoother. Moreover, for Poisson problem \eqref{eq:Poisson}, we use Geometric Agglomerated Algebraic Multigrid Solver GAMG with the Gauss-Seidel smoother. The required accuracy is 1e-6 at each time step.

\subsubsection{Mesh convergence}
\begin{figure}[]
\centering
 \begin{overpic}[width=0.4\textwidth]{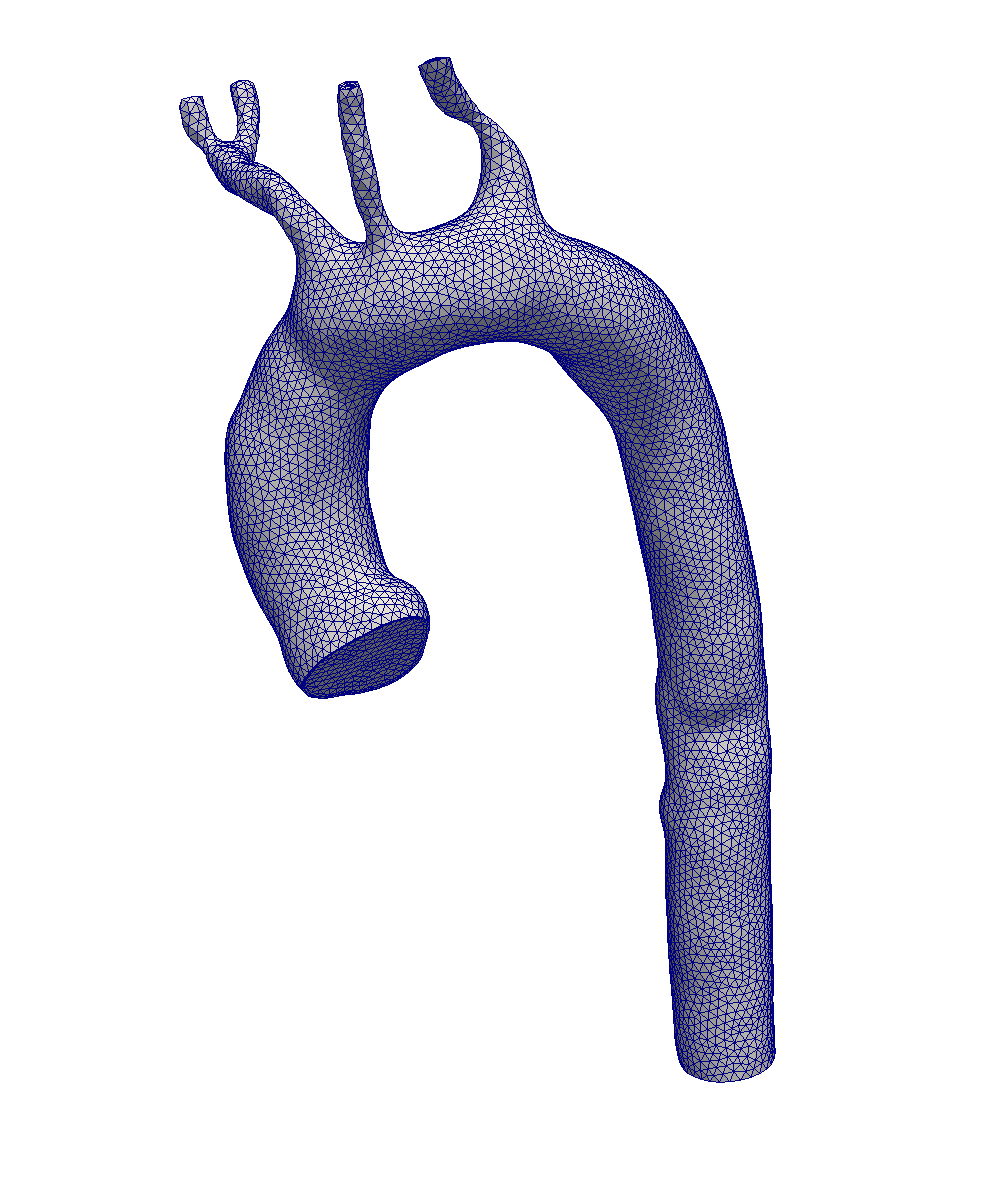}
\put(30,100){\small{a)}}
      \end{overpic}
 \begin{overpic}[width=0.4\textwidth]{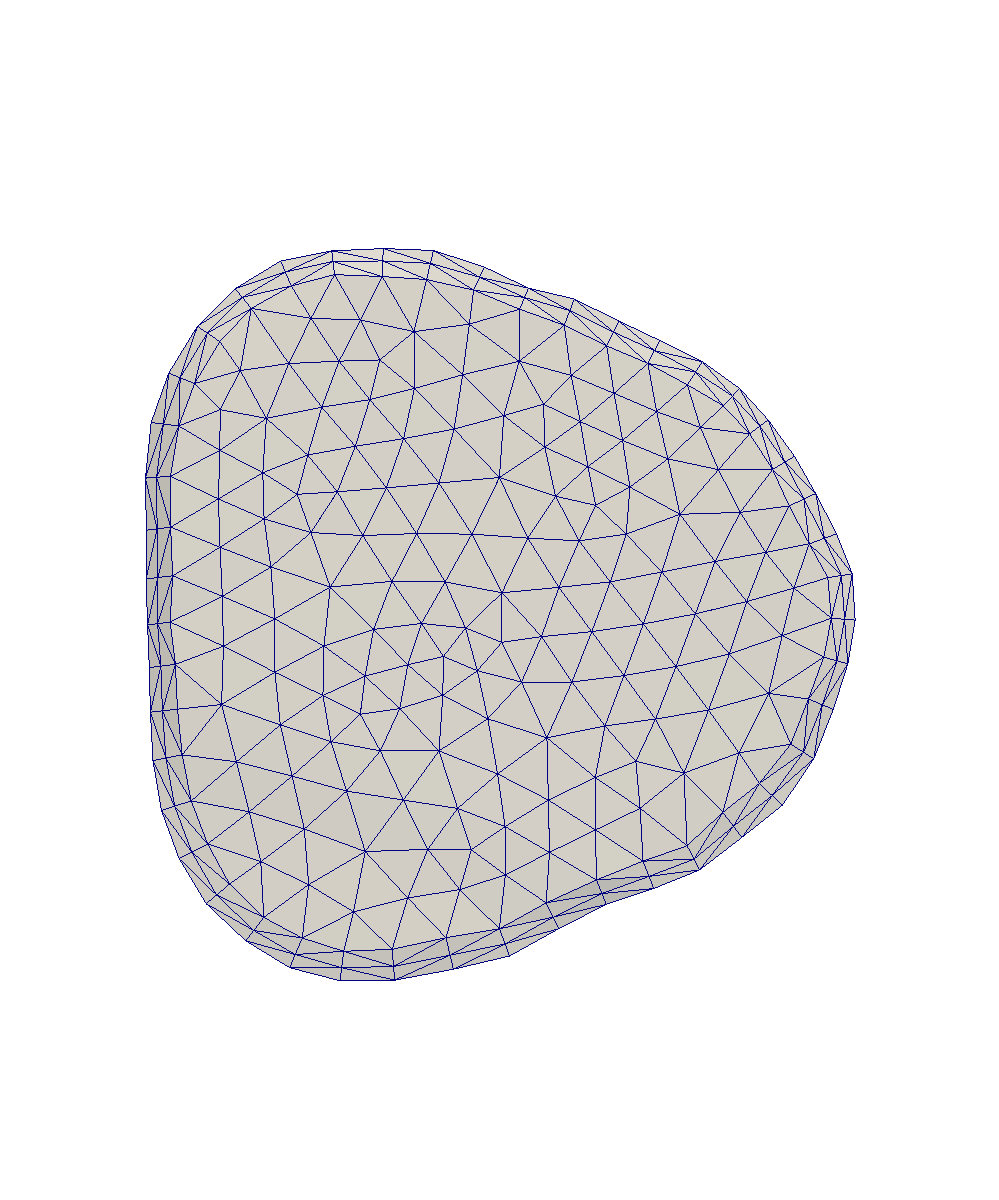}
\put(30,100){\small{b)}}
      \end{overpic}
\caption{View of the mesh $230k$: a) aortic wall, b) a section next to the aortic inlet.}\label{fig:mesh}
\end{figure}

In order to obtain grid independent solutions, we consider three meshes with tetrahedral elements. Table \ref{tab:mesh} reports name, minimum and maximum diameter, and number of cells for each mesh. Fig. \ref{fig:mesh} shows the mesh $230k$. All the meshes under consideration have very low values of average non-orthogonality (around $30 ^\circ$) and skewness (around 1). The estimation of the Reynolds number is based on the diameter computed by considering the inlet areas, i.e. the ascending aorta (\emph{ao}) section in the pre-surgery configuration and the outflow cannula section (\emph{oc}) in the post-surgery configuration, as circular areas. We have

\begin{equation}\label{eq:Re_pre}
Re = \dfrac{\dfrac{Q}{A_{ao}}\sqrt{\dfrac{4A_{ao}}{\pi}}}{\nu}
\end{equation}

\begin{equation}\label{eq:Re_post}
Re = \dfrac{\dfrac{PF}{A_{oc}}\sqrt{\dfrac{4A_{oc}}{\pi}}}{\nu}
\end{equation}
for the pre-surgery and post-surgery configuration, respectively. We carry out the mesh convergence study for the pre-surgery configuration because it is more critical with respect to the the post-surgery configuration being characterized by a greater Reynolds number $Re$ as showed in Table \ref{tab:Re}. Moreover, note that in the pre-surgery configuration the Reynolds number is time dependent, with  $0 \leq Re \leq 4200$.  

\begin{table}
\begin{center}
\begin{tabular}{|c|c|c|c|}
\hline
mesh name & $h_{min}$ & $h_{max}$ & No. of cells    \\
\hline
$230k$ & 5.8e-4 & 3e-3  & 228296 \\\hline
$415k$ & 5.6e-4 & 2.5e-3  & 414192 \\\hline
$2000k$ & 4.4e-4 & 1.5e-3  & 1993514 \\\hline
\end{tabular}
\end{center}
\caption{Name, minimum diameter $h_{min}$ , maximum diameter $h_{max}$, and number of cells for
all the meshes used for the convergence study.}\label{tab:mesh}
\end{table}

\begin{table}
\begin{center}
\begin{tabular}{|c|c|}
\hline
 & $Re$     \\
\hline
Pre-surgery &  [0, 4200]\\\hline
Post-surgery: test 1 &  1818\\\hline
Post-surgery: test 2 &  1862\\\hline
Post-surgery: test 3 &  1995\\\hline
Post-surgery: test 4 &  2217\\\hline
\end{tabular}
\end{center}
\caption{Reynolds number $Re$ for all the flow regimes under consideration.}\label{tab:Re}
\end{table}

\begin{figure}[]
\centering
 \begin{overpic}[width=0.4\textwidth]{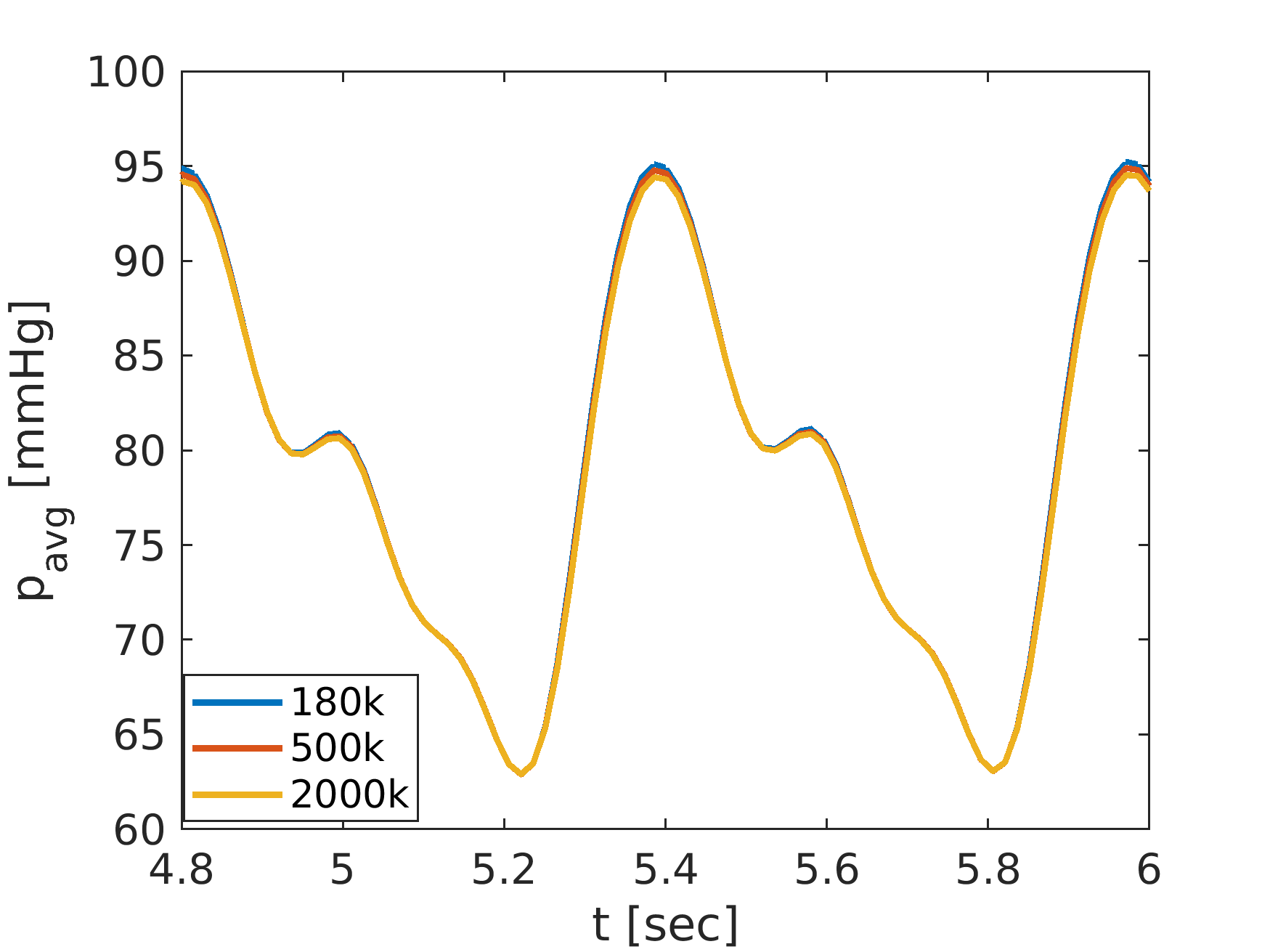}
\put(45,72){\small{a)}}
      \end{overpic}
 \begin{overpic}[width=0.41\textwidth]{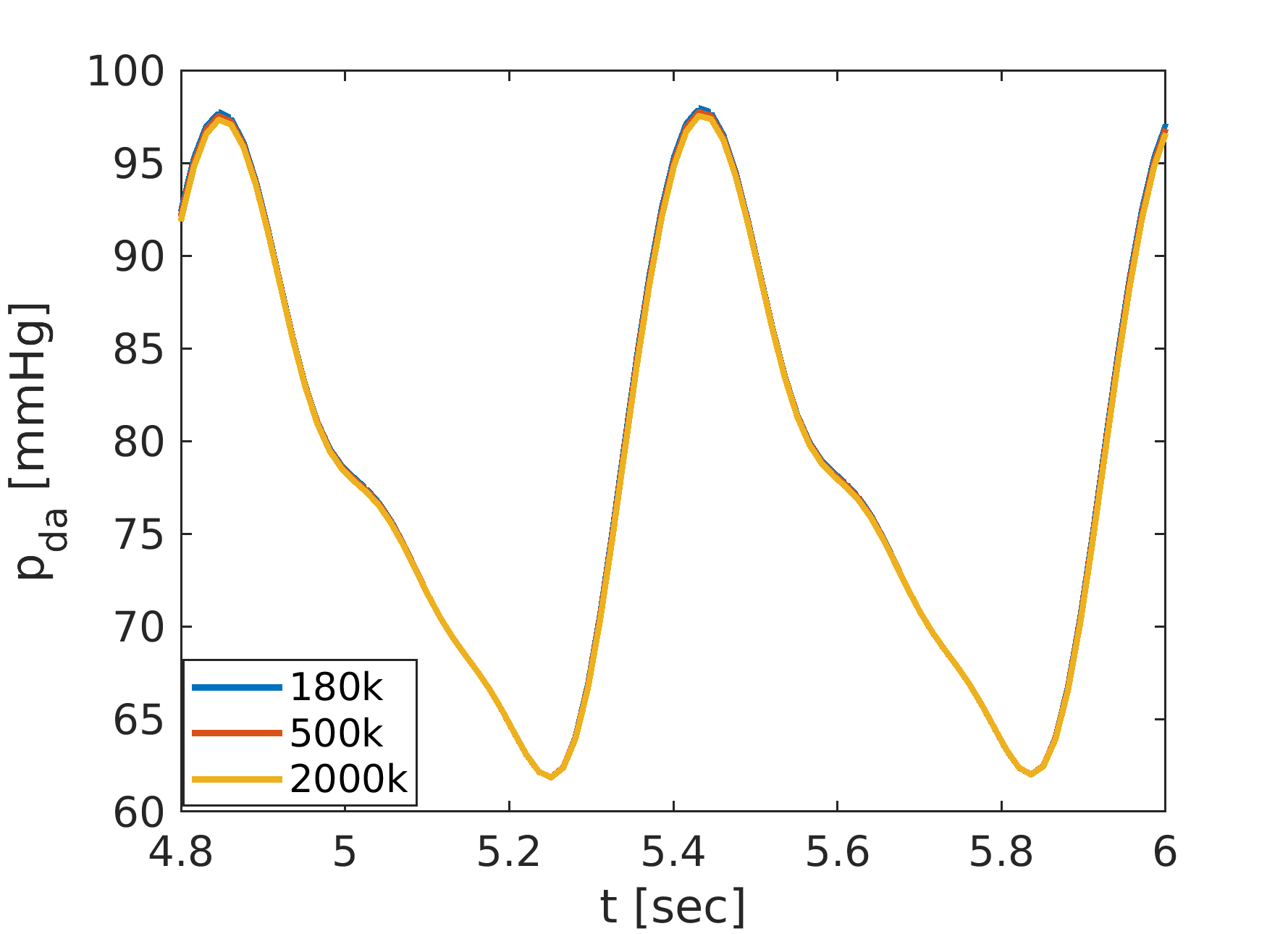}
\put(50,70){\small{b)}}
     \end{overpic}
\caption{Pre-surgery configuration: time evolution on two cardiac cycles of the volume averaged arterial pressure $p_{avg}$ as defined in \eqref{eq:pavg} (a) and the pressure related to the descending aorta cross-section $p_{da}$ (b) for the different meshes under consideration.}\label{fig:convergence}
\end{figure}

Fig. \ref{fig:convergence} compares the solution obtained with all the meshes reported in Table \ref{tab:mesh} both in terms of a global variable, the volume averaged arterial pressure, $p_{avg}$, defined as
\begin{equation}\label{eq:pavg}
p_{avg} = \dfrac{1}{\Omega} \int_\Omega p\ d\Omega,
\end{equation}
and in terms of a local variable, the descending aorta cross-section pressure, $p_{da}$. We let the simulations run till transient effects are passed, $t^\star \approx 8 \cdot T$. For a more quantitative comparison, we computed the Weighted Absolute Percentage Error (WAPE) $\varepsilon$ \cite{MAKRIDAKIS1993527} with respect to the solution obtained with the finer mesh $2000k$:
\begin{equation}\label{eq:WAPE}
\varepsilon = \dfrac{100}{n} \sum_{i=1}^{n} \left|\dfrac{X_i - X_i^{2000k}}{\overline{X^{2000k}}}\right| \%,
\end{equation}
where $n$ is the number of sampling points, $X_i$ is the solution related either meshes $230k$ and $415k$ at the $i$-th time step, $X_i^{2000k}$ is the solution related to the mesh $2000k$ at the $i$-th time step and $\overline{X^{2000k}}$ is the time-averaged solution related to the mesh $2000k$. For $p_{avg}$, we obtained $\varepsilon = 0.34 \%$ for the mesh $230k$ and $\varepsilon = 0.16 \%$ for the mesh $415k$. On the other hand, for $p_{da}$, we obtained $\varepsilon = 0.28 \%$ for the mesh $230k$ and $\varepsilon = 0.13 \%$ for the mesh $415k$. Thus, hereinafter, we will refer to the solutions computed by using the mesh $230k$. Regarding the post-surgery configuration, we choose a mesh with a similar refinement, having $200k$ cells, $h_{min} = 6.3e-4$ and $h_{max} = 3.4e-3$.

\subsubsection{Pre-surgery configuration}
The comparison between computational and experimental data is carried out in terms of systolic arterial pressure $PAS$, diastolic arterial pressure $PAD$ and average arterial pressure $PAM$. Computational estimates of such quantities are evaluated by simulations in the following way:
\begin{equation}\label{eq:PAS_c}
PAS = \max\limits_{t\in [0, T]} p_{avg},
\end{equation}

\begin{equation}\label{eq:PAD_c}
PAD = \min\limits_{t\in [0, T]} p_{avg},
\end{equation}

\begin{equation}\label{eq:PAM_c}
PAM = \dfrac{1}{T} \int_0^T p_{avg} dt.
\end{equation}


Fig. \ref{fig:convergence} (a) shows the temporal evolution of $p_{avg}$ (eq. \eqref{eq:pavg}). Table \ref{tab:comp_num_exp_pre} reports both numerical and experimental data marked by the abbreviations \emph{num} and \emph{exp}, respectively. We observe that the agreement is very good, within 11.7\% for $PAS$, 4\% for $PAD$ and, 2.4\% for $PAM$.

\begin{table}
\begin{center}
\begin{tabular}{|c|c|c|}
\hline
$PAS$ ($exp/num$) [mmHg] & $PAD$ ($exp/num$) [mmHg] & $PAM$ ($exp/num$) [mmHg]   \\\hline
108/95.4 & 66/63.4 & 78/79.9 \\\hline
\end{tabular}
\end{center}
\caption{Pre-surgery configuration: comparison between computational and experimental data.}\label{tab:comp_num_exp_pre}
\end{table}

\begin{figure}[h]
\centering
 \begin{overpic}[width=0.22\textwidth]{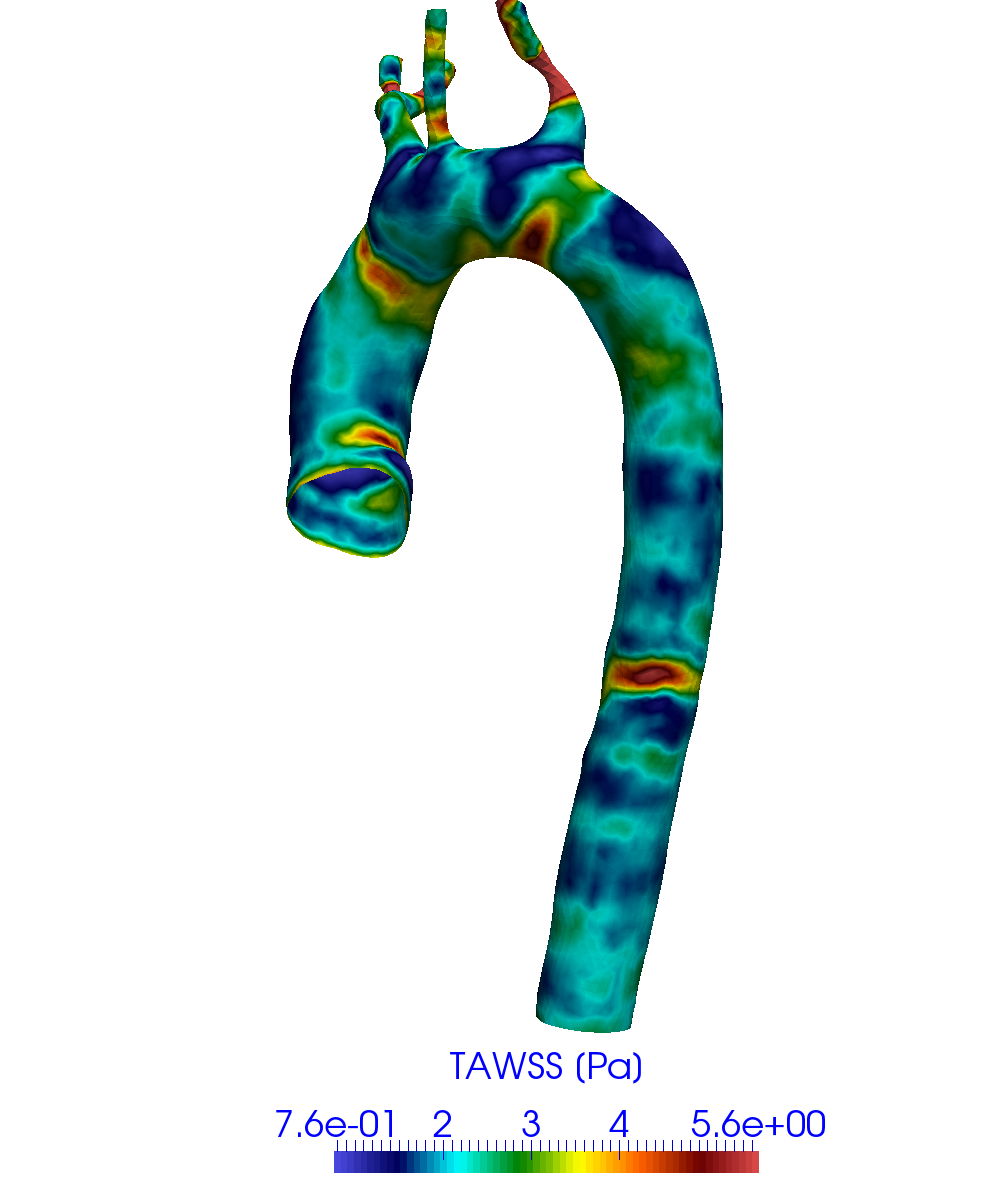}
      \end{overpic}
 \begin{overpic}[width=0.22\textwidth]{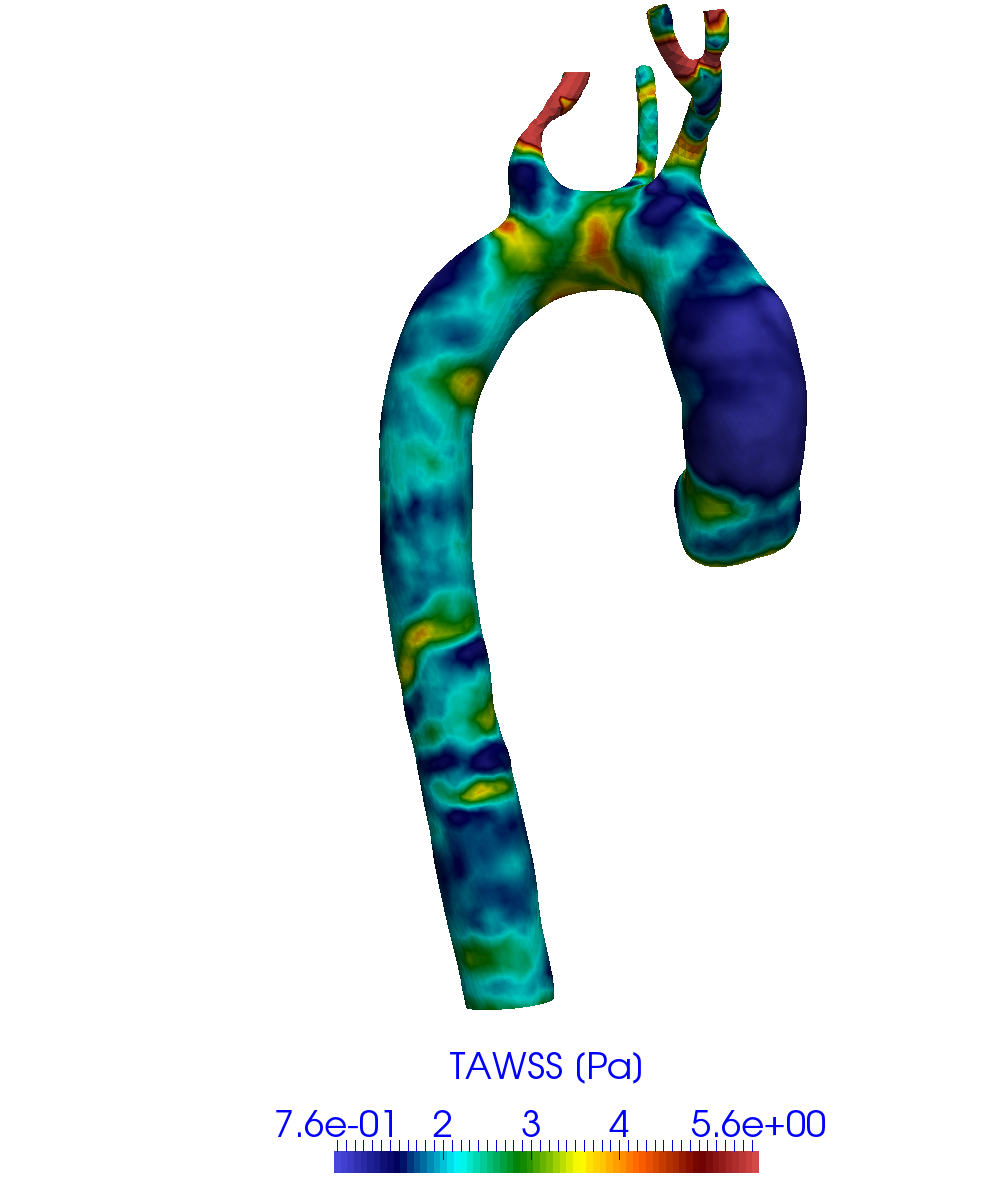}
     \end{overpic}
 \begin{overpic}[width=0.22\textwidth]{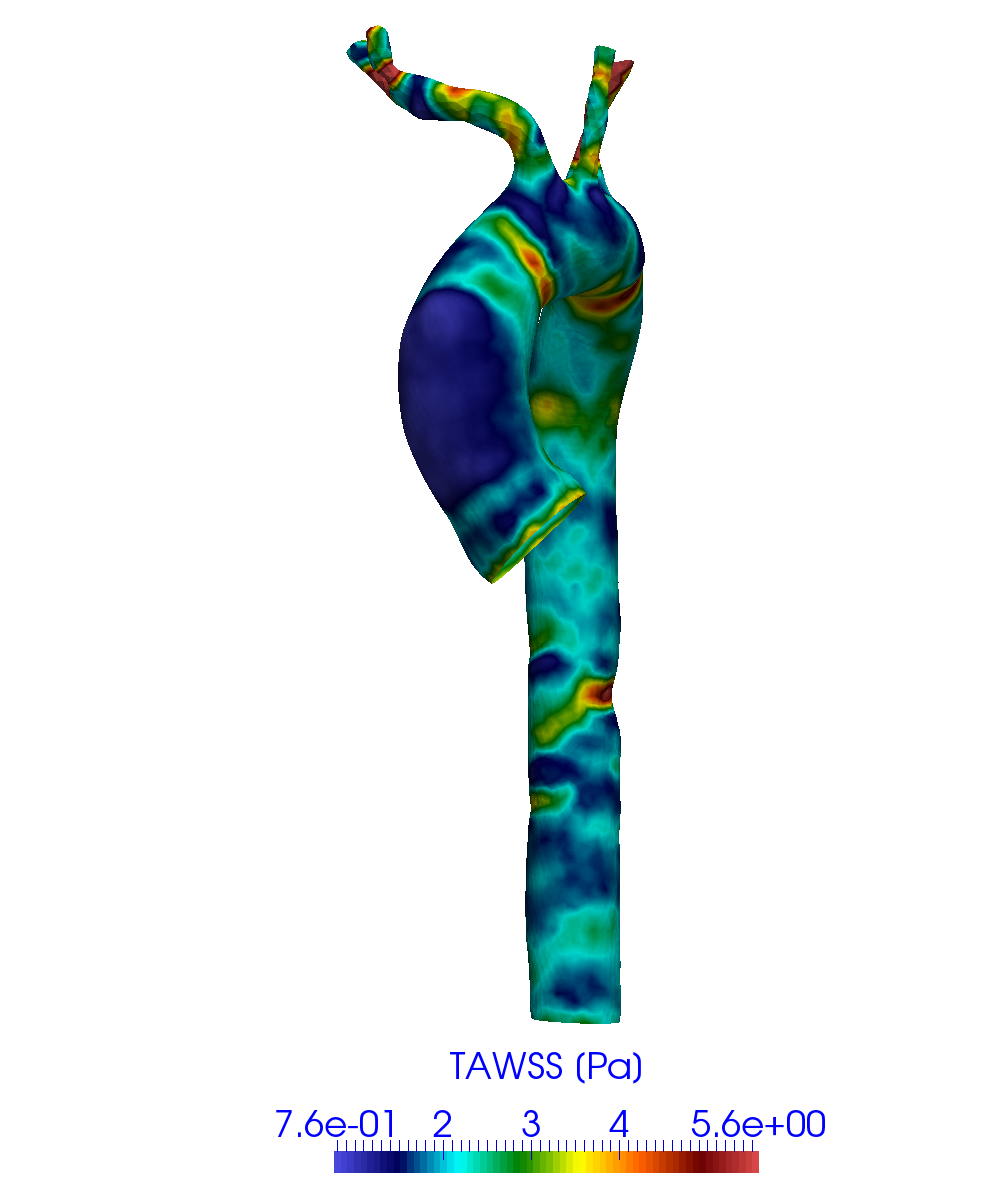}
     \end{overpic}
 \begin{overpic}[width=0.22\textwidth]{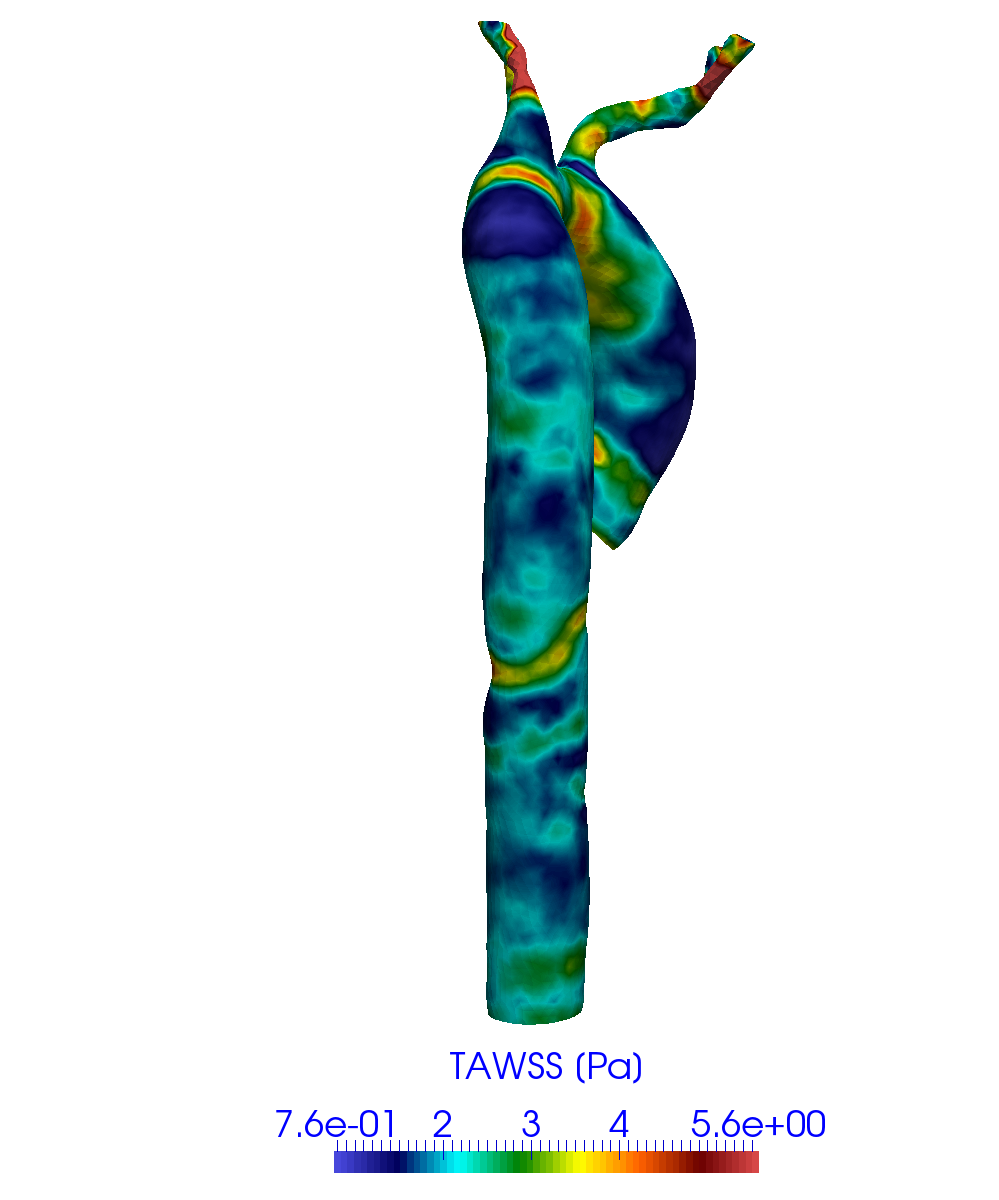}
     \end{overpic}
\caption{Pre-surgery configuration: TAWSS magnitude distribution on the entire wall of the model.}\label{fig:pre_TAWSS}
\end{figure}

Fig. \ref{fig:pre_TAWSS} displays the TAWSS magnitude distribution. Since in this case experimental data are not available, we just provide rough indications in order to justify the patterns obtained. Basically, we observe that peak values of TAWSS are localized in regions where narrowing of cross section happens or characterized by large curvature. On the other hand, regions characterized by lower TAWSS correspond to section enlargements. These results are expected by considering the classic findings for a straight cylindric vessel with steady Poiseuille flow. In this simplified case, $WSS \propto 1/d^3$, where $d$ is the pipe diameter. For biomedical experimental works that confirm such trend, the reader could see, e.g., \cite{Farag2018, Ooij2017}.


\begin{figure}[h]
\centering
 \begin{overpic}[width=0.4\textwidth]{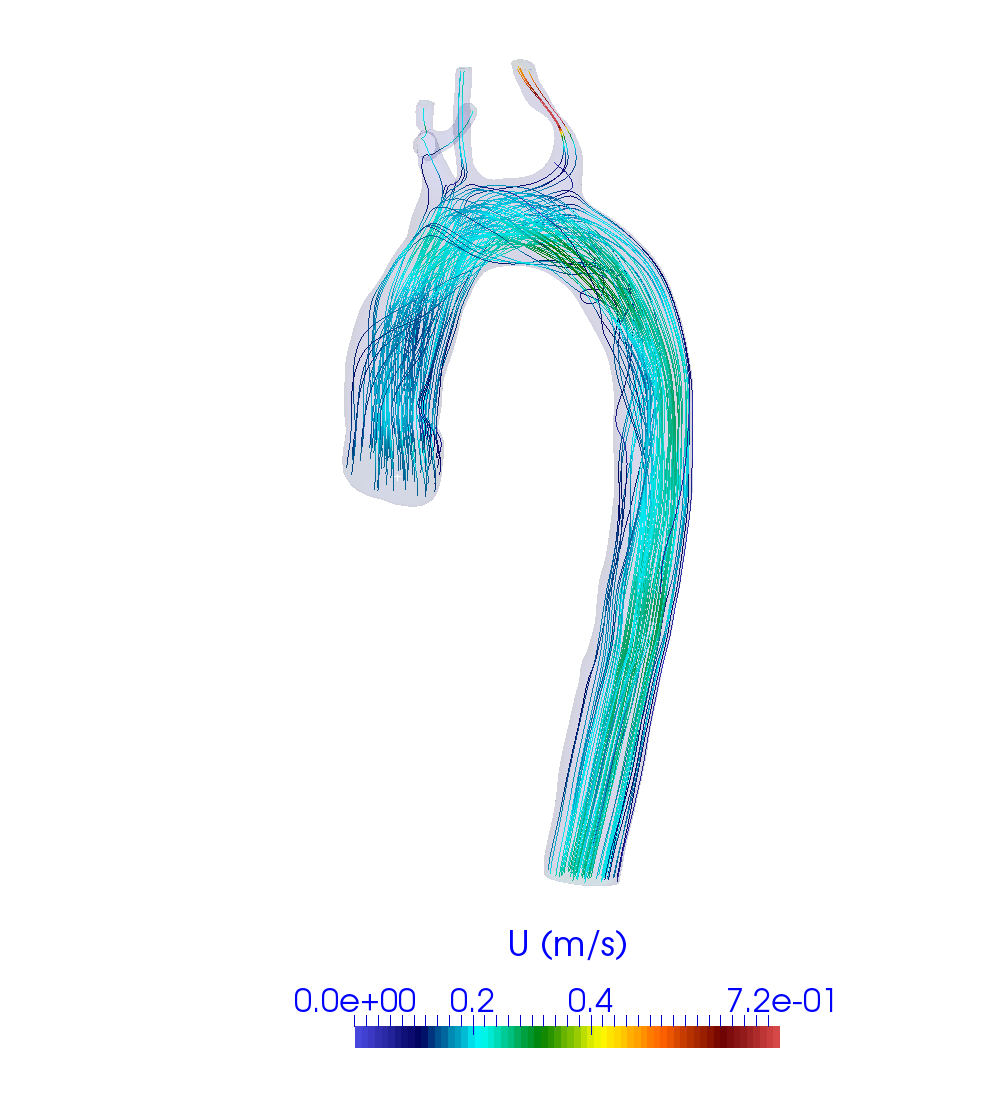}
\put(45,97){\small{a)}}
      \end{overpic}
 \begin{overpic}[width=0.38\textwidth]{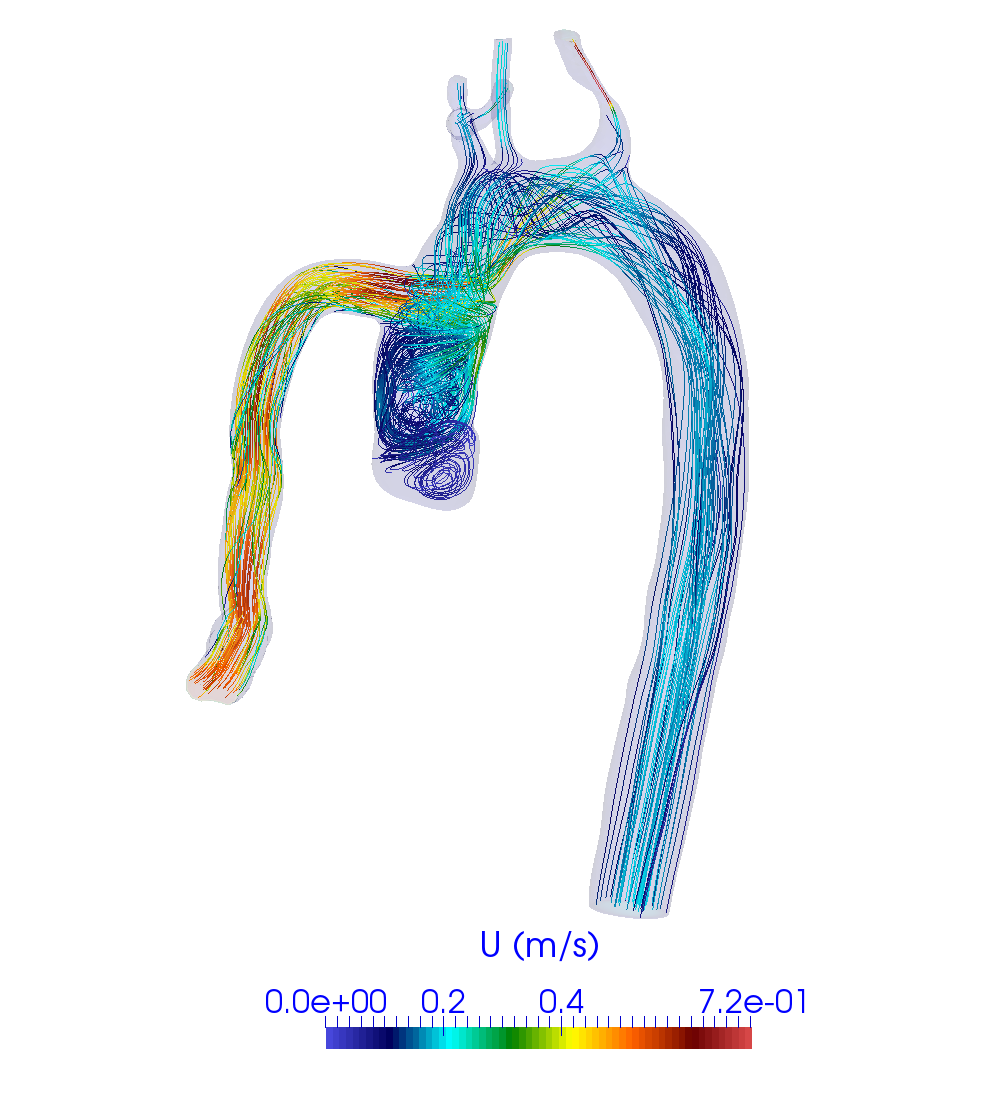}
\put(45,101){\small{b)}}
     \end{overpic}
\caption{Velocity streamlines rlated to the pre-surgery configuration (a) and the post-surgery configuration for $PF = 4.1$ l/min and $\omega = 5400$ rpm (b).}\label{fig:streamlines}
\end{figure}

Fig. \ref{fig:streamlines} a) depicts time averaged velocity streamlines. 
As expected, we note the generation of helical flow patterns in the aortic arch region (see, e.g. \cite{Lorenz2014}).

\subsubsection{Post-surgery configuration}
Unlike the pre-surgery case, in the post-surgery configuration, since the LVAD flow rate is continuous and not pulsatile, and the aortic valve is closed, the solution is steady in time. Therefore, the comparison between computational and experimental data is based on a value only, $PAM = PAD = PAS$. 

Table \ref{tab:comp_num_exp_post} reports both numerical (\emph{num}) and experimental (\emph{exp}) data for all the $PF$ values considered. We observe that the agreement is excellent, within 1\% in all the cases. 

\begin{table}
\begin{center}
\begin{tabular}{|c|c|}
\hline
$PF$ [l/min] & $PAM$ ($exp/num$) [mmHg]  \\\hline
4.1 &  78/78.6 \\\hline
4.2 & 90/90.6   \\\hline
4.5 & 100/100.7  \\\hline
5 &  83/82  \\\hline
\end{tabular}
\end{center}
\caption{Post-surgery configuration: comparison between computational and experimental data.}\label{tab:comp_num_exp_post}
\end{table}

\begin{figure}[h]
\centering
 \begin{overpic}[width=0.32\textwidth]{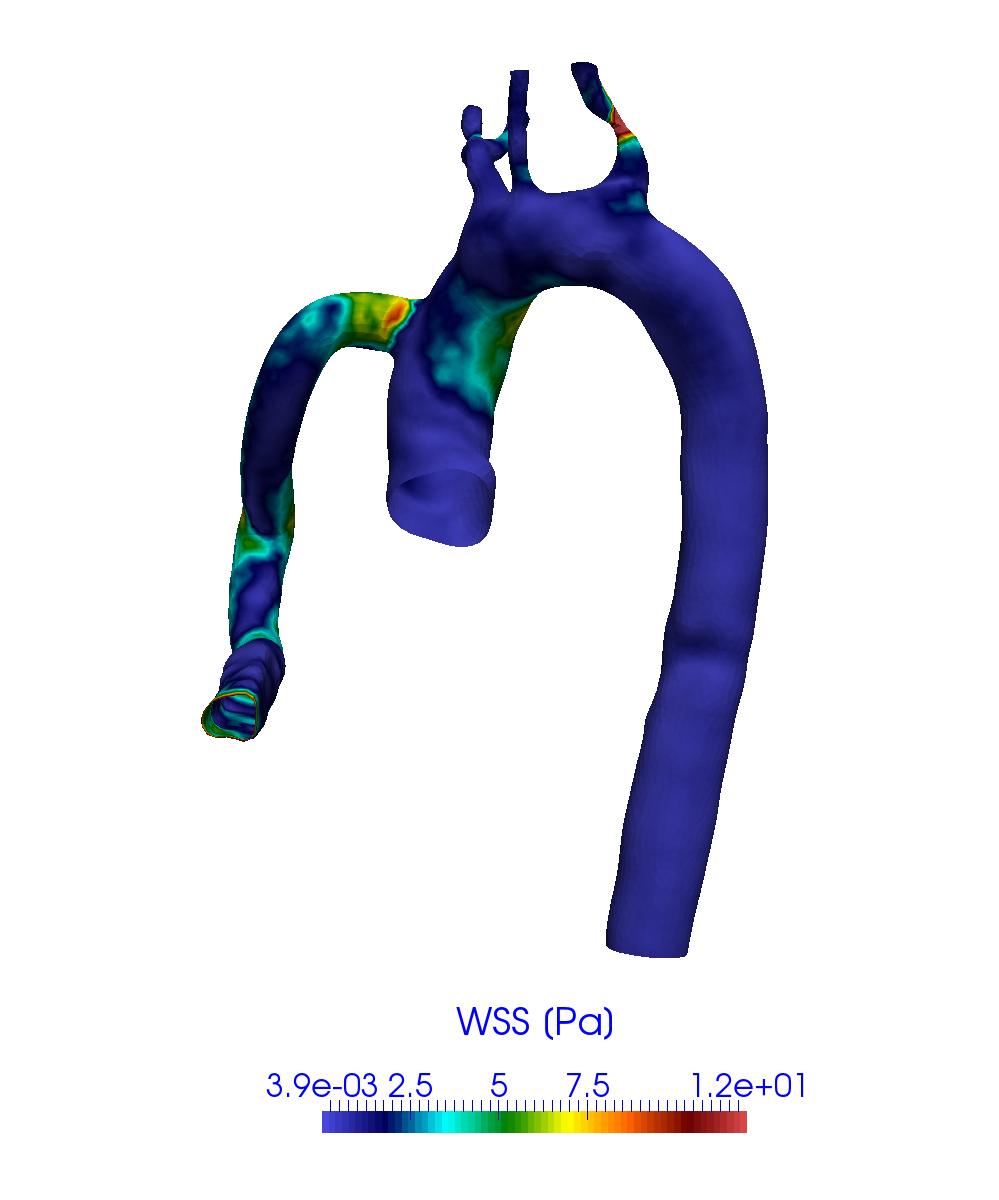}
      \end{overpic}
 \begin{overpic}[width=0.32\textwidth]{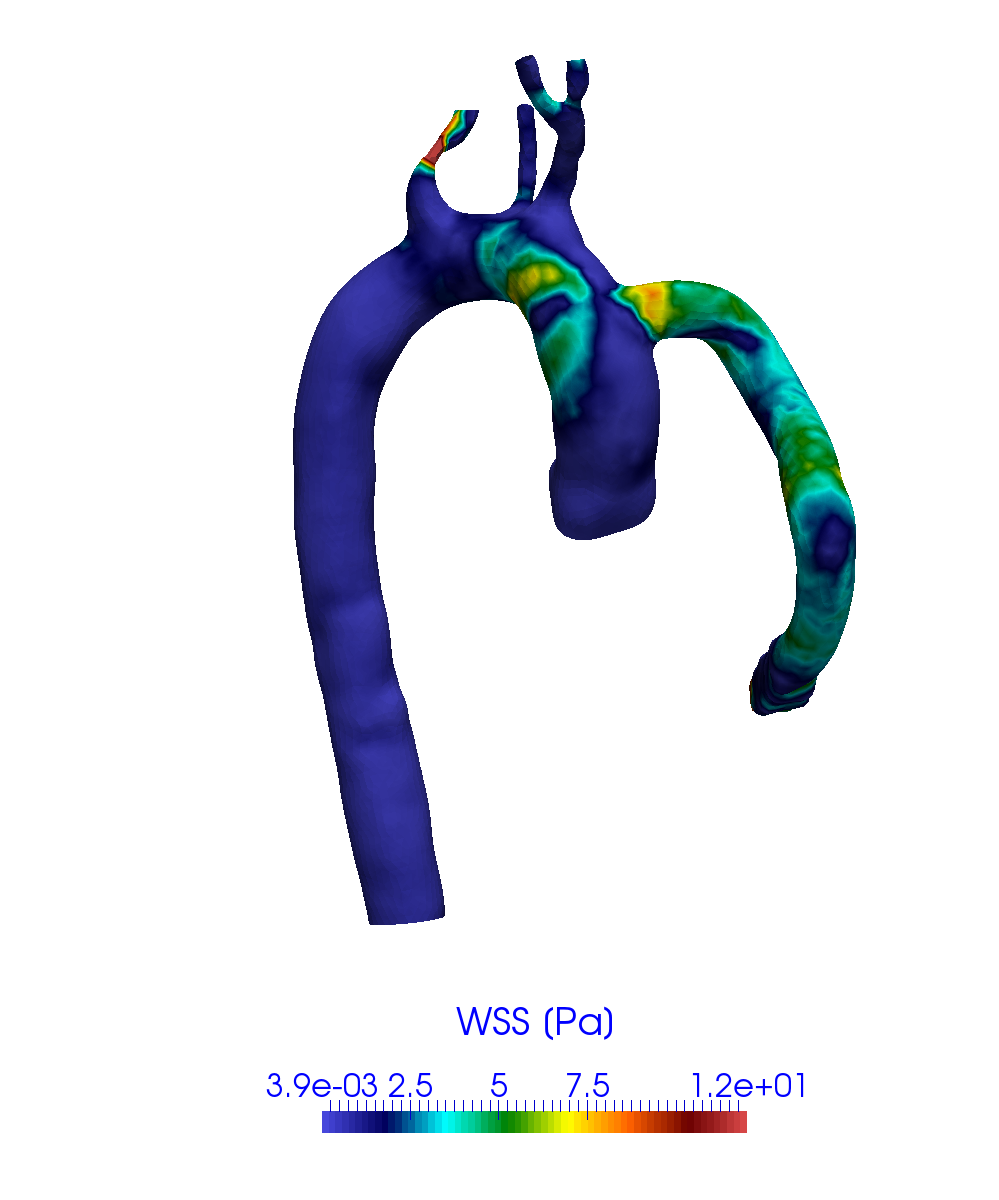}
     \end{overpic}
 \begin{overpic}[width=0.32\textwidth]{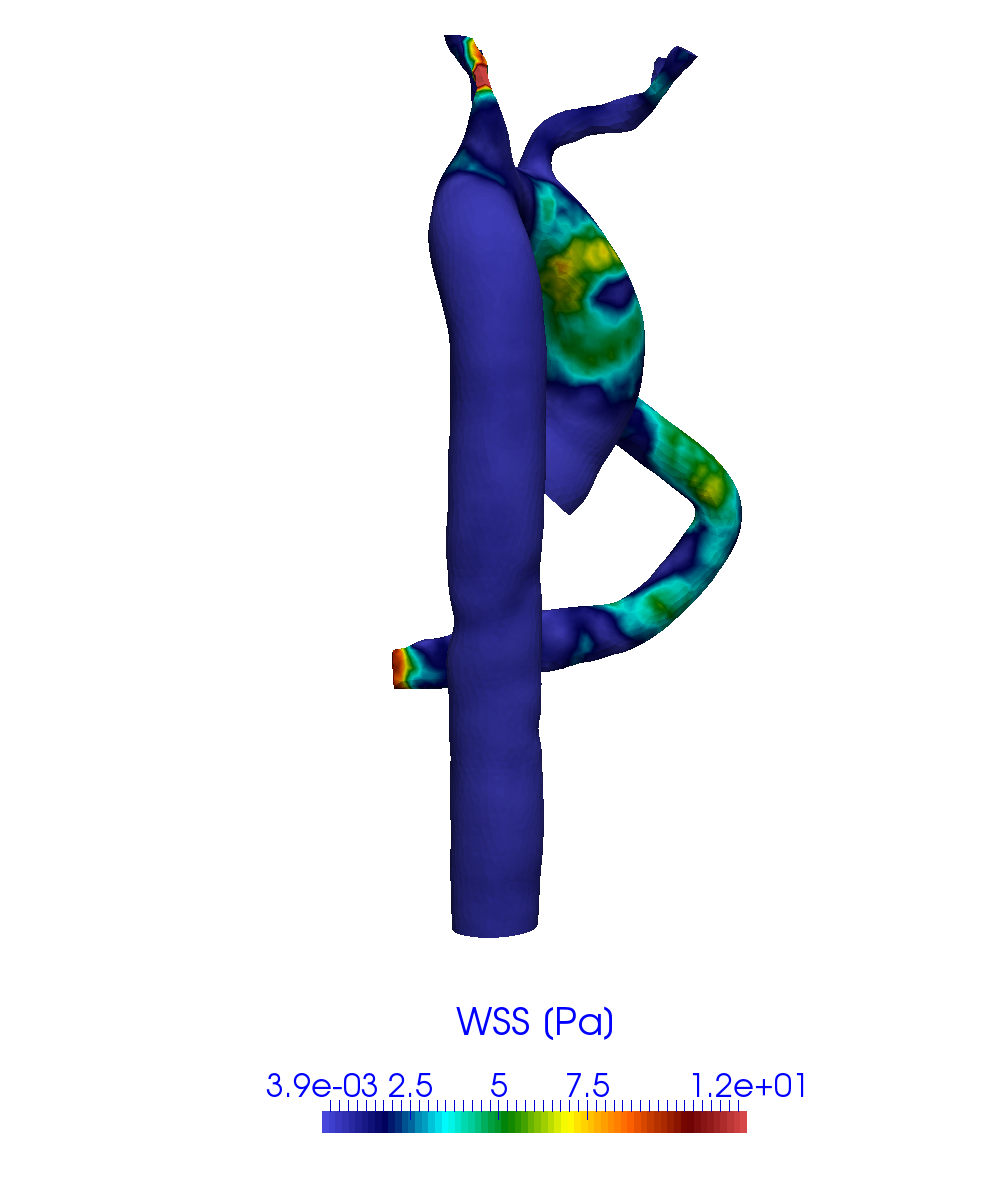}
      \end{overpic}
\caption{Post-surgery configuration: distribution of the WSS magnitude for $PF = 4.1$ l/min and $\omega = 5400$ rpm.}\label{fig:post1_WSS}
\end{figure}
\begin{figure}[h]
\centering
 \begin{overpic}[width=0.3\textwidth]{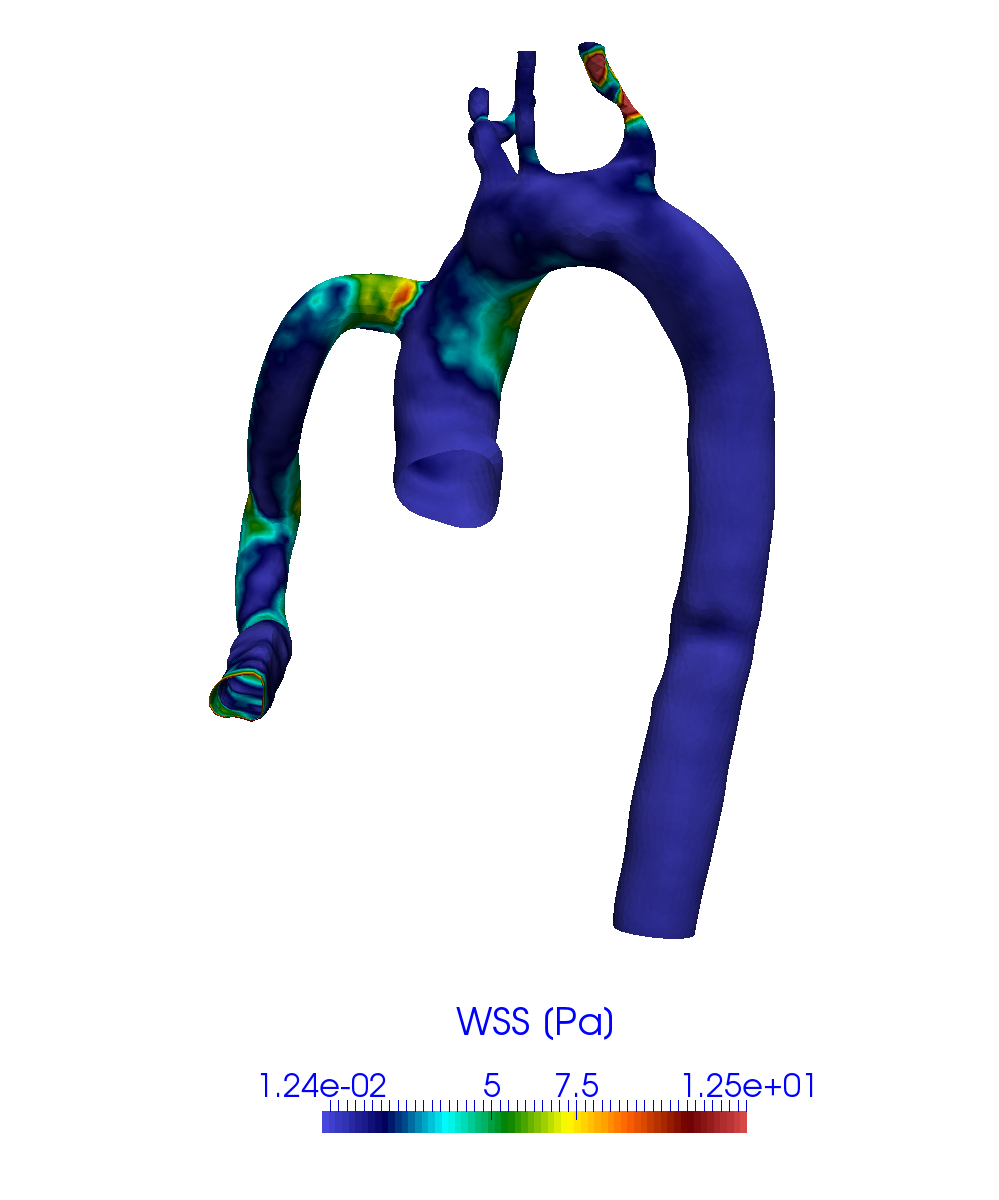}
      \end{overpic}
 \begin{overpic}[width=0.3\textwidth]{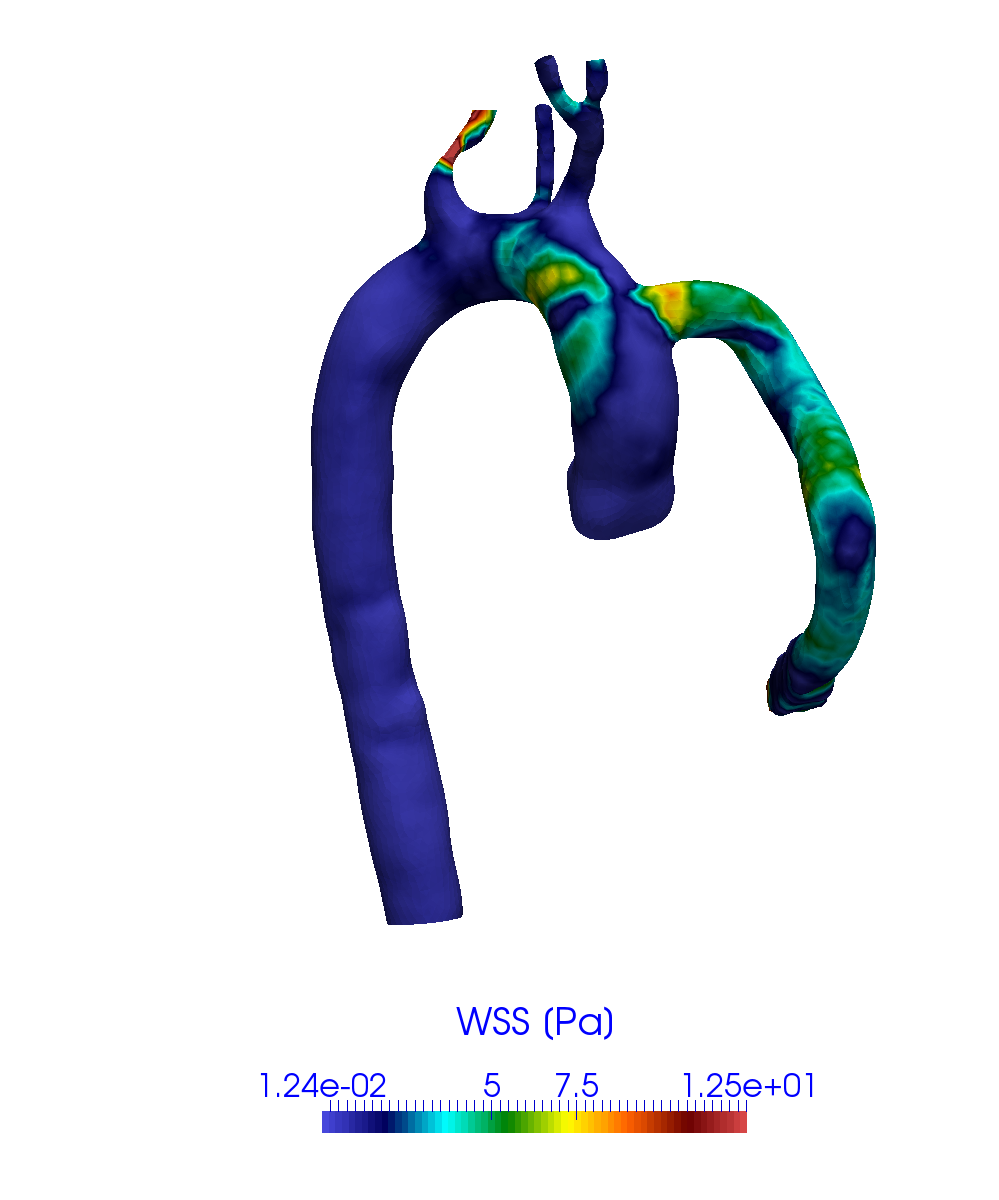}
     \end{overpic}
 \begin{overpic}[width=0.3\textwidth]{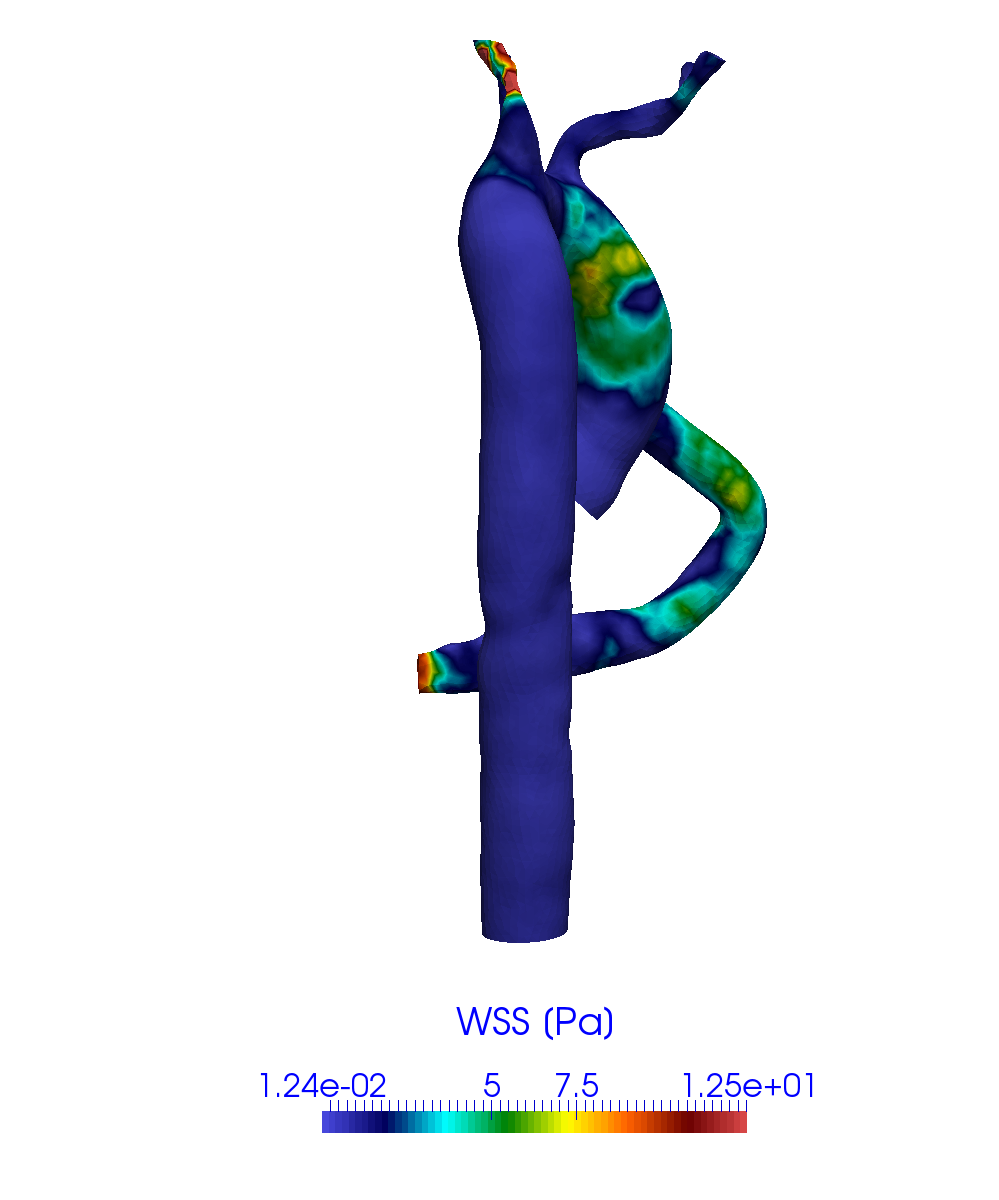}
      \end{overpic}
\caption{Post-surgery configuration: distribution of the WSS magnitude for $PF = 4.2$ l/min and $\omega = 5600$ rpm.}\label{fig:post2_WSS}
\end{figure}
\begin{figure}[h]
\centering
 \begin{overpic}[width=0.3\textwidth]{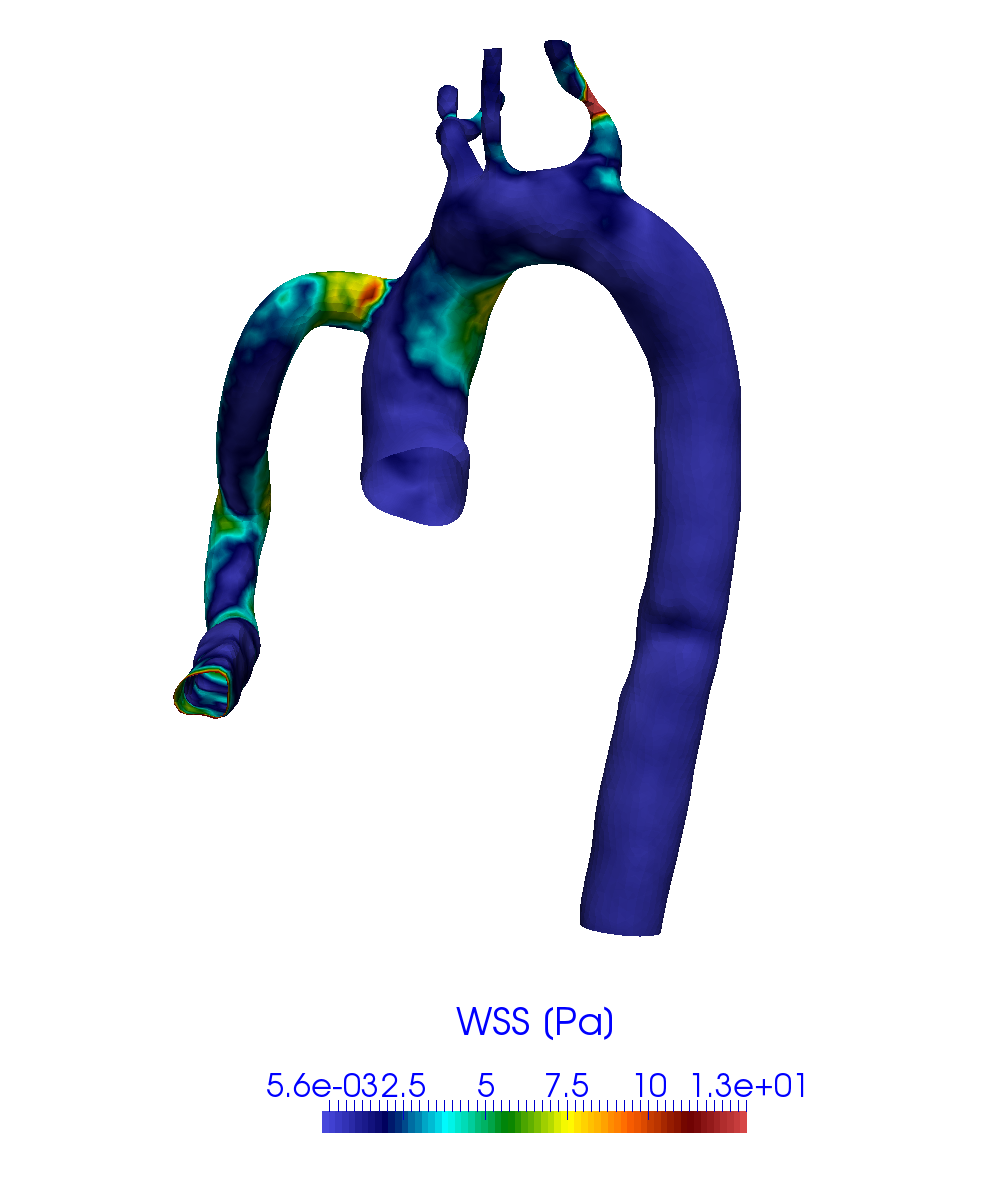}
      \end{overpic}
 \begin{overpic}[width=0.3\textwidth]{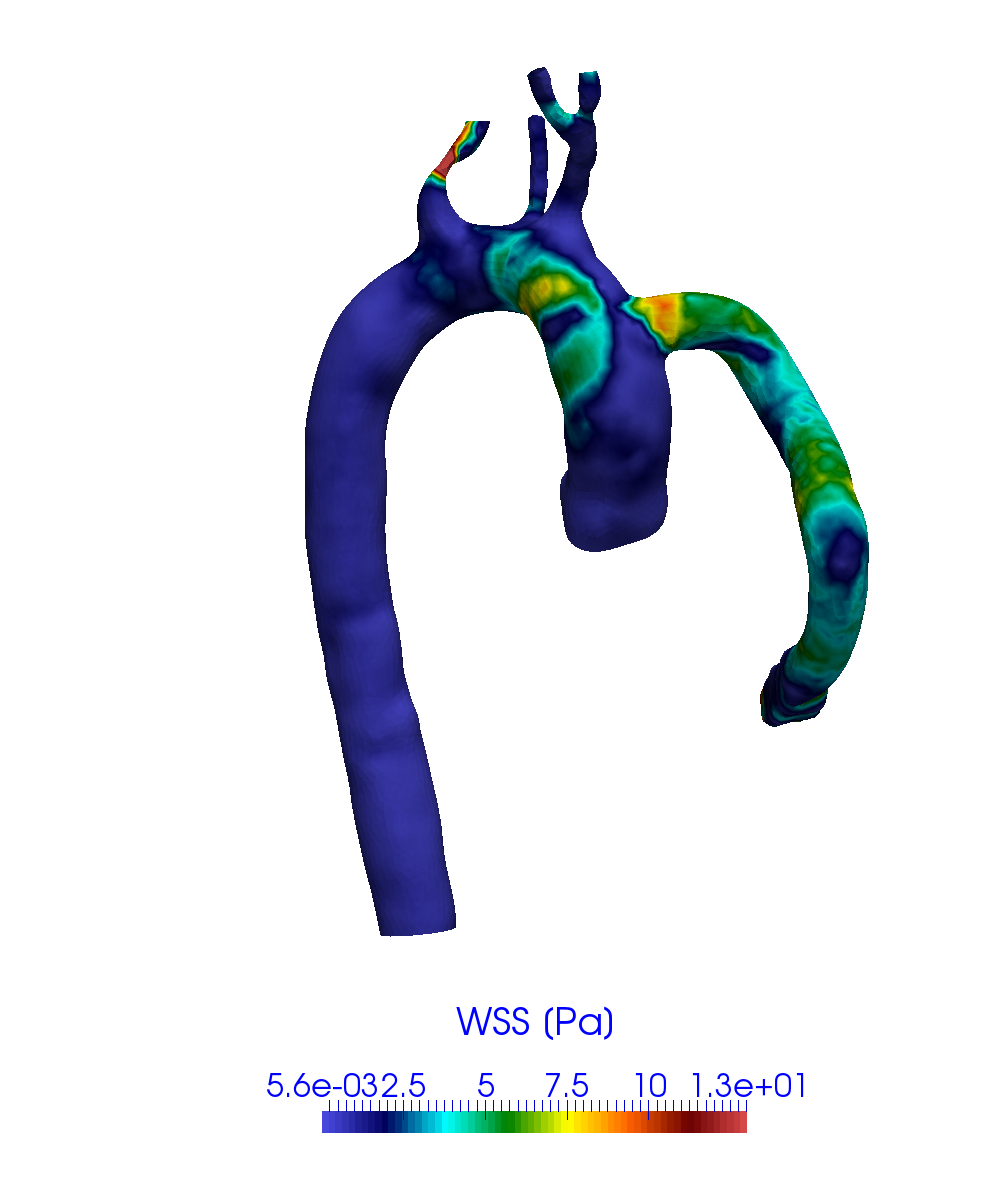}
     \end{overpic}
 \begin{overpic}[width=0.3\textwidth]{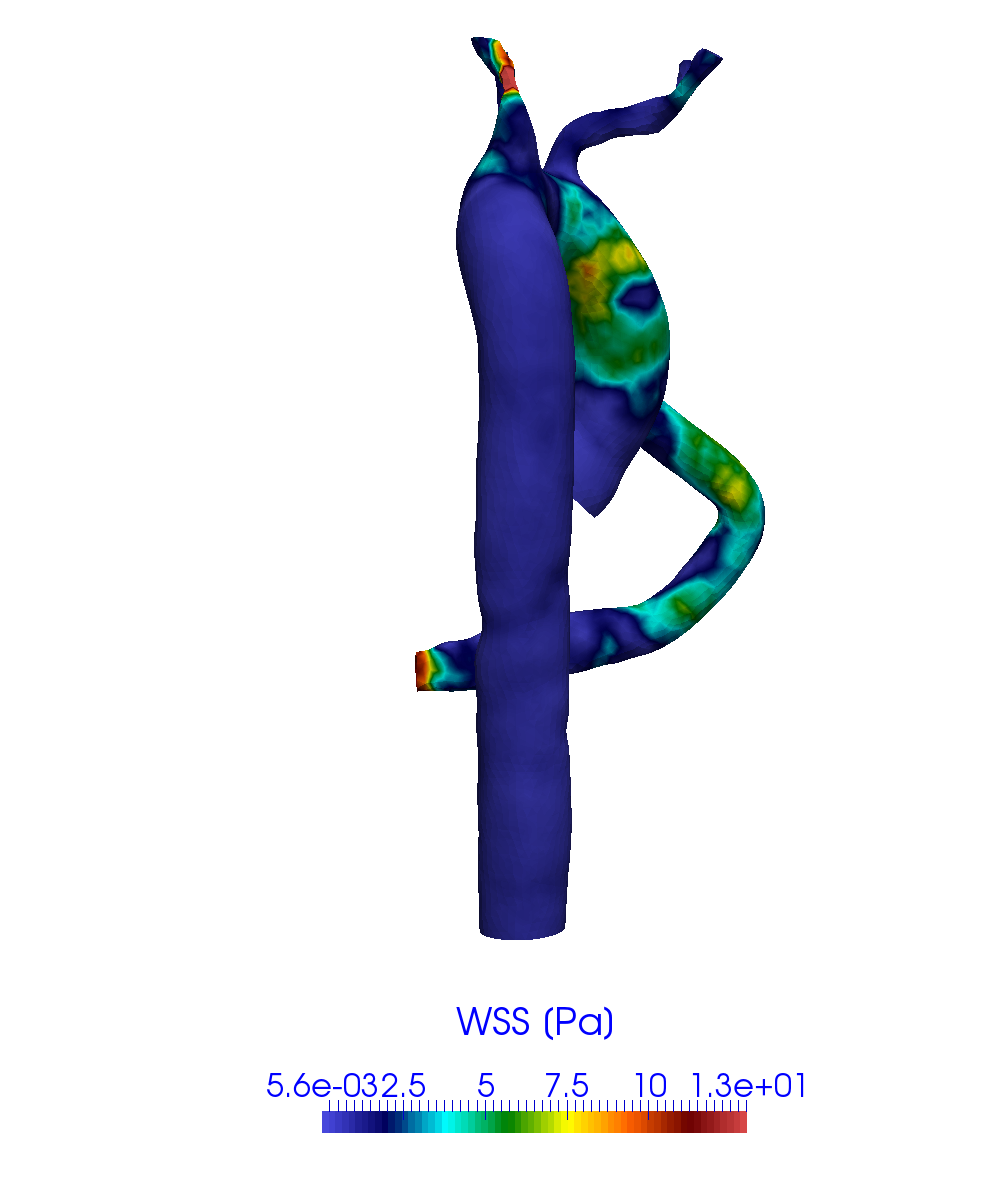}
      \end{overpic}
\caption{Post-surgery configuration: distribution of the WSS magnitude for $PF = 4.5$ l/min and $\omega = 6000$ rpm.}\label{fig:post3_WSS}
\end{figure}
\begin{figure}[h]
\centering
 \begin{overpic}[width=0.3\textwidth]{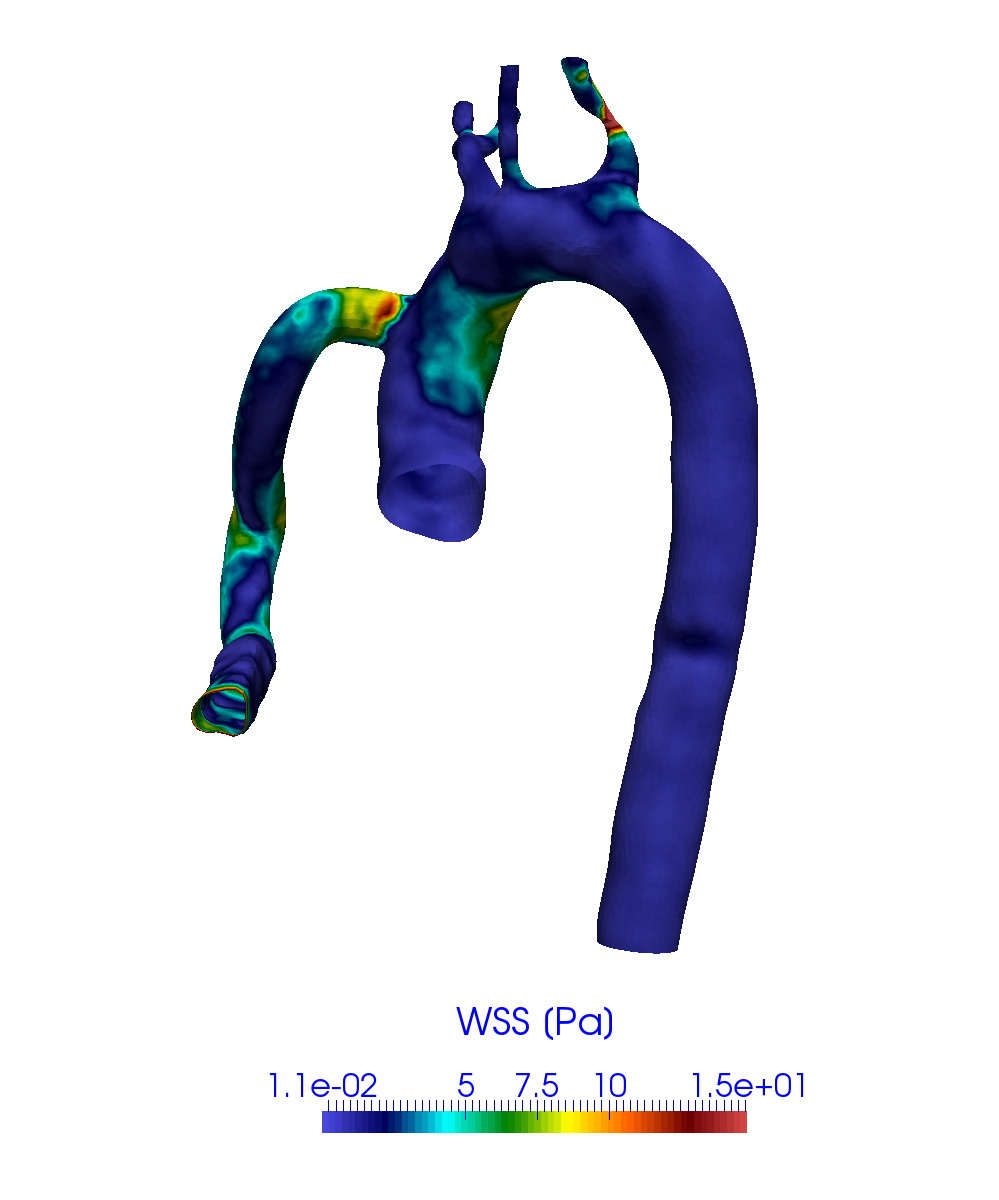}
      \end{overpic}
 \begin{overpic}[width=0.3\textwidth]{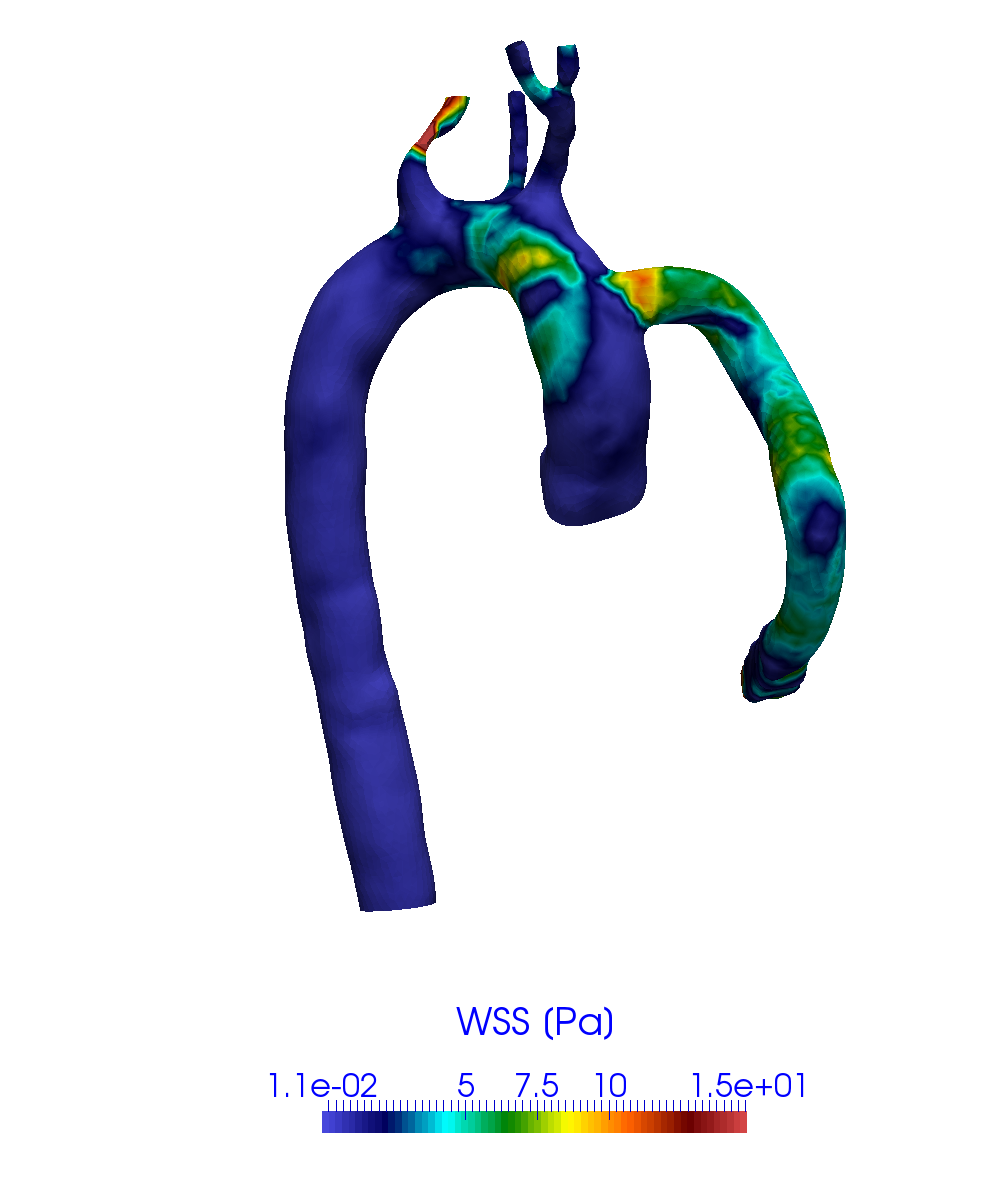}
     \end{overpic}
 \begin{overpic}[width=0.3\textwidth]{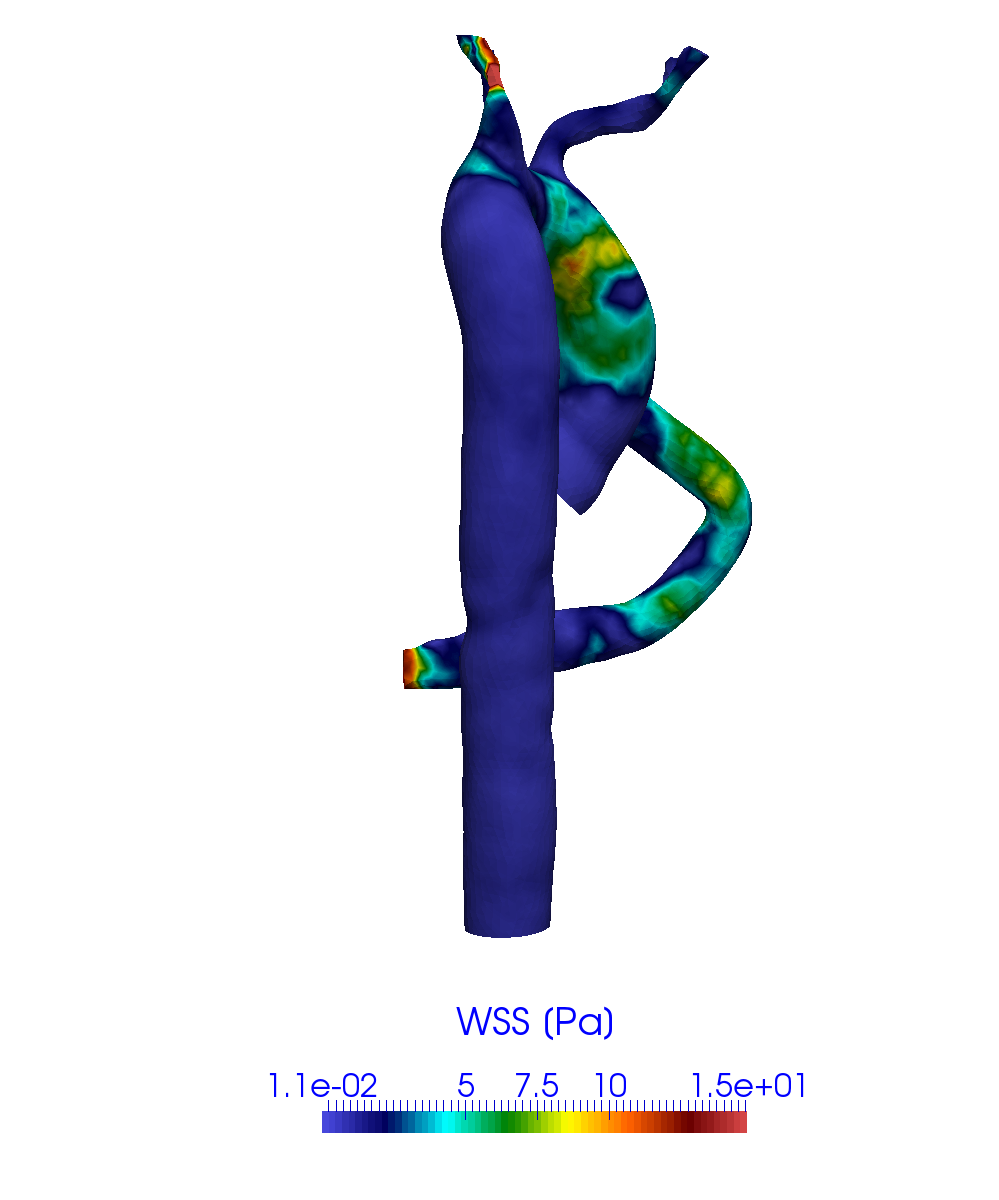}
      \end{overpic}
\caption{Post-surgery configuration: distribution of the WSS magnitude for $PF = 5$ l/min and $\omega = 5600$ rpm.}\label{fig:post4_WSS}
\end{figure}

Figs. \ref{fig:post1_WSS}-\ref{fig:post4_WSS} show the WSS distribution for the three configurations investigated. As for the pre-surgery configuration, even in this case experimental data related to WSS are not available but it is possible to provide some interesting observations to be compared with previous works. In all the cases, we observe that there are high WSS, significantly greater than that obtained in the pre-surgery configuration, on the posterior region of the aortic arch, in front of the anastomosis. This high WSS zone is associated with the impingement of the jet from the cannula. This result is in agreement with those observed by \cite{Aliseda2017, Caruso2015, Zhang2018}. Moreover we observe that elevated WSS also occur near the location of the outflow cannula, as found by \cite{Karmonik2012, Karmonik2014, Karmonik2014b, May-Newman2006, Inci2012, Brown2011, Zhang2018}.  On the contrary, on the most part of the aortic arch and descending aorta, very low WSS occurs. These patterns are critical from clinical viwpoints becase highly heterogeneous WSS distribution coupled with the presence of a small region of the aortic arch exposed to high WSS could be associated to the development of atherosclerosis \cite{Dolan2011, Malek1999}. Finally, we note that at increasing of $PF$ from $4.1$ to $5$ l/min, the peak value of WSS moves from $12$ to $15$ Pa by following an almost linear trend. 

Fig. \ref{fig:streamlines} b) displays the velocity streamlines for the Test 1. With respect to the pre-surgery configuration, we observe that in the ascending aorta, below the anastomosis location, retrograde
flow and recirculation zone are generated \cite{Caruso2015, Bonnemain2012, Karmonik2012, May-Newman2006, Aliseda2017}. In addition, we observe that velocity values in the outflow cannula are higher than those in aorta because of its smaller diameter. 

\subsection{ROM}\label{sec:res_ROM}
To train the ROM, the values of the LVAD flow rate, $PF$, are chosen using an equispaced distribution inside the range $PF \in [3, 5]$ that covers typical clinical values. Two sampling cases were considered. In the first case, we have 21 snapshots, and in the second one, 11 snapshots. Thus, the snapshots are collected every 0.1 in the first case and 0.2 in the latter one. For all the simulations, we use resistances and capacitances of Test 1 (see Tables \ref{tab:data_computed_pre}, \ref{tab:data_computed_post}, \ref{tab:bc_post}). By assuming that we vary $PF$, and consequently $\omega$, at a given $\Delta P = 75$ (see Table \ref{tab:delta_P}) and using the analytical fit \ref{eq:pump}, we obtain that the range $PF \in [3, 5]$ corresponds to $\omega \in [5076, 5720]$. Two new values of the $PF$ ($\omega$) in which the ROM has not been trained but which belongs to the range of the training space, $PF = 3.45$ ($\omega = 5200$) and $PF = 4.35$ ($\omega = 5484$), are used to evaluate the performance of the parametrized ROM. POD modes and coefficients are computed as explained in Section \ref{sec:ROM}.

\begin{figure}[h]
\centering
 \begin{overpic}[width=0.45\textwidth]{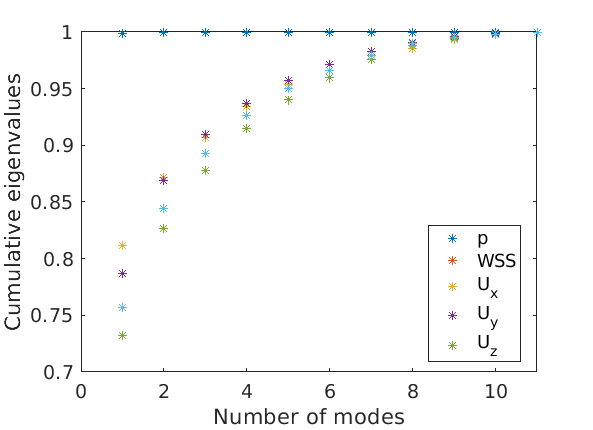}
\put(50,70){\small{a)}}
      \end{overpic}
 \begin{overpic}[width=0.45\textwidth]{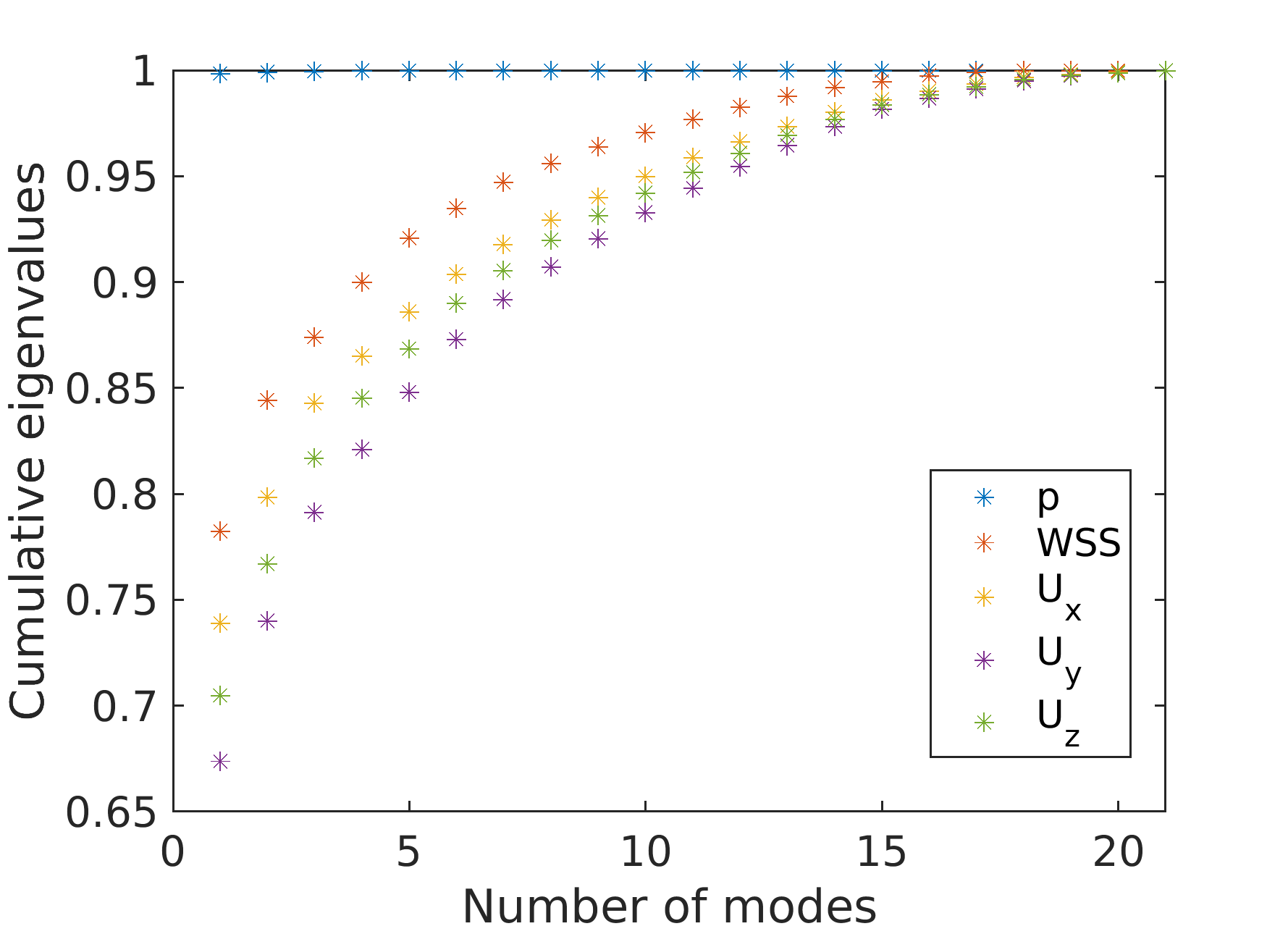}
\put(50,70){\small{b)}}
     \end{overpic}
\caption{Cumulative energy of the eigenvalues for pressure $p$, wall shear stress WSS, and velocity components, $u_x$, $u_y$ and $u_z$. The sampling frequency of the eigenvalues is 0.2 (a) and 0.1 (b).}\label{fig:cum_eigen}
\end{figure}

Figure \ref{fig:cum_eigen} shows the cumulative energy of the eigenvalues for pressure $p$, wall shear stress WSS, and velocity components, $u_x$, $u_y$ and $u_z$.
In order to retain the 99.9\% of the system's energy, when we consider 21 snapshots, 1 mode for $p$, 14 for WSS, 16 for $u_x$,  $u_y$ and $u_z$ are selected. On the other hand, when 11 snapshots are taken into account, 1 mode for $p$, 9 for WSS, 8 for $u_x$, $u_y$ and $u_z$ are selected. Moreover, to provide some quantitative results, the relative error in the $L^2$-norm, calculated as
\begin{equation}\label{eq:error1}
E_X = 100\dfrac{||X_{FOM} - X_{ROM}||_{L^2(\Omega)}}{||X_{FOM}||_{L^2(\Omega)}} \%
\end{equation}
where $X_{FOM}$ is the value of a particular field in the FOM model, and $X_{ROM}$ the one that is calculated using the ROM, is considered. In Tables \ref{tab:relative_errors1} and \ref{tab:relative_errors2}, one could observe that the differences between the two spaces are minimal, for both the values of $PF$ considered. Therefore, hereinafter results will be based on the database of 11 snapshots. 

\begin{table}
\begin{center}
\begin{tabular}{|c|c|c|c|c|c|}
\hline
 & $p$ & WSS &  $u_x$ & $u_y$ & $u_z$ \\\hline
$E_X$ (21 snapshosts) &  0.1\% & 4.1\% & 5.6\% & 7.9\% & 6.2\% \\\hline
$E_X$ (11 snapshosts) &  0.2\% & 4.1\% & 5\% & 7.8\% & 5.8\% \\\hline
\end{tabular}
\end{center}
\caption{$L^2$ norm relative errors for pressure $p$, wall shear stress WSS, and velocity components, $u_x$, $u_y$ and $u_z$, to varying of the number of snapshots collected for $PF = 3.45$ l/min.}\label{tab:relative_errors1}
\end{table}

\begin{table}
\begin{center}
\begin{tabular}{|c|c|c|c|c|c|}
\hline
 & $p$ & WSS &  $u_x$ & $u_y$ & $u_z$ \\\hline
$E_X$ (21 snapshosts) &  0.2\% & 9.6\% & 10.7\% & 14.5\% & 10.5\% \\\hline
$E_X$ (11 snapshosts) &  0.5\% & 7.2\% & 9.7\% & 13.5\% & 9.3\% \\\hline
\end{tabular}
\end{center}
\caption{$L^2$ norm relative errors for pressure $p$, wall shear stress WSS, and velocity components, $u_x$, $u_y$ and $u_z$, to varying of the number of snapshots collected for $PF = 4.35$ l/min.}\label{tab:relative_errors2}
\end{table}

\begin{figure}[h]
\centering
 \begin{overpic}[width=0.3\textwidth]{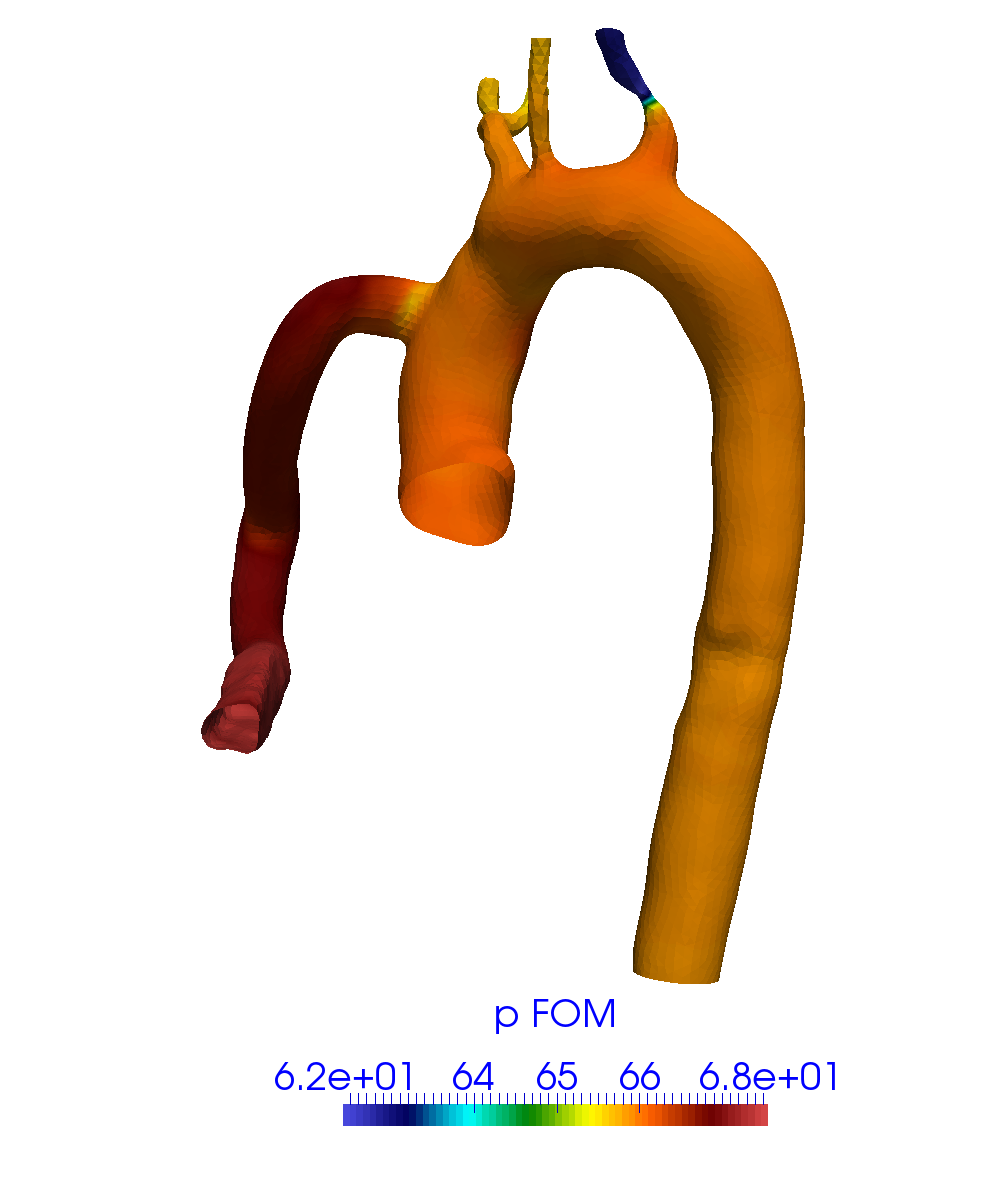}
      \end{overpic}
 \begin{overpic}[width=0.3\textwidth]{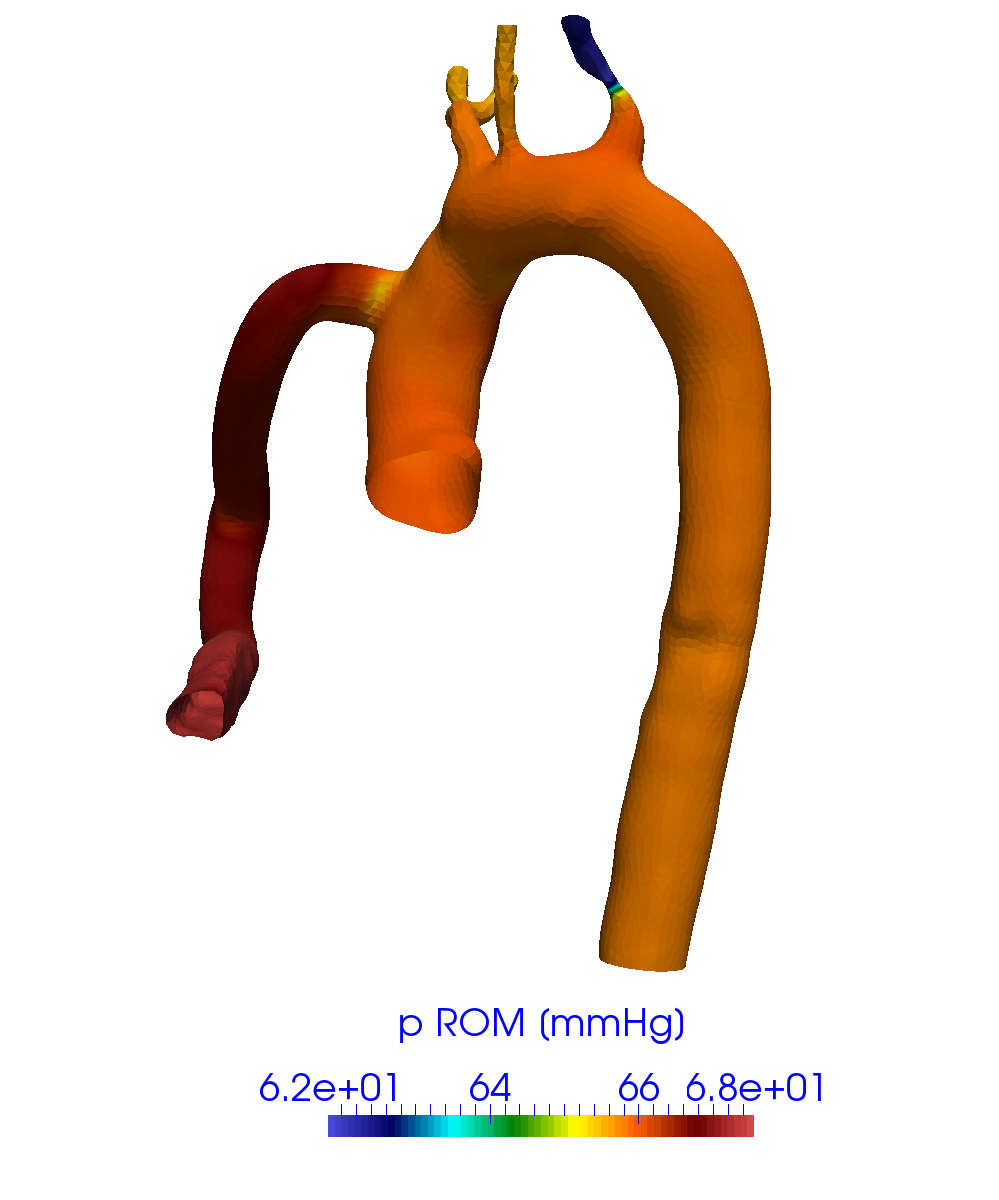}
     \end{overpic}
 \begin{overpic}[width=0.3\textwidth]{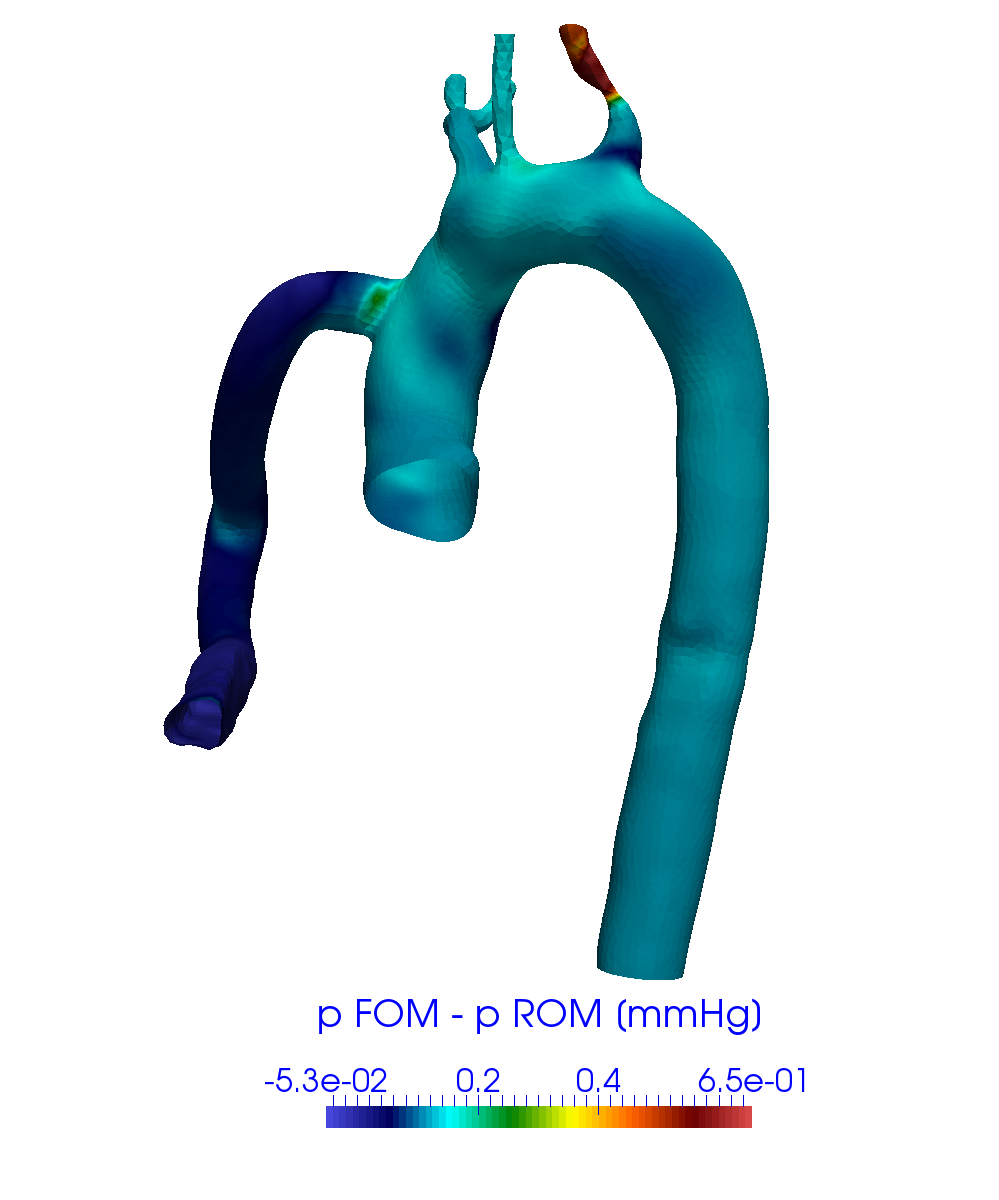}
      \end{overpic}
 \begin{overpic}[width=0.3\textwidth]{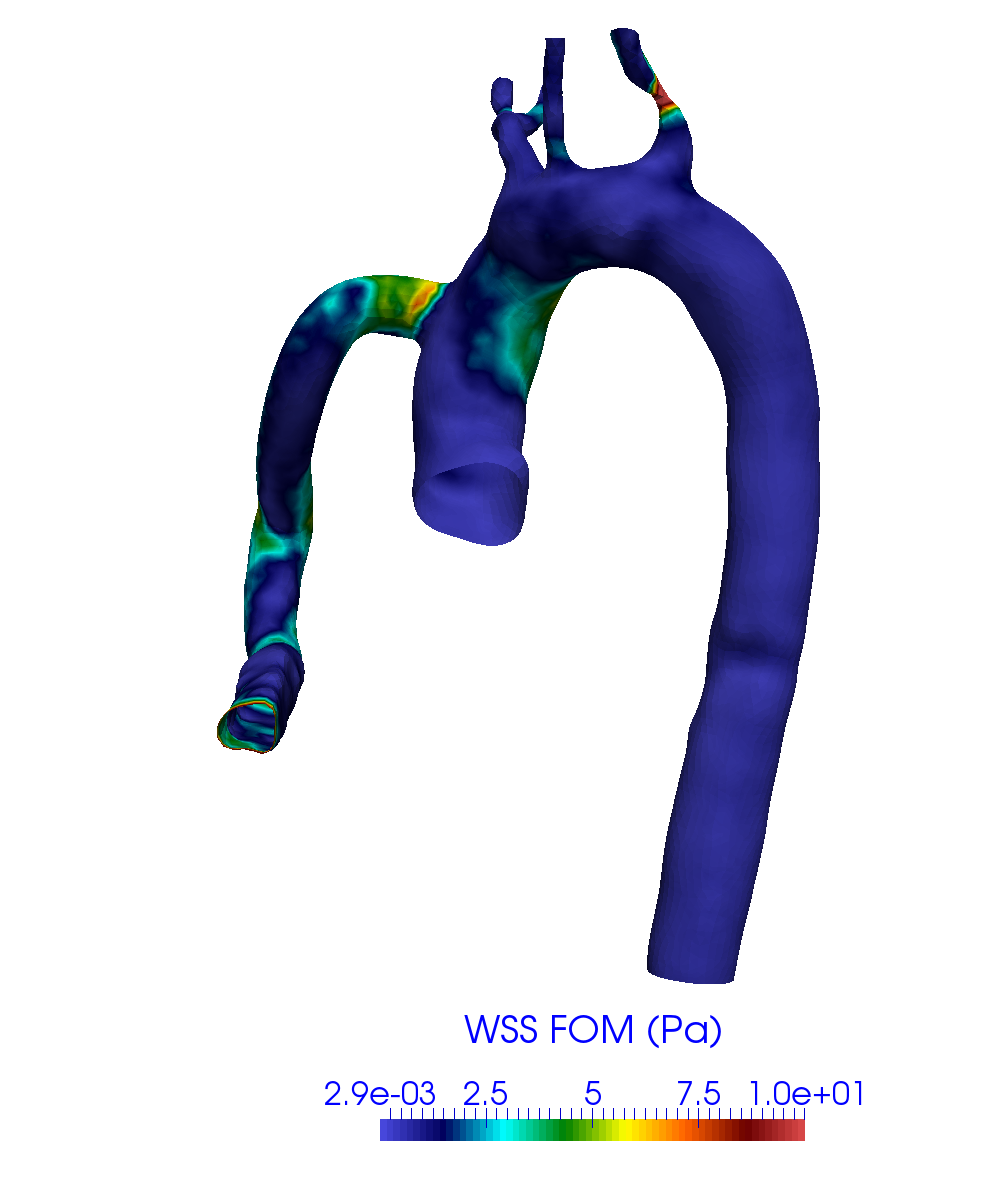}
      \end{overpic}
 \begin{overpic}[width=0.3\textwidth]{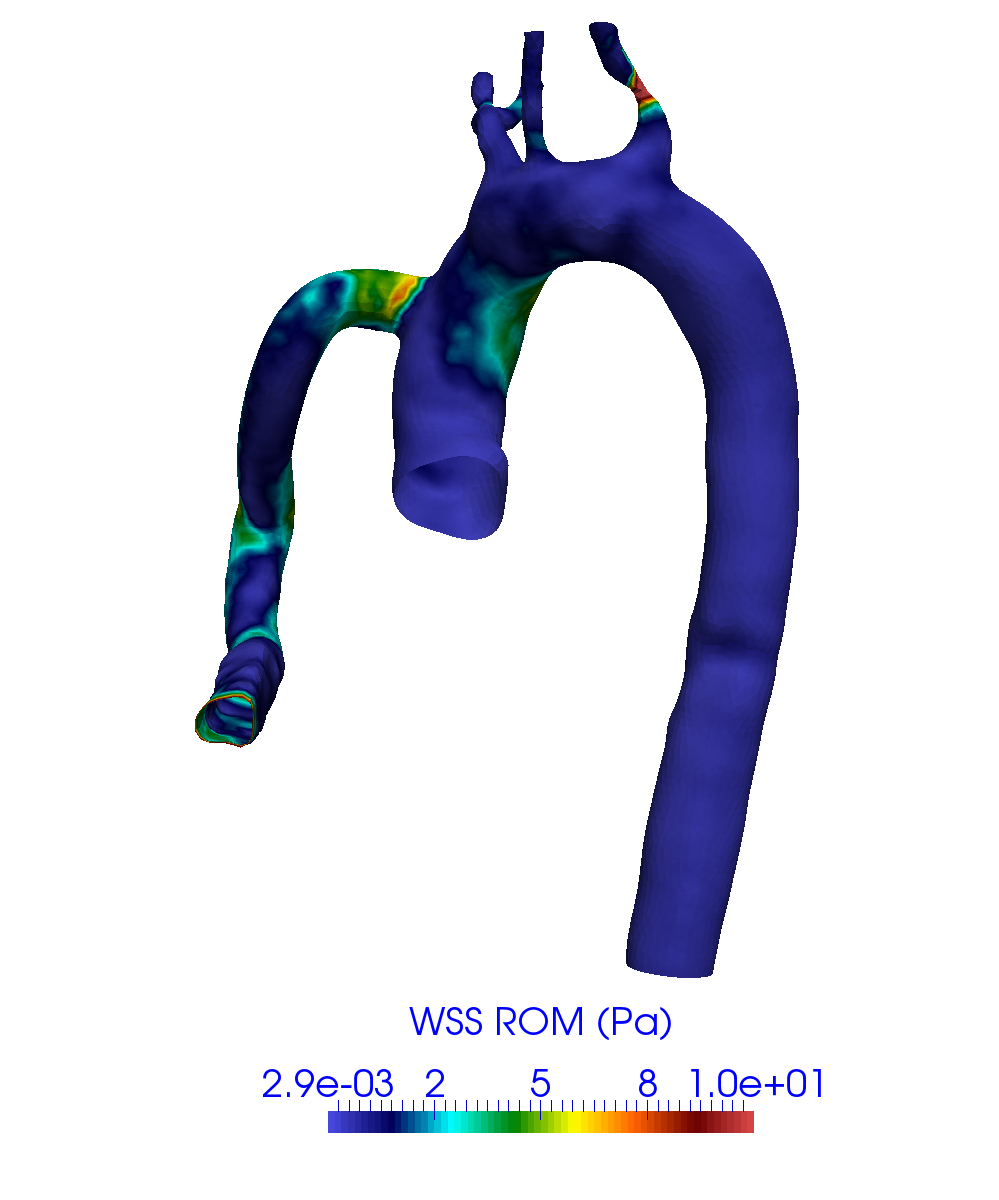}
     \end{overpic}
 \begin{overpic}[width=0.3\textwidth]{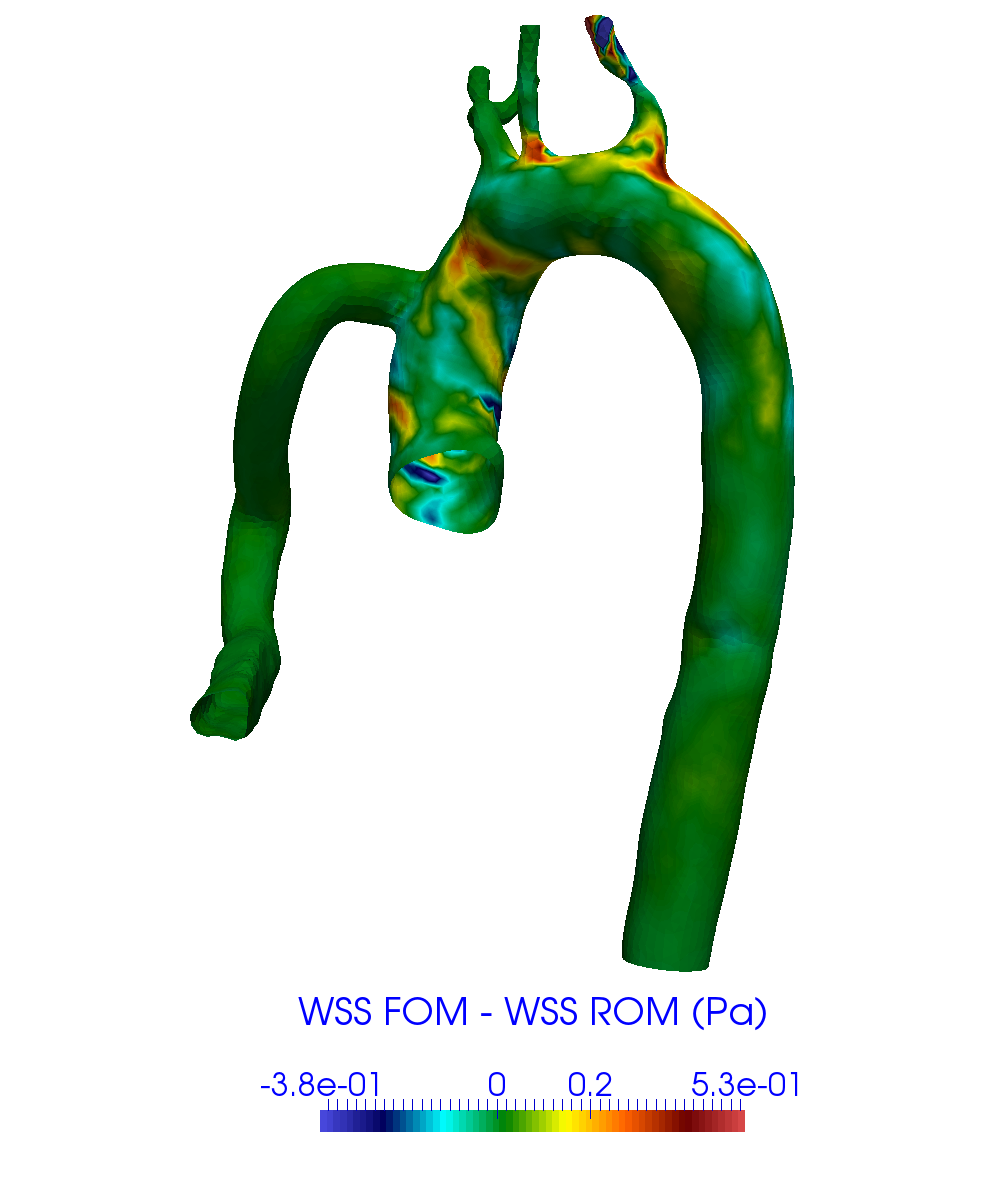}
      \end{overpic}
\caption{Comparison of the FOM/ROM pressure (1st row) and WSS (2nd row) at $PF = 3.45$ ($\omega = 5200$).}\label{fig:ROM_3_45}
\end{figure}

\begin{figure}[h]
\centering
 \begin{overpic}[width=0.3\textwidth]{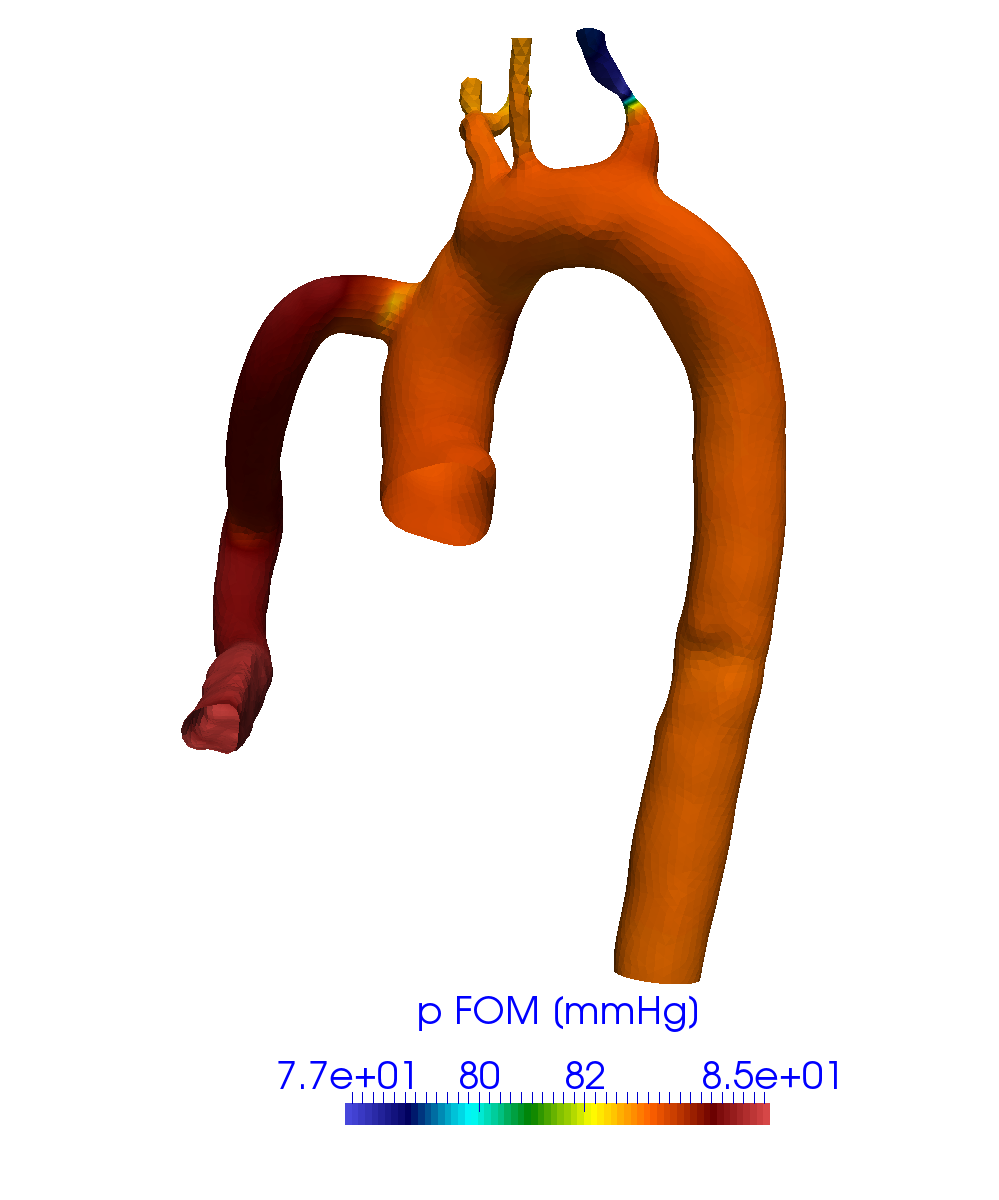}
      \end{overpic}
 \begin{overpic}[width=0.3\textwidth]{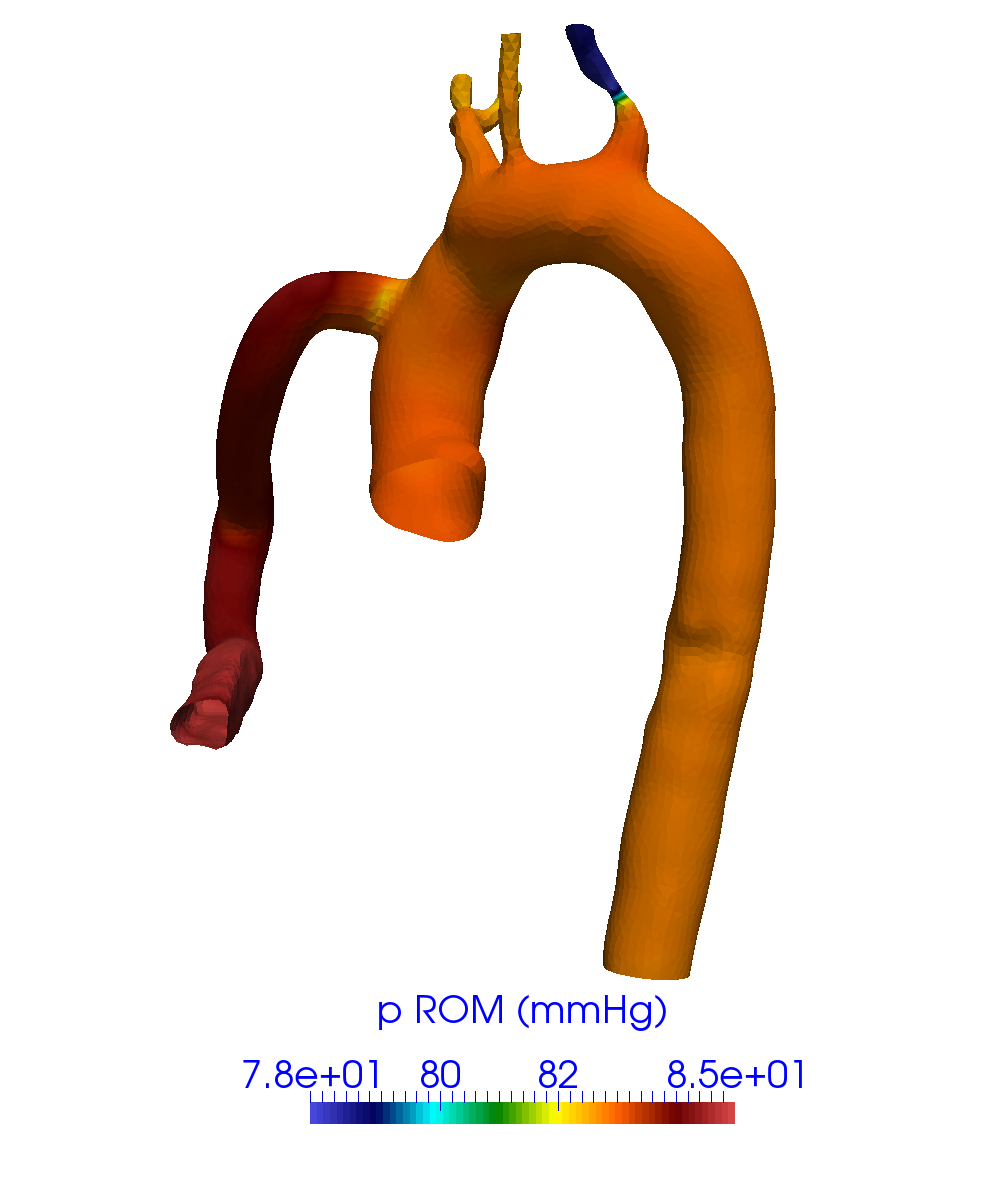}
     \end{overpic}
 \begin{overpic}[width=0.3\textwidth]{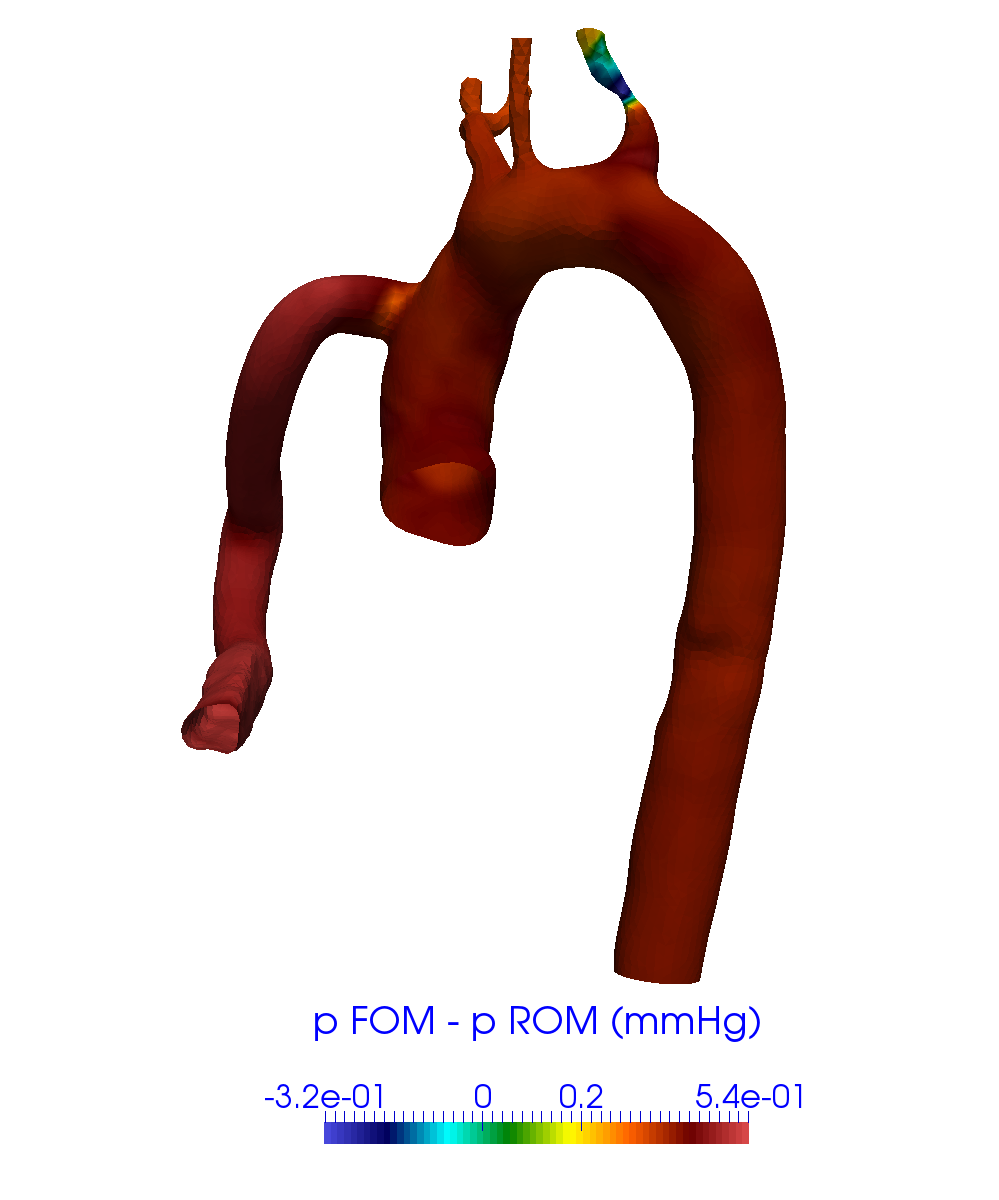}
      \end{overpic}
 \begin{overpic}[width=0.3\textwidth]{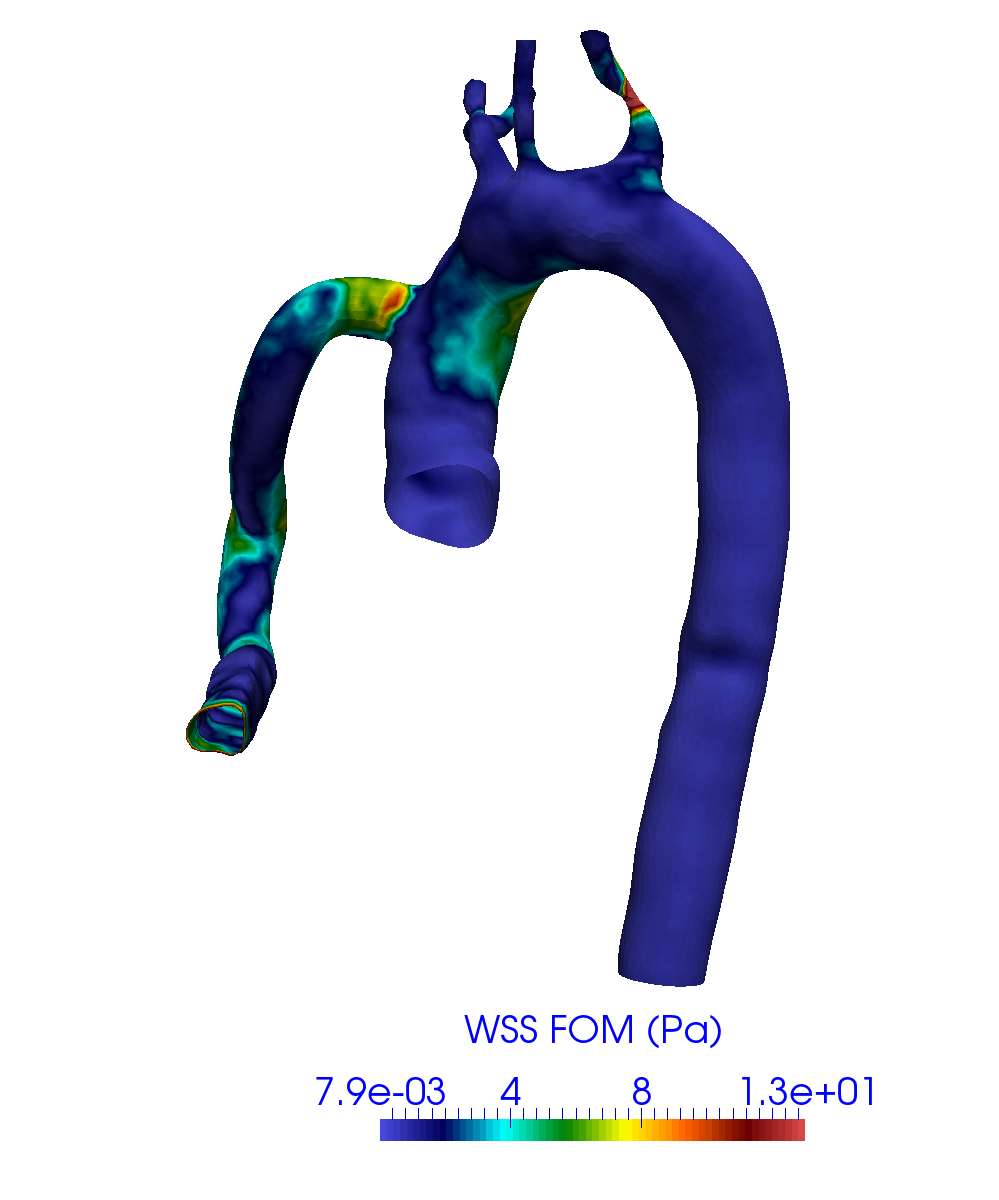}
      \end{overpic}
 \begin{overpic}[width=0.3\textwidth]{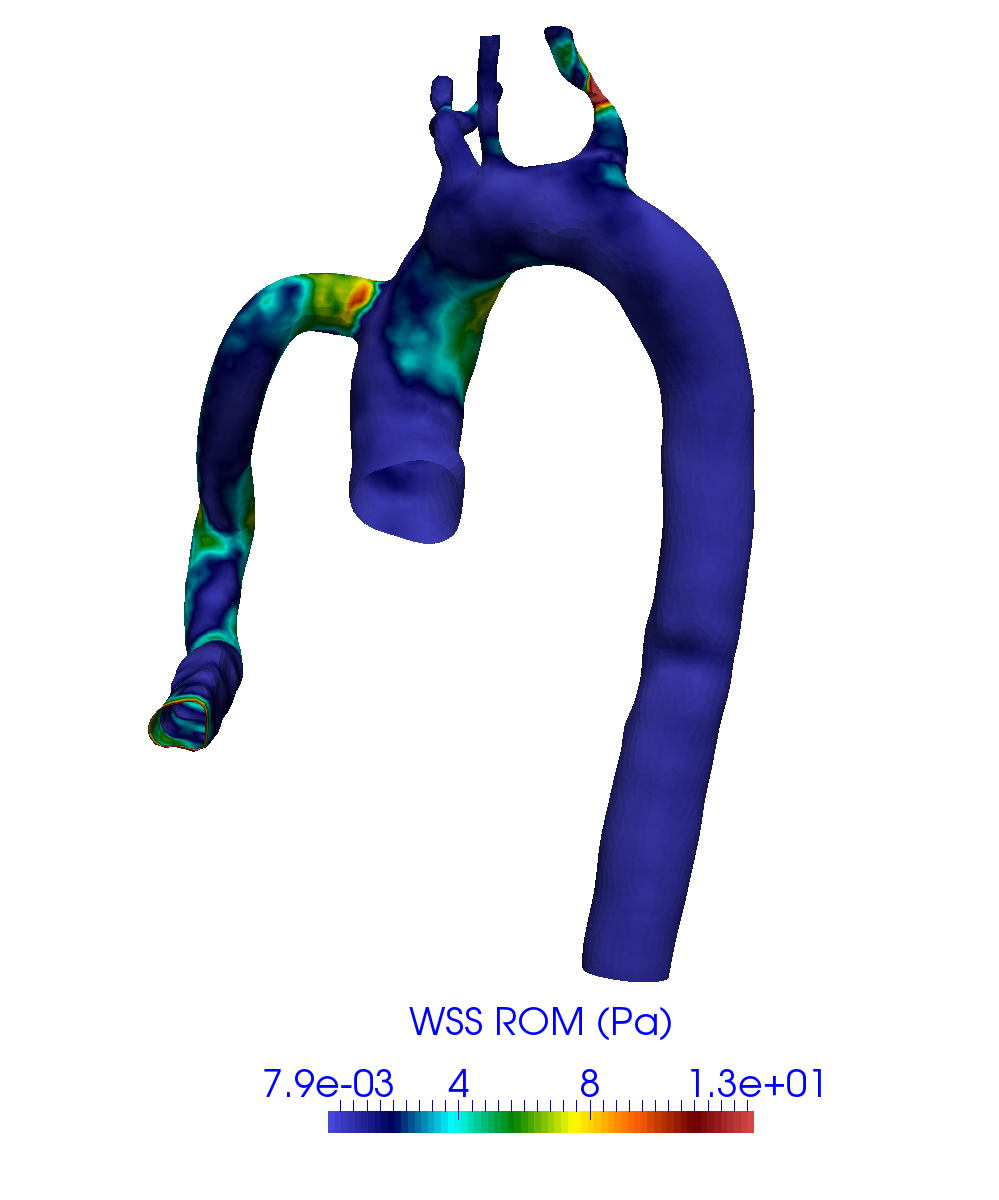}
     \end{overpic}
 \begin{overpic}[width=0.3\textwidth]{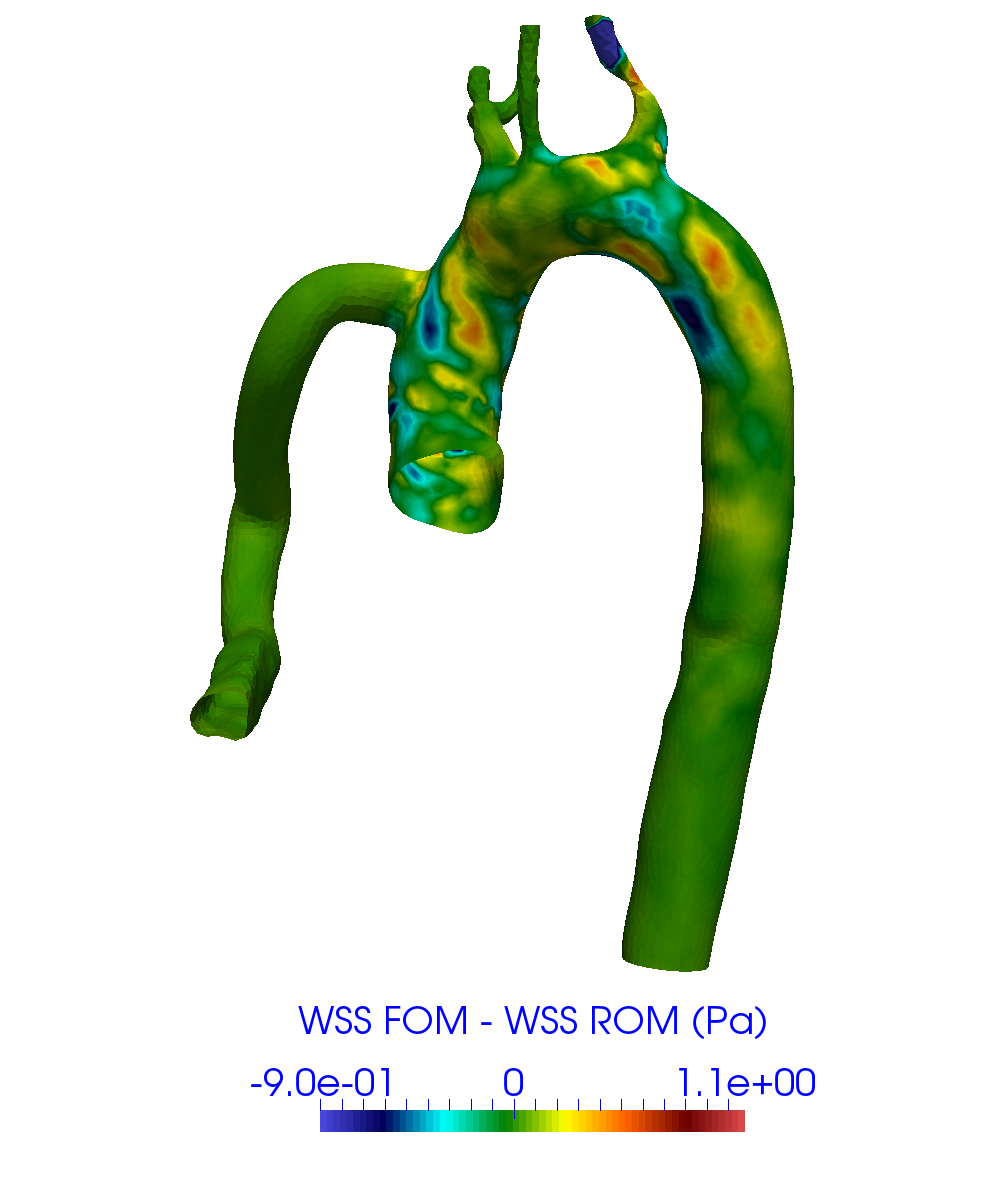}
      \end{overpic}
\caption{Comparison of the FOM/ROM pressure (1st row) and WSS (2nd row) at $PF = 4.35$ ($\omega = 5484$).}\label{fig:ROM_4_35}
\end{figure}

Figure \ref{fig:ROM_3_45} and \ref{fig:ROM_4_35} display a comparison between FOM and ROM for $p$ and WSS fields, for both $PF$ ($\omega$) values under consideration. The comparison indicates that the ROM is able to provide a good reconstruction for both variables.

\begin{figure}[h]
\centering
 \begin{overpic}[width=0.35\textwidth]{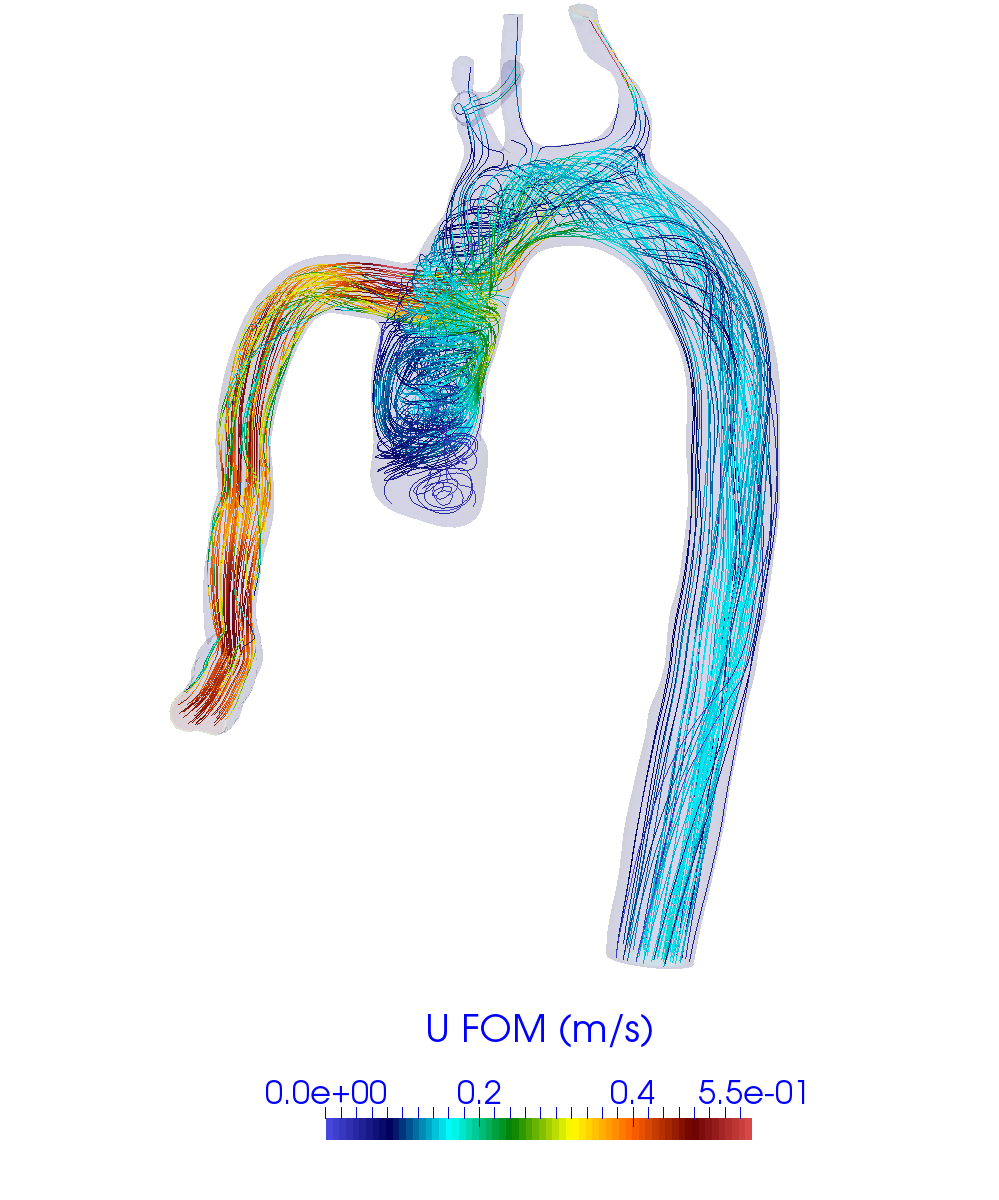}
      \end{overpic}
 \begin{overpic}[width=0.35\textwidth]{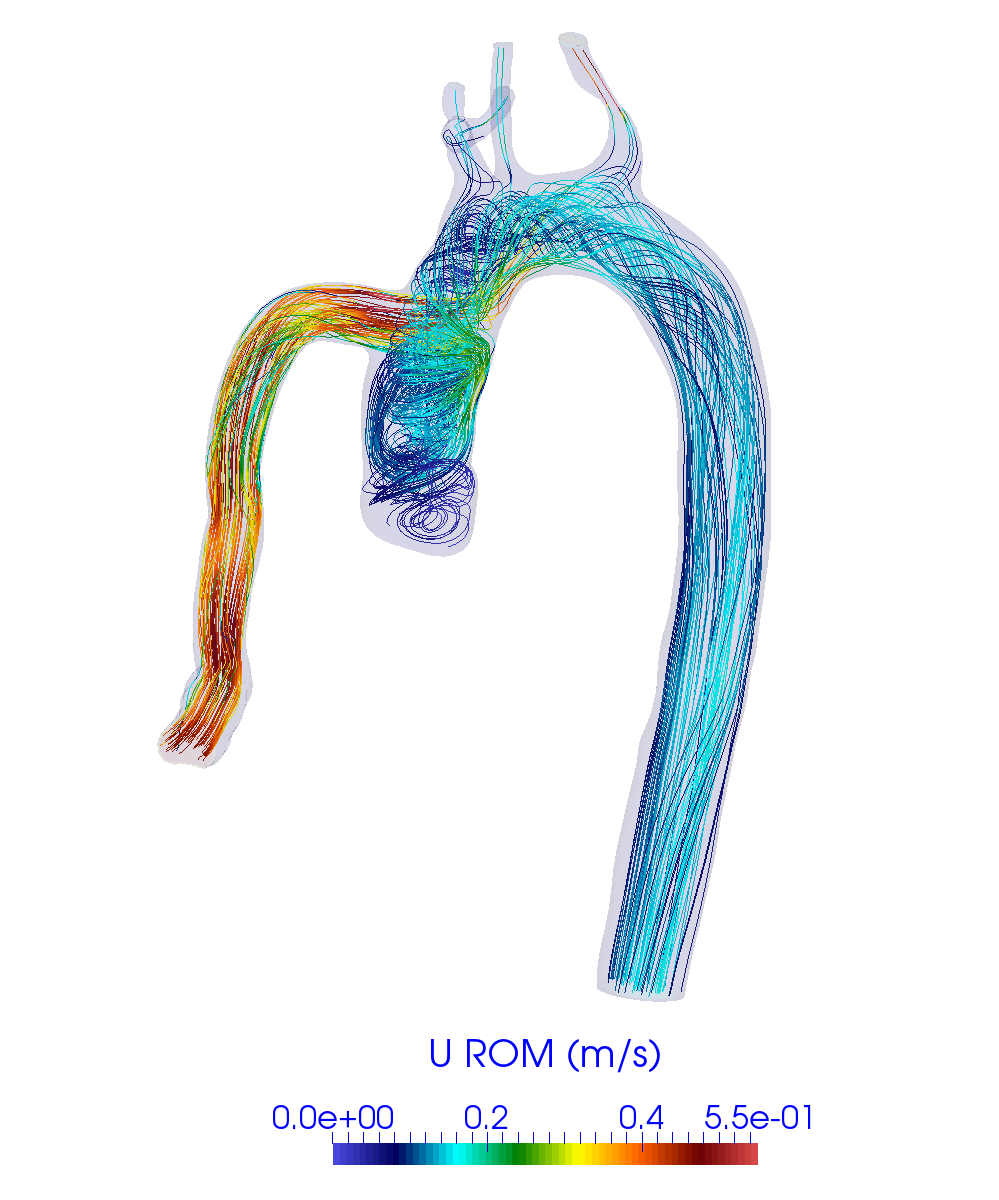}
     \end{overpic}\\
 \begin{overpic}[width=0.35\textwidth]{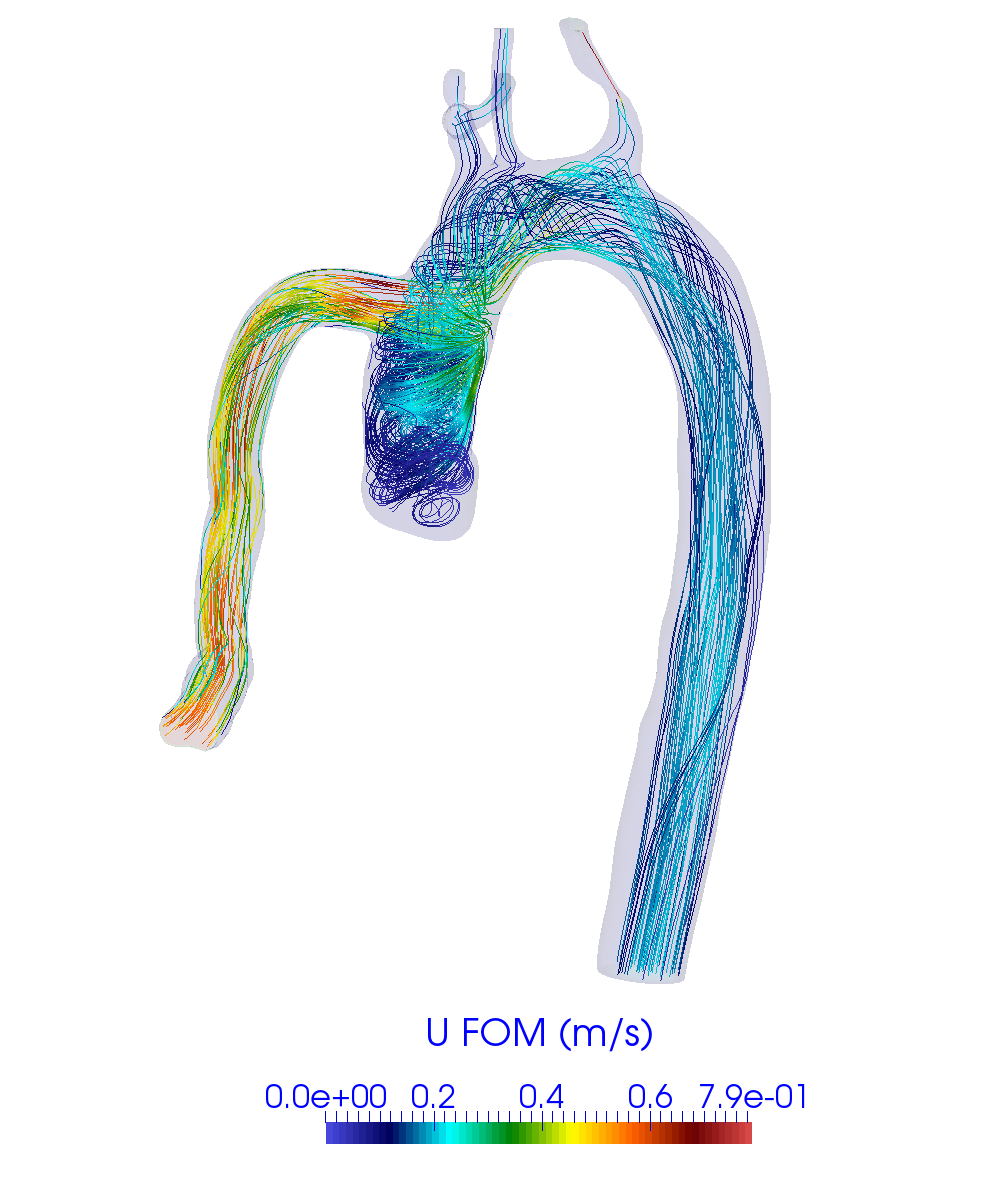}
      \end{overpic}
 \begin{overpic}[width=0.35\textwidth]{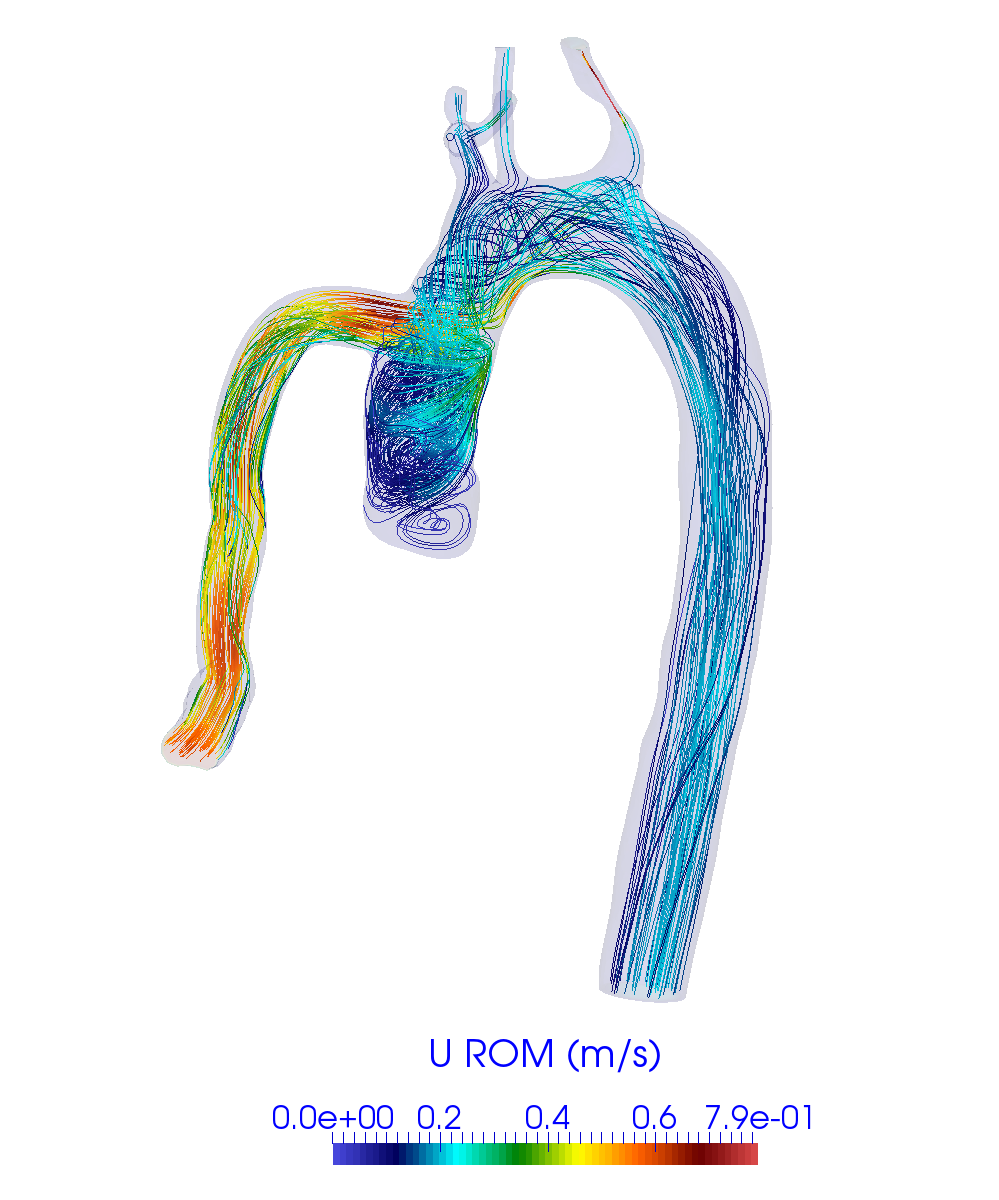}
     \end{overpic}
\caption{Comparison of the FOM/ROM velocity streamlines at $PF = 3.45$ ($\omega = 5200$) (1st row) and $PF = 4.35$ ($\omega = 5484$) (2nd row).}\label{fig:U_1}
\end{figure}

\begin{figure}[h]
\centering
 \begin{overpic}[width=0.3\textwidth]{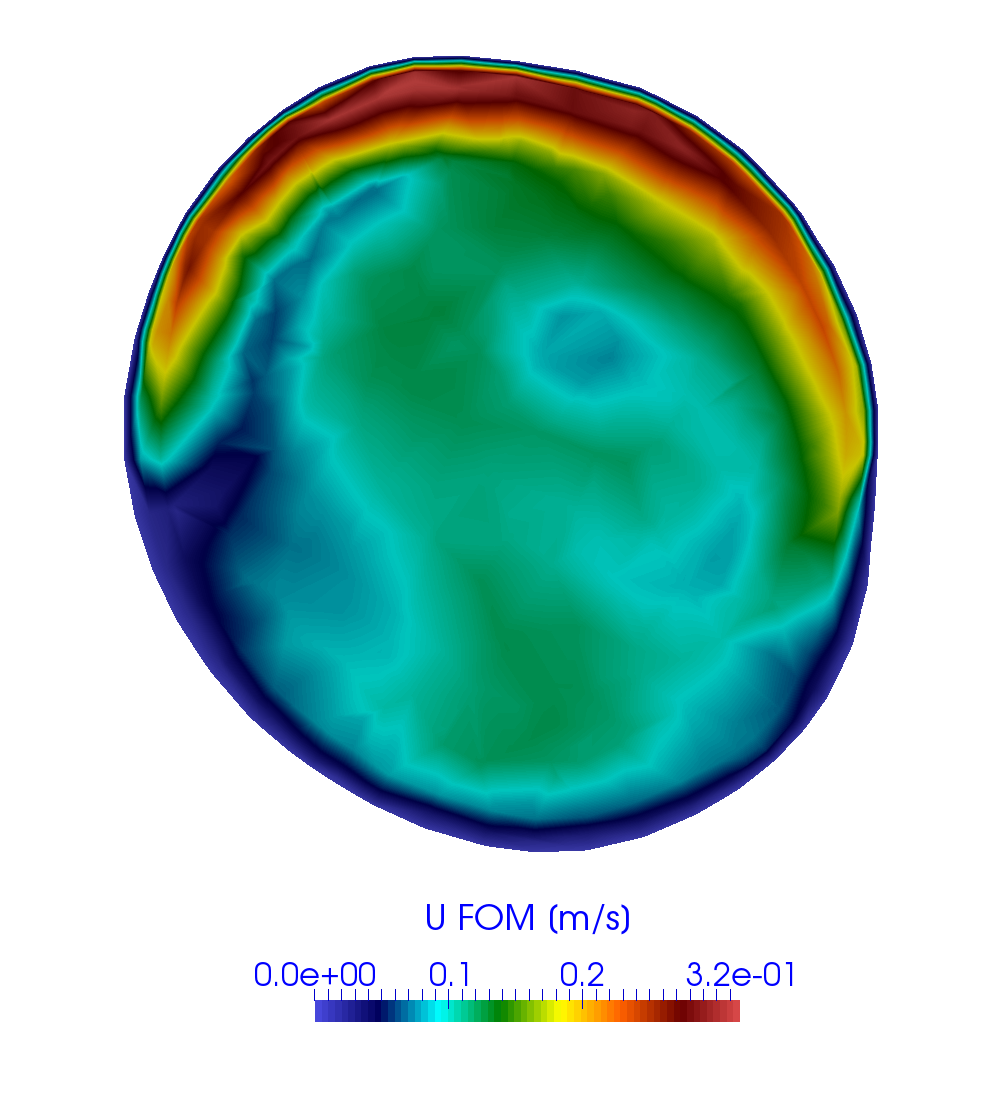}
      \end{overpic}
 \begin{overpic}[width=0.3\textwidth]{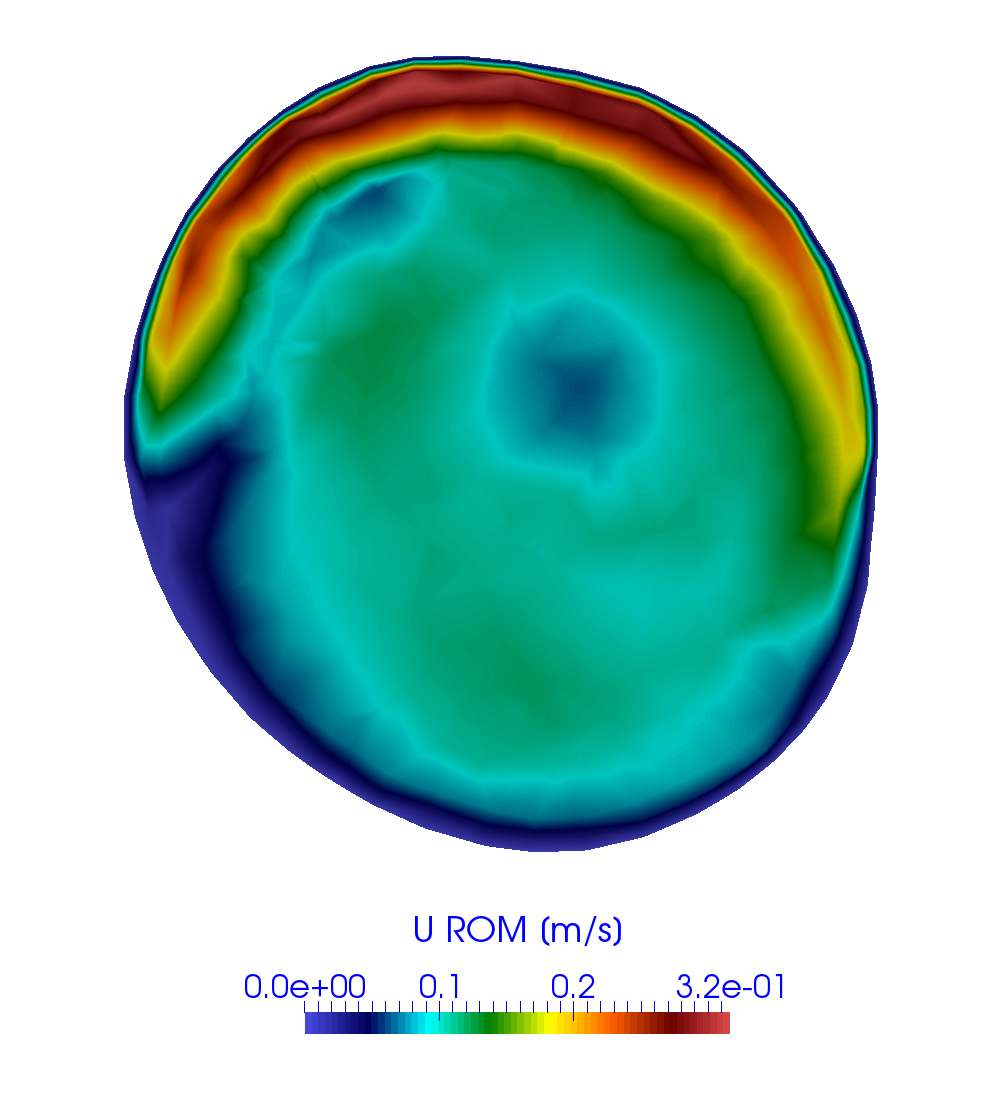}
     \end{overpic}
 \begin{overpic}[width=0.3\textwidth]{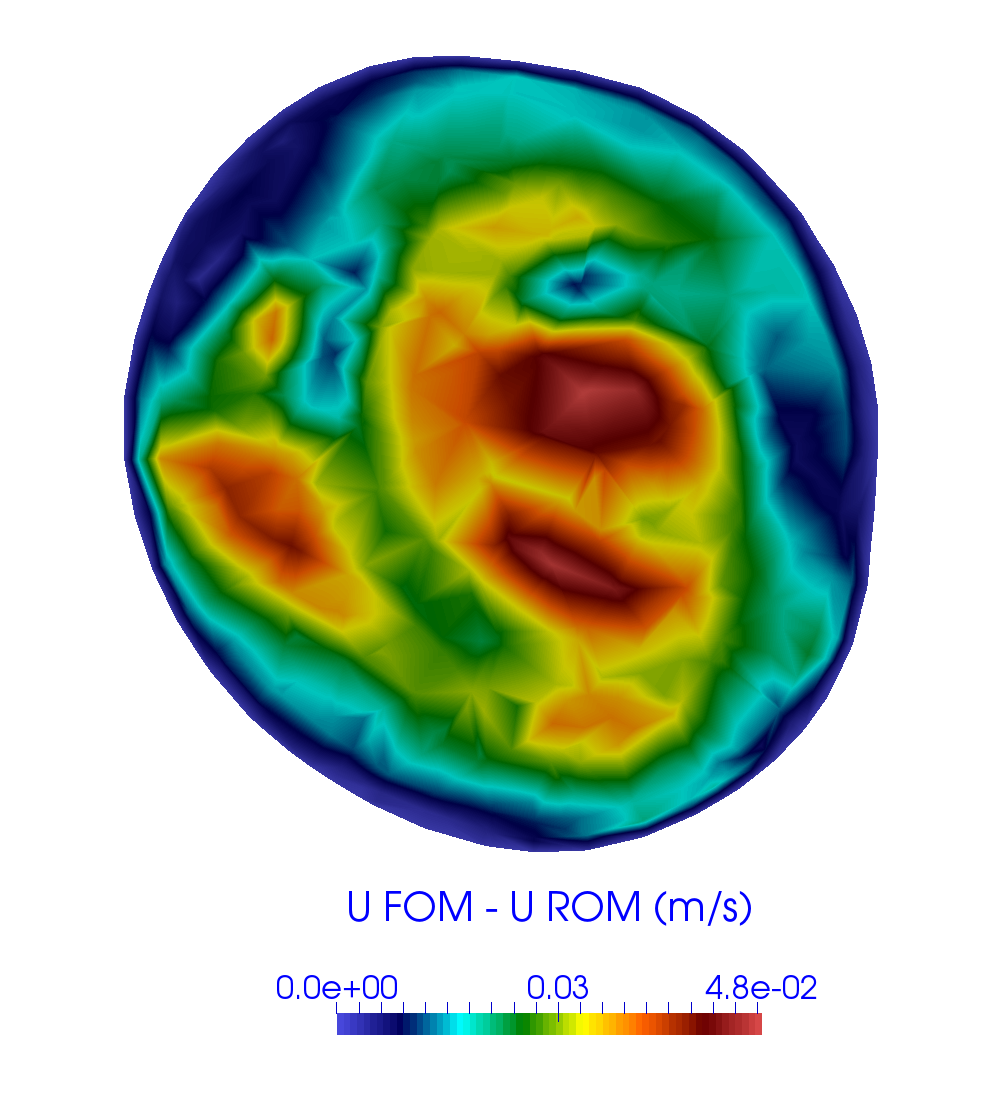}
      \end{overpic}
 \begin{overpic}[width=0.31\textwidth]{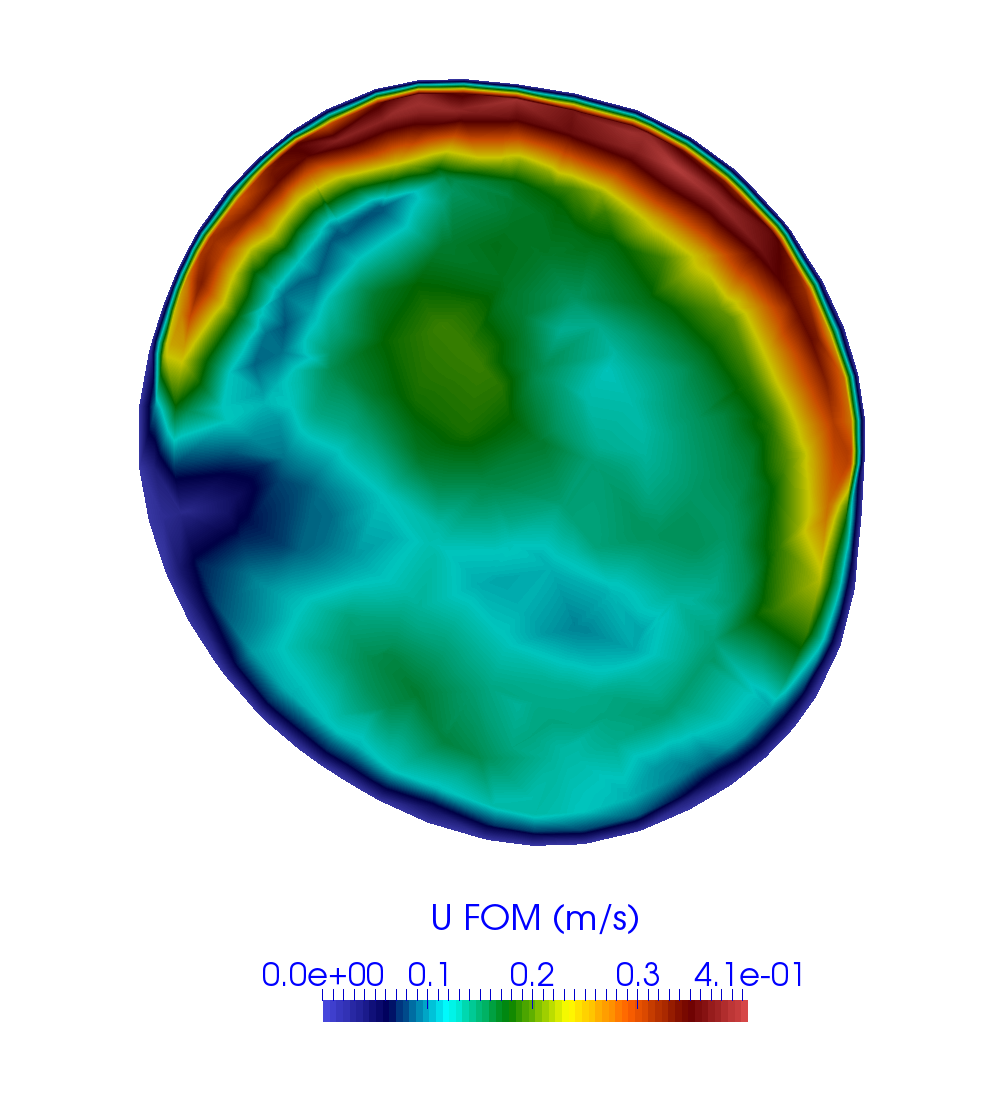}
      \end{overpic}
 \begin{overpic}[width=0.275\textwidth]{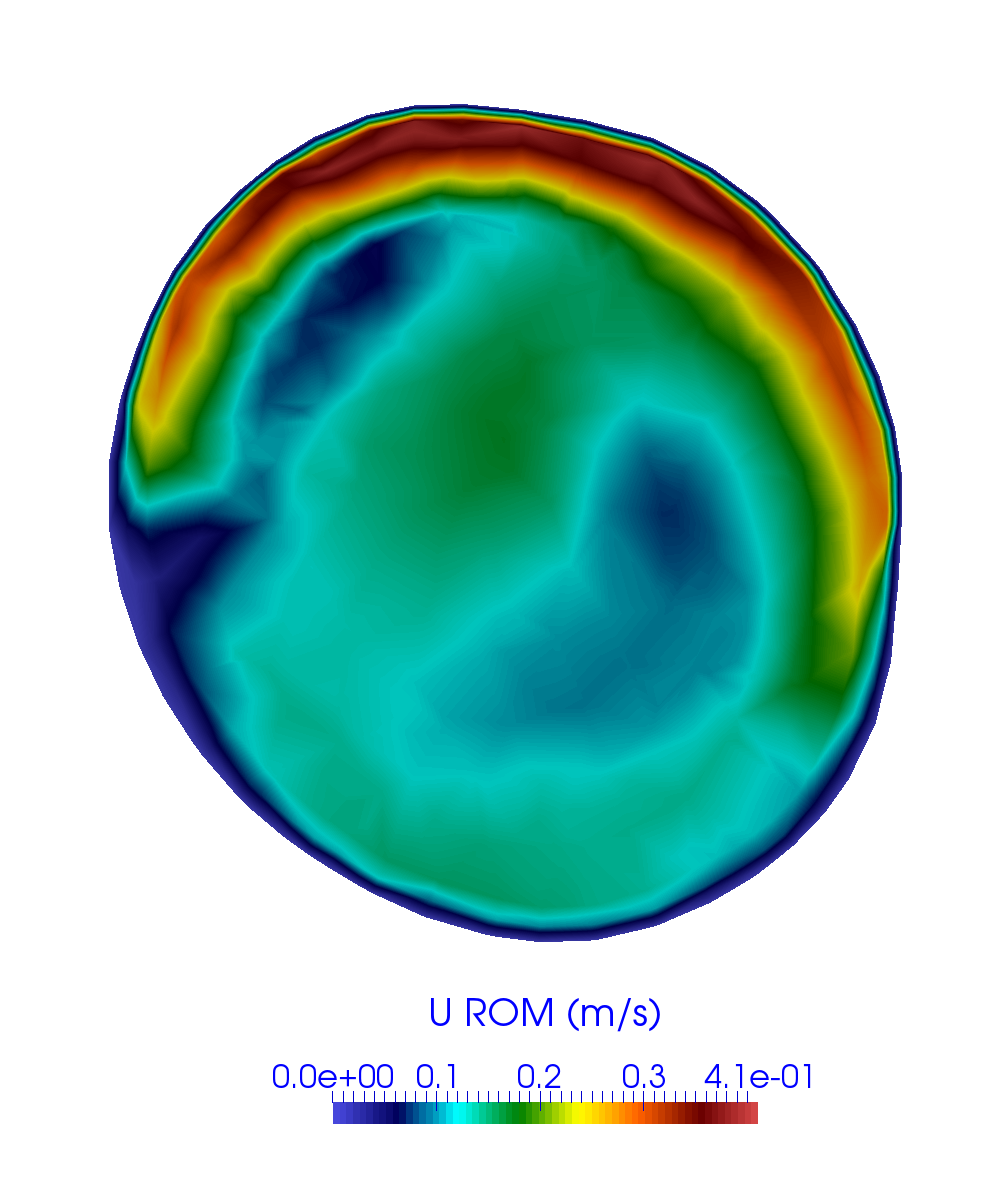}
     \end{overpic}
 \begin{overpic}[width=0.325\textwidth]{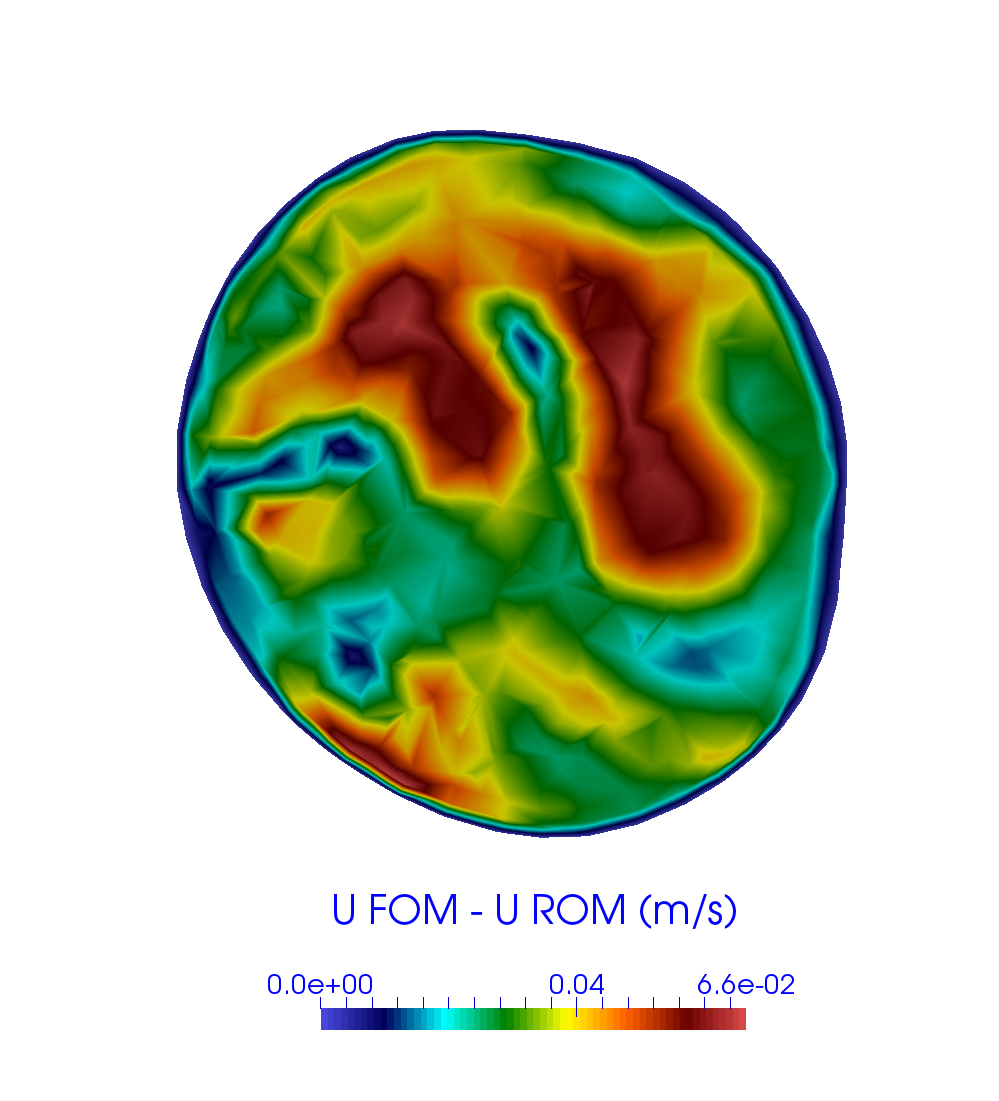}
      \end{overpic}
\caption{Comparison of the FOM/ROM velocity field related to a section of the ascending aorta next to the anastomosis location at $PF = 3.45$ ($\omega = 5200$) (1st row) and $PF = 4.35$ ($\omega = 5484$) (2nd row).}\label{fig:U_2}
\end{figure}

Figure \ref{fig:U_1} displays the velocity streamlines obtained both with FOM and ROM, and for both $PF$ ($\omega$) values under consideration. In order to further investigate the flow field, in Figure \ref{fig:U_2} a comparison between FOM and ROM for the velocity field related to a section of the ascending aorta next to the anastomosis location is showed. As observed for $p$ and WSS fields, the ROM also performs well for the velocity.

The CPU time of the FOM model is 9600s and the one of the ROM is 40s. This corresponds to a speed-up of $\approx 240$, that demonstrates the fact that it is possible to use the ROM in the place of the FOM in order to obtain accurate simulations with a significant reduction of the computational cost.

\section{Conclusion and perspectives}\label{sec:conclusion}
In this work, a parametrized non-intrusive ROM using PODI method is used for the investigation of patient-specific aortic blood flow in presence of a LVAD device. The FOM is represented by the incompressible Navier-Stokes equations discretized by using a FV technique, coupled with three-element Windkessel models to enforce outlet boundary conditions. CT images of a patient are considered for the reconstruction of the geometry as well as RCH and ECHO data are exploit for the individualization of the three-element Windkessel models coefficients used to enforce boundary conditions. Therefore, a complete patient-specific framework is presented.

In order to showcase the features of our approach, we have successfully validated the FOM both for pre-surgery and post-surgery configuration by comparing numerical and experimental data. Then, the ROM developed is used to carry out a parametric study with respect to the LVAD flow rate. We show that the ROM provides accurate solutions with a significant reduction of the computational cost, up to at least two orders of magnitudes.

As a follow-up of the present work, we are going to investigate the influence of the LVAD device on the left and right ventricle flow patterns as well as their interaction. We are also interested in efficiently handling geometrical parametrization (e.g. in order to consider different anastomosis angles, or different designs of the outflow cannula) in the context of patient-specific geometries, extending e.g. the work carried out in \cite{BallarinFaggianoManzoniQuarteroniRozzaIppolitoAntonaScrofani2016} to different applications and different model reduction techniques.

Furthermore, we are designing an application that will be available online and allow computation to be run from standard web browsers. 

\section{Acknowledgements}\label{sec:acknowledgements}
We acknowledge the support provided by the European Research Council
Executive Agency by the Consolidator Grant project AROMA-CFD "Advanced
Reduced Order Methods with Applications in Computational Fluid Dynamics" -
GA 681447, H2020-ERC CoG 2015 AROMA-CFD and INdAM-GNCS 2020 project "Tecniche Numeriche Avanzate per Applicazioni Industriali".

\bibliographystyle{amsplain_mod}
\bibliography{bib/bibfile}
\end{document}